\documentclass[11pt,a4paper,final, dvipsnames]{amsart}

\input xy
\xyoption{all}

\DeclareMathAlphabet{\mathpzc}{OT1}{pzc}{m}{it}

\usepackage{amsfonts,amsmath,amssymb,indentfirst,mathrsfs,amscd}
\usepackage[mathscr]{eucal}
\usepackage{graphicx}
\usepackage{nicefrac}
\usepackage{amsthm}

\usepackage{color}
\usepackage[utf8]{inputenc}
\usepackage[spanish, english]{babel}
\usepackage{booktabs}
\usepackage{multirow}
\usepackage{verbatim}

\usepackage{calrsfs}
\DeclareMathAlphabet{ \mathcal}{OMS}{zplm}{m}{n}

\newcommand{\noop}[1]{}

\usepackage[all]{xy}
\usepackage{graphicx}
\usepackage{stmaryrd}
\usepackage{enumitem}

\usepackage{empheq}
\usepackage{color}
\usepackage{xcolor}

\author[]{\'Angel Gonz\'alez-Prieto}
\address{Instituto de Ciencias Matem\'aticas CSIC-UAM-UC3M-UCM. C.\ Nicol\'as Cabrera, 13-15, 28049. Madrid, Spain.}
\address{Escuela T\'ecnica Superior de Ingenieros de Sistemas Inform\'aticos, Universidad Polit\'ecnica de Madrid, Calle Alan Turing s/n (Carretera de Valencia Km 7), 28031 Madrid, Spain}\email{angel.gonzalez.prieto@icmat.es}
\title[]{Motivic theory of representation varieties via\\Topological Quantum Field Theories}

\keywords{}

\DeclareMathOperator{\img}{Im \,}           %Image%
\DeclareMathOperator{\coker}{coker\,}
\DeclareMathOperator{\Hom}{Hom\,}           %Hom%
\DeclareMathOperator{\tr}{Tr\,}             %Trace Tr%
\DeclareMathOperator{\Aut}{Aut\,}           %Aut%

\DeclareMathOperator{\Gr}{Gr}

\usepackage{amsmath}
\begin{document}

\newtheorem{thm}{Theorem}[section]
\newtheorem{prop}[thm]{Proposition}
\newtheorem{lem}[thm]{Lemma}
\newtheorem{cor}[thm]{Corollary}
\newtheorem{conjecture}{Conjecture}
\newtheorem*{theorem*}{Theorem}

\theoremstyle{definition}
\newtheorem{defn}[thm]{Definition}
\newtheorem{ex}[thm]{Example}
\newtheorem{as}{Assumption}

\theoremstyle{remark}
\newtheorem{rmk}[thm]{Remark}

\theoremstyle{remark}
\newtheorem*{prf}{Proof}

\newcommand{\iacute}{\'{\i}} %i con acento%
\newcommand{\norm}[1]{\lVert#1\rVert} %norma%

\newcommand{\lto}{\longrightarrow}
\newcommand{\hra}{\hookrightarrow}

\newcommand{\suchthat}{\;\;|\;\;}
\newcommand{\dbar}{\overline{\partial}}

\newcommand{\cA}{\mathcal{A}}
\newcommand{\cC}{\mathcal{C}}
\newcommand{\cD}{\mathcal{D}}
\newcommand{\cE}{\mathcal{E}}
\newcommand{\cF}{\mathcal{F}}
\newcommand{\cG}{\mathcal{G}} %Gauge group%
\newcommand{\cI}{\mathcal{I}} %Gauge group%
\newcommand{\cO}{\mathcal{O}} %Holomorphic functions sheaf%
\newcommand{\cM}{\mathcal{M}} %Moduli space% %moduli of parabolic bundles% %moduli of U(p,q) bundles%
\newcommand{\cN}{\mathcal{N}} %Space of minimal points of the Morse function%
\newcommand{\cP}{\mathcal{P}} %Moduli of K(D) pairs%
\newcommand{\cQ}{\mathcal{Q}} %Moduli of K(D) pairs%
\newcommand{\cS}{\mathcal{S}} %Moduli of solutions of Hitchin's equations, contructed by Konno%
\newcommand{\cU}{\mathcal{U}} %Moduli of stable U(p,q) parabolic Higgs bundles%
\newcommand{\cJ}{\mathcal{J}}
\newcommand{\cX}{\mathcal{X}}
\newcommand{\cT}{\mathcal{T}}
\newcommand{\cV}{\mathcal{V}}
\newcommand{\cW}{\mathcal{W}}
\newcommand{\cB}{\mathcal{B}}
\newcommand{\cR}{\mathcal{R}}
\newcommand{\cH}{\mathcal{H}}
\newcommand{\cZ}{\mathcal{Z}}
\newcommand{\D}{\bar{B}}

\newcommand{\Ker}{\textrm{Ker}\,}
\renewcommand{\coker}{\textrm{coker}\,}

\newcommand{\ext}{\mathrm{ext}} % an extension%
\newcommand{\x}{\times}

\newcommand{\mM}{\mathscr{M}} %Meromorphic function sheaf%

\newcommand{\CC}{\mathbb{C}} %Complex numbers%
\newcommand{\QQ}{\mathbb{Q}} %Rational numbers%
\newcommand{\FF}{\mathbb{F}} %Rational numbers%
\newcommand{\PP}{\mathbb{P}} %projective space%
\newcommand{\HH}{\mathbb{H}} %Hypercohomology, quaternions..%
\newcommand{\RR}{\mathbb{R}} %Real numbers%
\newcommand{\ZZ}{\mathbb{Z}} %Integer numbers%
\newcommand{\NN}{\mathbb{N}} %Natural numbers%
\newcommand{\DD}{\mathbb{D}} %Natural numbers%

\renewcommand{\lg}{\mathfrak{g}} %Lie algebra of G%
\newcommand{\lh}{\mathfrak{h}} %Lie algebra of H%
\newcommand{\lu}{\mathfrak{u}} %Lie algebra of U%
\newcommand{\la}{\mathfrak{a}} %Lie algebra of A%
\newcommand{\lb}{\mathfrak{b}} %Lie algebra of B%
\newcommand{\lm}{\mathfrak{m}} %Lie algebra of M%
\newcommand{\lgl}{\mathfrak{gl}} %Lie algebra of GL%
\newcommand{\lZ}{\mathfrak{Z}} %Almost-TQFT%

\newcommand{\too}{\longrightarrow}
\newcommand{\imat}{\sqrt{-1}} %i%
\newcommand{\tinyclk}{{\scriptscriptstyle \Taschenuhr}} %Tiny Clock Symbol from ifsym package%
\newcommand\restr[2]{\left.#1\right|_{#2}}
\newcommand\rtorus[1]{{\mathbb{T}^{#1}}}
\newcommand\actPartial{\overline{\partial}}
\newcommand\handle[2]{\mathcal{A}^{#1}_{#2}}

\newcommand\pastingArea[2]{S^{#1-1} \times \bar{B}^{#2-#1}}
\newcommand\pastingAreaPlus[2]{S^{#1} \times \bar{B}^{#2-#1-1}}
\newcommand\coordvector[2]{\left.\frac{\partial}{\partial {#1}}\right|_{#2}}

% Categories
\newcommand\Sets{\textbf{Set}}
\newcommand\Cat{\textbf{Cat}}
\newcommand\Top{\textbf{Top}}
\newcommand\TopHLC{\textbf{Top}_{hlc}}
\newcommand\TopS{\textbf{Top}_\star}
\newcommand\Diff{\textbf{Diff}}
\newcommand\Diffc{\textbf{Diff}_c}

\newcommand\CBord[3]{\mathbf{Bd}_{{#1 #3}}^{#2}}
\newcommand\Bord[1]{\CBord{#1}{}{}}
\newcommand\EBord[2]{\CBord{#1}{#2}{}}
\newcommand\Bordo[1]{\CBord{#1}{or}{}}
\newcommand\Bordp[1]{\mathbf{Bdp}_{{#1}}}
\newcommand\Bordpar[2]{\mathbf{Bd}_{{#1}}(#2)}
\newcommand\Bordppar[2]{\mathbf{Bdp}_{{#1}}(#2)}

\newcommand\Tubo[1]{\mathbf{Tb}_{#1}^0}
\newcommand\Tub[1]{\mathbf{Tb}_{#1}}
\newcommand\ETub[2]{\mathbf{Tb}_{#1}^{#2}}
\newcommand\Tubp[1]{\mathbf{Tbp}_{#1}}
\newcommand\Tubpo[1]{\mathbf{Tbp}_{#1}^0}
\newcommand\Tubppar[2]{\mathbf{Tbp}_{#1}(#2)}

\newcommand\Embc{\textbf{Emb}_c}
\newcommand\EEmbc[1]{\textbf{Emb}_c^{#1}}
\newcommand\Embpc{\textbf{Embp}_c}
\newcommand\Embparc[1]{\textbf{Emb}_c(#1)}
\newcommand\Embpparc[1]{\textbf{Embp}_c(#1)}

\newcommand\cSp{\cS_p}
\newcommand\cSpar[1]{\cS_{#1}}

\newcommand\Diffpc{\textbf{Diff}_c}

\newcommand\PVar[1]{\mathbf{PVar}_{#1}}
\newcommand\PVarC{\mathbf{PVar}_{\CC}}
\newcommand\Sr[1]{\mathrm{S}{#1}}

\newcommand\Bordpo[1]{\mathbf{Bdp}_{{#1}}^{or}}
\newcommand\ClBordp[1]{\mathbf{l}\CBord{#1}{}{}}
\newcommand\CTub[3]{\mathbf{Tb}_{{#1 #3}}^{#2}}
\newcommand\CTubp[1]{\CTub{#1}{}{}}
\newcommand\CTubpp[1]{\mathbf{Tbp}_{#1}{}{}}
\newcommand\CTubppo[1]{\mathbf{Tbp}_{#1}^0}
\newcommand\CClose[3]{\mathbf{Cl}_{{#1 #3}}^{#2}}
\newcommand\CClosep[1]{\CClose{#1}{}{}}
\newcommand\Obj[1]{\mathrm{Obj}(#1)}
\newcommand\Mor[1]{\mathrm{Mor}(#1)}
\newcommand\Vect[1]{{#1}\textrm{-}\mathbf{Vect}}
\newcommand\Vecto[1]{{#1}\textrm{-}\mathbf{Vect}_0}
\newcommand\Mod[1]{{#1}\textrm{-}\mathbf{Mod}}
\newcommand\Modt[1]{{#1}\textrm{-}\mathbf{Mod}_t}
\newcommand\Rng{\mathbf{Ring}}
\newcommand\Grp{\mathbf{Grp}}
\newcommand\Grpd{\mathbf{Grpd}}
\newcommand\Grpdo{\mathbf{Grpd}_0}
\newcommand\HS[1]{\mathbf{HS}^{#1}}
\newcommand\MHS[1]{\mathbf{MHS}}
\newcommand\PHS[2]{\mathbf{HS}^{#1}_{#2}}
\newcommand\MHSq{\mathbf{HS}}
\newcommand\Sch{\mathbf{Sch}}
\newcommand\Sh[1]{\mathbf{Sh}\left(#1\right)}
\newcommand\QSh[1]{\mathbf{QSh}\left(#1\right)}
\newcommand\Var[1]{\mathbf{Var}_{#1}}
\newcommand\Varc{\mathbf{Var}}
\newcommand\CVar{\Varc}
\newcommand\PHM[2]{\cM_{#1}^p(#2)}
\newcommand\MHM[1]{\cM_{#1}}
\newcommand\HM[2]{\textrm{HM}^{#1}(#2)}
\newcommand\HMW[1]{\textrm{HMW}(#1)}
\newcommand\VMHS[1]{VMHS({#1})}
\newcommand\geoVMHS[1]{VMHS_g({#1})}
\newcommand\goodVMHS[1]{\mathrm{VMHS}_0({#1})}
\newcommand\Par[1]{\mathrm{Par}({#1})}
\newcommand\K[1]{\mathrm{K}#1}
\newcommand\Ko[1]{\left(\mathrm{K}{#1}\right)_0}
\newcommand\KM[1]{\mathrm{K}\MHM{#1}}
\newcommand\KMo[1]{\mathrm{K}{\MHM{#1}}_0}
\newcommand\Ab{\mathbf{Ab}}
\newcommand\CPP{\cP\cP}
\newcommand\Bim[1]{{#1}\textrm{-}\mathbf{Bim}}
\newcommand\Span[1]{\mathrm{Span}({#1})}
\newcommand\Spano[1]{\mathrm{Span}^{op}({#1})}
\newcommand\GL[1]{\mathrm{GL}_{#1}}
\newcommand\SL[1]{\mathrm{SL}_{#1}}
\newcommand\PGL[1]{\mathrm{PGL}_{#1}}
\newcommand\Rep[1]{\mathfrak{X}_{#1}}
\newcommand\KVarrel[1]{\K{\mathbf{Var}_{\downarrow#1}}}
\newcommand\Varrel[1]{\mathbf{Var}_{\downarrow#1}}
\newcommand\KVar[1]{\K{\mathbf{Var}_{#1}}}
\newcommand\KVarc{\K{\mathbf{Var}}}

\newcommand\Qtm[1]{\mathcal{Q}_{#1}}
\newcommand\sQtm[1]{\mathcal{Q}_{#1}^0}
\newcommand\Fld[1]{\mathcal{F}_{#1}}
\newcommand\Id{\mathrm{Id}}

% Hodge theory
\newcommand\DelHod[1]{e\left(#1\right)}
\newcommand\RDelHod{e}
\newcommand\eVect{\mathcal{E}}
\newcommand\e[1]{\eVect\left(#1\right)}
\newcommand\intMor[2]{\int_{#1}\,#2}

% Representation Varieties
\newcommand\Dom[1]{\mathcal{D}_{#1}}

\newcommand\Xf[1]{{X}_{#1}}					% Free non-parabolic
\newcommand\Xs[1]{\mathfrak{X}_{#1}}							% Surface group non-parabolic

\newcommand\Xft[2]{\overline{{X}}_{#1, #2}} % Free tr=2
\newcommand\Xst[2]{\overline{\mathfrak{X}}_{#1, #2}}			% Surface group tr=2

\newcommand\Xfp[2]{X_{#1, #2}}			% Free J+
\newcommand\Xsp[2]{\mathfrak{X}_{#1, #2}}						% Surface group J+

\newcommand\Xfd[2]{X_{#1; #2}}			% Free diagonal
\newcommand\Xsd[2]{\mathfrak{X}_{#1; #2}}						% Surface group diagonal

\newcommand\Xfm[3]{\mathcal{X}_{#1, #2; #3}}		% Free mixed
\newcommand\Xsm[3]{X_{#1, #2; #3}}					% Surface group mixed

\newcommand\XD[1]{#1^{D}}
\newcommand\XDh[1]{#1^{\delta}}
\newcommand\XU[1]{#1^{UT}}
\newcommand\XP[1]{#1^{U}}
\newcommand\XPh[1]{#1^{\upsilon}}
\newcommand\XI[1]{#1^{\iota}}
\newcommand\XTilde[1]{#1^{\varrho}}
\newcommand\Xred[1]{#1^{r}}
\newcommand\Xirred[1]{#1^{ir}}

\newcommand\Char[1]{\cR_{#1}}
\newcommand\Chars[1]{\cR_{#1}}
\newcommand\CharW[1]{\mathscr{R}_{#1}}

% Derived Category
\newcommand\Ch[1]{\textrm{Ch}\,{#1}}
\newcommand\Chp[1]{\textrm{Ch}^+{#1}}
\newcommand\Chm[1]{\textrm{Ch}^-{#1}}
\newcommand\Chb[1]{\textrm{Ch}^b{#1}}
\newcommand\Der[1]{\textrm{D}{#1}}
\newcommand\Dp[1]{\textrm{D}^+{#1}}
\newcommand\Dm[1]{\textrm{D}^-{#1}}
\newcommand\Db[1]{\textrm{D}^b{#1}}

% TQFT
\newcommand\Gs{\cG}
\newcommand\Gq{\cG_q}
\newcommand\Gg{\cG_c}
\newcommand\Zs[1]{Z_{#1}}
\newcommand\Zg[1]{Z^{gm}_{#1}}
\newcommand\cZg[1]{\cZ^{gm}_{#1}}

% Miscelany
\newcommand\RM[2]{R\left(\left.#1\right|#2\right)}
\newcommand\RMc[3]{R_{#1}\left(\left.#2\right|#3\right)}
\newcommand\set[1]{\left\{#1\right\}}
\newcommand{\Stab}{\textrm{Stab}\,} %\Stabilizer%
\renewcommand{\tr}{\textrm{tr}\,}             %Trace Tr%
\newcommand\EuChS{E}             %Trace Tr%
\newcommand\EuCh[1]{E\left(#1\right)}             %Trace Tr%

\newcommand{\Ann}{\textrm{Ann}\,}
\newcommand{\Rad}{\textrm{Rad}\,}  
\newcommand\supp[1]{\mathrm{supp}{(#1)}}
\newcommand\coh[1]{\left[#1\right]}
\newcommand\Bt{\Theta}
\newcommand\Be{B_e}
\newcommand\re{\textrm{Re}\,}
\newcommand\imag{\textrm{Im}\,}
\newcommand\Kahc{\textbf{K\"ah}_c}

% Mixed Hodge modules
\newcommand\Ccs[1]{C_{cs}(#1)}
\newcommand\can{\textrm{can}}
\newcommand\var{\textrm{var}}
\newcommand\Perv[1]{\textrm{Perv}(#1)}
\newcommand\DR[1]{\textrm{DR}{#1}}
\renewcommand\Gr[2]{\textrm{Gr}_{#1}^{#2}\,}
\newcommand\ChV{\textrm{Ch}\,}
\newcommand\VerD{^{\textrm{Ve}}\DD}
\newcommand\DHol[1]{\textrm{D}^b_{\textrm{hol}}(\cD_{#1})}
\newcommand\RegHol[1]{\Mod{\cD_{#1}}_{\textrm{rh}}}
\newcommand\Drh[1]{\textrm{D}^b_{\textrm{rh}}(\cD_{#1})}
\newcommand\Dcs[2]{\textrm{D}^b_{\textrm{cs}}({#1}; {#2})}
\newcommand{\rat}{\mathrm{rat}}
\newcommand{\dmod}{\mathrm{Dmod}}

\hyphenation{mul-ti-pli-ci-ty}

\hyphenation{mo-du-li}

\begin{abstract}
In this paper, we use lax monoidal TQFTs as an effective computational method for motivic classes of representation varieties. In particular, we perform the calculation for parabolic $\SL{2}(\CC)$-representation varieties over a closed orientable surface of arbitrary genus and any number of marked points with holonomies of Jordan type. This technique is based on a building method of lax monoidal TQFTs of physical inspiration that generalizes the construction of Gonz\'alez-Prieto, Logares and Mu\~noz.
\end{abstract}
\null
\vspace{-1.1cm}
\maketitle

\vspace{-0.8cm}

%%%%%%%%%%%%%%%%%%%%%%%%%%%%%%%%%%%%%%%%%%%%%%%%%%%%%%%%%%%%%%%%
%%%%%%%%%%%%%%%%%%% SECTION: INTRODUCTION %%%%%%%%%%%%%%%%%%%%%%
%%%%%%%%%%%%%%%%%%%%%%%%%%%%%%%%%%%%%%%%%%%%%%%%%%%%%%%%%%%%%%%%

\section{Introduction}
\let\thefootnote\relax\footnotetext{\noindent \emph{2010 Mathematics Subject Classification}. Primary:
 14C15. %Motivic theory
 Secondary:
 57R56, % TQFT
 14D07, % Variations of Hodge structures
 14D21. % Applications of moduli in mathematical physics

\emph{Key words and phrases}: Grothendieck ring, TQFT, representation varieties.}

Let $G$ be a complex algebraic group, let $M$ be a compact manifold and let $A \subseteq M$ be a non-empty finite set of points. The set of representations of the fundamental groupoid $\Pi(M, A)$ into $G$ has a natural structure of complex algebraic variety. This variety is called the $G$-representation variety of $(M, A)$ and it is denoted by $\Rep{G}(M, A) = \Hom(\Pi(M,A), G)$. In the case that $M$ is connected and $A$ has a single point, this variety is usually shorten $\Rep{G}(M) = \Hom(\pi_1(M), G)$.

Even more, we can consider a parabolic structure $Q$ on $M$. It is given by a finite set of codimension two submanifolds $S_1, \ldots, S_r \subseteq M$ and of conjugacy classes $\lambda_1, \ldots, \lambda_r \subseteq G$, called the holonomies. In that case, we can also define the parabolic $G$-representation variety, $\Rep{G}(M, A, Q)$, as the set of representations $\rho: \Pi(M - \cup S_i, A) \to G$ such that $\rho(\gamma_i) \in \lambda_i$ if $\gamma_i$ is a loop around $S_i$. As in the non-parabolic case, it has a natural structure of a complex variety.

There is a natural action of $G$ on $\Rep{G}(M, A, Q)$ by conjugation. Thus, if $G$ is also a reductive group, we can also consider the GIT quotient $\Char{G}(M, A, Q) = \Rep{G}(M, A, Q) \sslash G$. This algebraic variety is the so-called parabolic $G$-character variety and it can be shown to be the moduli space of such parabolic representations.

The understanding of the topological and algebraic structure of these character varieties is an active area of research. A starting point is to compute their Betti numbers or, equivalently, their Poincar\'e polynomial for the case $M=\Sigma$, a closed oriented manifold, a single basepoint and no parabolic structure. In this direction, the non-abelian Hodge theory becomes a useful tool. Roughly speaking, this theory shows that three different moduli spaces on $\Sigma$ are diffeomorphic: the character variety, the moduli space of flat connections and the moduli space of Higgs bundles. For a precise description of these spaces and these equivalences, see \cite{SimpsonI, SimpsonII}, \cite{Simpson:1992} and \cite{Corlette:1988}.

In this way, the topology of $\Char{G}(\Sigma)$ can be studied via a natural perfect Morse-Bott function on the moduli space of Higgs bundles. As a result, the Poincar\'e polynomial of character varieties has been computed for $G=\SL{2}(\CC)$ \cite{Hitchin}, for $G=\SL{3}(\CC)$ \cite{Gothen} and for $G=\GL{4}(\CC)$ \cite{GP-Heinloth-Schmitt}.  In general, in \cite{Schiffmann:2016} and \cite{Mozgovoy-Schiffmann} (see also \cite{Mellit}) a combinatorial formula is given for arbitrary $G=\GL{n}(\CC)$ provided that $n$ and the degree of the bundle are coprime. In the parabolic case, it has been computed for $G=\SL{2}(\CC)$ and generic semi-simple conjugacy classes in \cite{Boden-Yokogawa} and for $G=\GL{3}(\CC),\SL{3}(\CC)$ in \cite{GP-Gothen-Munoz}. 

However the equivalences from non-abelian Hodge theory are not holomorphic, so the complex structures on the character variety and on the moduli spaces of flat connections and Higgs bundles do not agree. For this reason, the study of specific algebraic invariants of character varieties turns out to be important. An relevant invariant of a complex variety $X$ is its Deligne-Hodge polynomial, $\DelHod{X} \in \ZZ[u^{\pm 1}, v^{\pm 1}]$, which is an alternating sum of the Hodge numbers of $X$ as in the spirit of an Euler characteristic.

Two different approaches have been used to compute the Deligne-Hodge polynomials of character varieties over a closed orientable surface. The first one, called the arithmetic method, was introduced in \cite{Hausel-Rodriguez-Villegas:2008}. It is based on a theorem of Katz inspired by the Weil conjectures. Using it, several computations of Deligne-Hodge polynomials have been done, as in \cite{Hausel-Rodriguez-Villegas:2008} for $G = \GL{n}(\CC)$ and a parabolic structure with a single marked point and holonomy a primitive root of unit, and for $G=\SL{n}(\CC)$ in \cite{Mereb}. However, in the parabolic case, very little has been advanced, being the major achievement the computation of these polynomials for $G = \GL{n}(\CC)$ and a generic parabolic structure in \cite{Hausel-Letellier-Villegas}.

In order to overcome this problem, a new geometric method was introduced in \cite{LMN}. The idea is to stratify the representation variety into simpler pieces for which the Deligne-Hodge polynomial can be computed easily and to add all the contributions. From this polynomial, the one of the character variety follows by analyzing the identifications that take place in the GIT quotient. Using this method, explicit expressions of these polynomials have been computed for the case $G=\SL{2}(\CC)$. They were computed for at most one marked point and closed orientable surfaces of genus $g = 1, 2$ in \cite{LMN}, for $g = 3$ in \cite{MM:2016} and for arbitrary genus in \cite{MM}. In \cite{LM}, the case of two marked points is accomplished for $g=1$. Finally, in \cite{Baraglia-Hekmati:2016}, a mix between the arithmetic and the geometric method was used, computing explicit expressions of the Deligne-Hodge polynomials for $G = \SL{2}(\CC), \SL{3}(\CC)$, arbitrary genus and no parabolic structure. Despite its success, this method requires the use of very specific stratifications that are not clear how to generalize to other groups or more general parabolic structures.

Motivated by this problem, in \cite{GPLM-2017} a new categorical approach to this problem was introduced as part of the PhD Thesis \cite{Gonzalez-Prieto:Thesis} of the author. In that paper, given a complex algebraic group $G$, it was constructed a lax monoidal Topological Quantum Field Theory (TQFT for short) that computes the Deligne-Hodge polynomials of $G$-representation varieties. Recall that a lax monoidal TQFT is a lax monoidal symmetric functor $\Zs{G}: \Bordppar{n}{\Lambda} \to \Mod{\K{\MHSq}}$. Here $\K{\MHSq}$ is the Grothendieck ring (aka.\ the $K$-theory) of the category of Hodge structures, $\Bordppar{n}{\Lambda}$ is the category of parabolic bordisms of pairs and $\Mod{\K{\MHSq}}$ is the category of $\K{\MHSq}$-modules. The crucial property of this TQFT is that, if $W$ is a closed $n$-dimensional manifold, $A \subseteq W$ is a finite set and $Q$ is a parabolic structure on $W$, then $\Zs{G}(W, A, Q): \K{\MHSq} \to \K{\MHSq}$ is a $\K{\MHSq}$-linear morphism given by the multiplication by the image in $K$-theory of the Hodge structure on the cohomology of $\Rep{G}(W, A, Q)$.

The aim of this paper is to improve $\Zs{G}$ for computing, not only the Hodge structure, but the whole virtual class of the representation variety in the Grothendieck ring of algebraic varieties. It will give rise to a lax monoidal TQFT, $\Zs{G}: \Bordppar{n}{\Lambda} \to \Mod{\KVar{k}}$, with $\KVar{k}$ the Grothendieck ring of algebraic varieties over an arbitrary field $k$. Moreover, we will show that $\Zs{G}$ can be used to give an effective recursive method of computation of these virtual classes, even in the parabolic case. As an application, we will compute it for $G=\SL{2}(\CC)$, arbitrary genus and any number of marked points of Jordan type. This recursive nature for character varieties has also being explored in the literature, as in and \cite{Carlsson-Rodriguez-Villegas}, \cite{Diaconescu:2017}, \cite{Hausel-Letellier-Villegas:2013} and \cite{Mozgovoy:2012}. On the other hand, the computation of the virtual classes of representation varieties seems to be a harder problem, being a major contribution the work \cite{Wyss:2017} for quiver varieties.

%In Section \ref{sec:hodge-theory} of this paper, we give a brief review of Hodge theory, with special attention to Saito's theory of mixed Hodge modules (Section \ref{subsec:mixed-hodge-modules}). This will be very useful since a fundamental piece in the construction of $\Zs{G}$ is the category of mixed Hodge modules over an algebraic variety $X$, denoted $\MHM{X}$. In Section \ref{sec:monodromy-as-mhm}, we review in deeper detail the relation of mixed Hodge modules with variations of mixed Hodge structures. Using this interplay, we give a concrete incarnation of some mixed Hodge modules as monodromies of locally trivial fibrations that will be necessary for computational purposes.
 
In Section \ref{sec:TQFT-physical-constr} of this paper, we extend the results of \cite{GPLM-2017} and describe a method of construction of TQFTs from two simpler pieces: a field theory and a quantisation. The idea is to consider an auxiliar category $\cC$ with pullbacks and final object $\star \in \cC$, that is going to play the role of a category of fields (in the physical sense). Let $\Embc$ be the category whose objects are compact manifolds (maybe with boundary) but only open embeddings between them.  A functor $F: \Embc \to \cC$ with the Seifert-van Kampen property (see Definition \ref{def:seifert-vankampen}) gives rise to a `field theory', which is a functor $\Fld{F} \to \Span{\cC}$.

On the other hand, in Section \ref{sec:C-algebras} we introduce an algebraic object called a $\cC$-algebra. Roughly speaking a $\cC$-algebra $\cA$ is a collection of algebras $\cA_c$ for all $c \in \cC$ that preserves the functorial structure of $\cC$. From a $\cC$-algebra, in Section \ref{sec:quantisation} we show how to construct a functor $\Qtm{\cA}: \Span{\cC} \to \Mod{\cA_\star}$ that plays the role of a `quantisation' of the fields of $\cC$. Composing the functors from the field theory part and the quantisation, we obtain a functor $Z = \Qtm{\cA} \circ \Fld{F}: \Bord{n} \to \Span{\cC} \to \Mod{\cA_\star}$ that we will show it is a lax monoidal TQFT. The same construction can be done if we consider an extra structure on the category of bordisms given by a sheaf, as finite configurations of points or parabolic structures (see Section \ref{sec:TQFT-physical-constr}). Therefore, putting together all this information, we obtain the following result (Theorem \ref{thm:physical-constr-sheaf}).

\vspace{-0.1cm}

\begin{theorem*}
Let $\cC$ be a category with final object $\star$ and pullbacks and let $\cS$ be a sheaf. Given a functor $F: \EEmbc{\cS} \to \cC$ with the Seifert-van Kampen property and a $\cC$-algebra $\cA$, there exists a lax monoidal Topological Quantum Field Theory over $\cS$, $
	Z_{F, \cA}: \EBord{n}{\cS} \to \Modt{\cA_\star}$.
\end{theorem*}

\vspace{-0.2cm}
In Section \ref{sec:grothendieck-rings}, we show how the slice categories of algebraic varieties, $\Varrel{X}$, can be used as the quantisation part of a lax monoidal TQFT. To be precise, we shall take $\cC = \Var{k}$, the category of $k$-algebraic varieties, and we will show that the Grothendieck rings of these slice categories can be made into a $\Var{k}$-algebra, $\KVarrel{}$, in a natural way. In Section \ref{sec:Hodge-monodromy-covering-spaces}, we shall explain how this $\Var{k}$-algebra can be understood as a promotion of the Hodge monodromy representation introduced in \cite{LMN} for the case of coverings. Some preliminary computations are also provided that will be useful in Section \ref{sec:sl2-repr-var}.

Section \ref{sec:repr-varieties} of this paper is devoted to representation varieties. In Section \ref{sec:TQFT-for-repr} we review the construction of the functor $\Zs{G}$ of \cite{GPLM-2017} in the context of the physical construction introduced in this paper. In this case, the auxiliary category of fields will be $\cC = \Var{k}$, the functor $F$ will be the representation variety functor and the quantisation will be induced by the $\Var{k}$-algebra $\KVarrel{}$. However, even though $\Zs{G}$ computes the virtual classes of representation varieties, for computational purposes it is better to consider a modified TQFT, called the geometric TQFT, and denoted $\Zg{G}$. It is built in Section \ref{subsec:piecewise-alg-var} by means of a procedure called reduction. The advantage of $\Zg{G}$ is that, in general, $\Zs{G}(S^1, \star) = \KVarrel{G}$ is not a finitely generated module, but $\Zg{G}(S^1, \star)$ is.

As an application, in Section \ref{sec:sl2-repr-var} we show how $\Zg{G}$ can be used to give an effective method of computation of virtual classes on representation varieties.
We focus on the case $G = \SL{2}(\CC)$ and we will allow arbitrary many marked points with holonomies in $\Lambda = \left\{[J_+], [J_-], \left\{-\Id\right\}\right\}$. Here, $[J_+], [J_-] \subseteq \SL{2}(\CC)$ are the set of Jordan-type matrices of traces $2$ and $-2$ respectively.
For that purpose, in Section \ref{subsec:core-submodule} we give an explicit finite set of generators for $\Zg{}(S^1, \star)$ in terms of monodromies of coverings.

Observe that any closed surface $\Sigma$, endowed with $A \subseteq \Sigma$ finite and $Q$ a parabolic structure, is a composition of the bordisms depicted in the Figure \ref{fig:img-tubes-pairs-short} below, with $\lambda = [J_+], [J_-]$ or $\left\{-\Id\right\}$. Hence, in order to compute the virtual class of $\Rep{\SL{2}(\CC)}(\Sigma, A, Q)$, it is enough to compute the homomorphisms $\Zg{}(D), \Zg{}(D^\dag), \Zg{}(L)$ and $\Zg{}(L_\lambda)$.

\begin{figure}[h]
\label{fig:img-tubes-pairs-short}
	\begin{center}
	\includegraphics[scale=0.4]{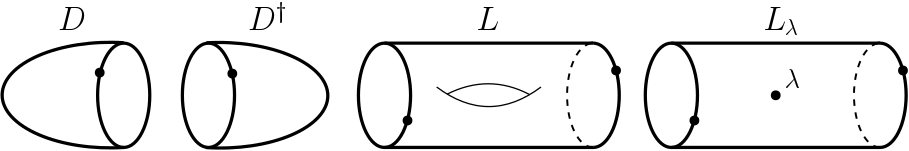}
	\end{center}
%	\vspace{-0.5cm}
	\caption{Basic bordisms for surfaces.}
\vspace{-0cm}
\end{figure}

Towards this aim, Section \ref{sec:discs-tubes} is devoted to the computation of the morphisms $\Zg{}(D)$, $\Zg{}(D^\dag)$ and $\Zg{}(L_{-\Id})$. In Section \ref{sec:tube-J+}, we compute the more involved homomorphisms $\Zg{}(L_{[J_+]})$ and $\Zg{}(L_{[J_-]})$. Finally, in Section \ref{sec:genus-tube} we will show how the results of \cite{MM} can be reformulated to give $\Zg{}(L)$. Putting together all this information, we obtain the main result of this paper.

\begin{theorem*}
Let $Q$ be a parabolic structure on $\Sigma_g$ with $r$ marked points, of which $\ell$ are of type $[J_-]$ or $\set{-\Id}$. Then, virtual class of $\Rep{\SL{2}(\CC)}(\Sigma_g, Q)$ in the Grothendieck ring of algebraic varieties is
\begin{itemize}
	\item If $r=0$,
\begin{align*}
\coh{\Rep{\SL{2}(\CC)}(\Sigma_g)} =& \,{\left(q^2 - 1\right)}^{2g - 1} q^{2g - 1} +
\frac{1}{2} \, {\left(q -
1\right)}^{2g - 1}q^{2g -
1}(q+1){\left({2^{2g} + q - 3}\right)} 
\\ &+ \frac{1}{2} \,
{\left(q + 1\right)}^{2g + r - 1} q^{2g - 1} (q-1){\left({2^{2g} +q -1}\right)} + q(q^2-1)^{2g-1}.
\end{align*}
	\item If $r>0$ and $\ell$ is even,
\begin{align*}
	\coh{\Rep{\SL{2}(\CC)}(\Sigma_g, Q)} =&\, {\left(q^2 - 1\right)}^{2g + r - 1} q^{2g - 1} +
\frac{1}{2} \, {\left(q -
1\right)}^{2g + r - 1}q^{2g -
1}(q+1){\left({2^{2g} + q - 3}\right)} 
\\ &+ \frac{\left(-1\right)^{r}}{2} \,
{\left(q + 1\right)}^{2g + r - 1} q^{2g - 1} (q-1){\left({2^{2g} +q -1}\right)}.
\end{align*}
	\item If $r=\ell = 1$,
\begin{align*}
	\coh{\Rep{\SL{2}(\CC)}(\Sigma_g, Q)} =& \,{\left(q - 1\right)}^{2g + r - 1} (q+1)q^{2g - 1}{{\left( {\left(q + 1\right)}^{2 \,
g + r-2}+2^{2g-1}-1\right)} } \\
&+ \left(-1\right)^{r + 1}2^{2g - 1}  {\left(q + 1\right)}^{2g + r
- 1} {\left(q - 1\right)} q^{2g - 1} + q(q^2-1)^{2g-1}.
\end{align*}
	\item If $\ell>1$ is odd,
\begin{align*}
	\coh{\Rep{\SL{2}(\CC)}(\Sigma_g, Q)} = &\, {\left(q - 1\right)}^{2g + r - 1} (q+1)q^{2g - 1}{{\left( {\left(q + 1\right)}^{2 \,
g + r-2}+2^{2g-1}-1\right)} } \\
&+ \left(-1\right)^{r + 1}2^{2g - 1}  {\left(q + 1\right)}^{2g + r
- 1} {\left(q - 1\right)} q^{2g - 1}.
\end{align*}
\end{itemize}
Here, $q = \coh{\CC} \in \KVar{\CC}$ denotes the Lefschetz motif.
\end{theorem*}

At this point, several research lines can be addressed as continuation of this work. The first one is to study how these computations of virtual classes of representation varieties can be used to give the corresponding ones for the associated character varieties. This gap is filled by the paper \cite{GP-2018b} (see also \cite{Gonzalez-Prieto:Thesis}).

Moreover, a wider class of punctures can be considered, as of semi-simple type. However, in that case, new generators of $\Zg{}(S^1, \star)$ appear to track the semi-simple holonomies. These generators also exhibit a complex interaction phenomenon that requires a subtler analysis. This objective goes beyond the scope of the present paper and it has been postponed to the upcoming paper \cite{GP-2019}.

A step further, it would be interesting to use the construction method described here to extend the TQFT for character varieties across the non-abelian Hodge correspondence. This will allow us to build analogous TQFTs for the moduli spaces of flat bundles and Higgs bundles that would help to understand the whole hyperk\"ahler structure.

Finally, at a long term, we want to study the relation between $\Zs{G}$ and $\Zs{{{^L}G}}$, where ${^L}G$ is the Langlands dual group of $G$. The reason are the conjectures formulated in \cite{Hausel:2005} and \cite{Hausel-Rodriguez-Villegas:2008} in the context of the mirror symmetry conjectures and the geometric Landglands proglam (see \cite{Beilinson-Drinfeld}). They conjecture the existence of very astonishing symmetries between the Deligne-Hodge polynomials for character varieties over $G$ and ${^L}G$. We hope that the kind of ideas introduced in this paper will help to shed some light over these problems in the future.

\subsection*{Acknowledgements}

The author wants to thank David Ben-Zvi, Georgios Charalambous, Tam\'as Hausel, Marta Mazzocco, Ignasi Mundet, Peter Newstead and Nils Prigge for very useful conversations. I would like to specially thank Thomas Wasserman for very interesting and fruitful conversations about TQFTs. Finally, I would like to express my highest gratitude to my PhD advisors Marina Logares and Vicente Mu\~noz for their invaluable help, support and encouragement throughout the development of this paper.

The author acknowledges the hospitality of the Facultad de Ciencias Matem\'aticas at Universidad Complutense de Madrid and the Centre for Mathematical Sciences at the University of Plymouth, in which part of this work was completed. The author has been supported by a "\!la Caixa" scholarship LCF/BQ/ES15/10360013 for PhD studies in Spanish Universities, by MINECO (Spain) Project MTM2015-63612-P and by a Severo Ochoa Postdoctoral Fellowship SOLAUT\_00030167 at Instituto de Ciencias Matem\'aticas.

%%%%%%%%%%%%%%%%%%%%%%%%%%%%%%%%%%%%%%%%%%%%%%%%%%%%%%%%%%%%%%%%
%%%%%%%%%%%%%%%%%% SECTION: LTQFTp %%%%%%%%%%%%%%%%%%%%%%%%%%%%%%
%%%%%%%%%%%%%%%%%%%%%%%%%%%%%%%%%%%%%%%%%%%%%%%%%%%%%%%%%%%%%%%%

\section{Topological Quantum Field Theories}
\label{sec:TQFT-physical-constr}

In this section, we will describe a general recipe for constructing lax monoidal TQFTs from two simpler pieces of data: one of geometric nature (a field theory) and one of algebraic nature (a quantisation). It was first described in \cite{GPLM-2017} in the context of representation varieties. However, for the purposes of this paper we need a more general setting.

The idea is a very natural construction that, in fact, has been widely used in the literature in some related form (see for example \cite{Ben-Zvi-Gunningham-Nadler, Ben-Zvi-Nadler:2016, Freed-Hopkins-Lurie-Teleman:2010, Haugseng}), sometimes referred to as the `push-pull construction'. In this paper, we recast this construction to identify the requiered input data in a simple way that will be useful for applications. For example, it fits perfectly with the Grothendieck ring of algebraic varieties as in Section \ref{sec:grothendieck-rings}.

Let us consider the category $\Embc$ whose objects are compact differentiable manifolds, maybe with boundary. Given compact manifolds $M_1$ and $M_2$ of the same dimension, a morphism $M_1 \to M_2$ in $\Embc$ is a class of smooth tame embeddings $f: M_1 \to M_2$. The tameness condition means that, for some union of connected components $N \subseteq \partial M_1$, $f$ extends to a smooth embedding $\tilde{f}: M_1 \cup_N \left(N \times [0,1]\right) \to M_2$ with $\tilde{f}^{-1}(\partial M_2) = \partial M_1 - N$. 

Two such a tame embeddings $f, f': M_1 \to M_2$ are declared equivalent if there exists an ambient diffeotopy $h$ between them i.e.\ a smooth map $h: M_2 \times [0, 1] \to M_2$ such that $h_t = h(-,t)$ is a diffeomorphism of $M_2$ for all $0 \leq t \leq 1$, $h_0 = \Id_{M_2}$, $h_t|_{\partial M_2} = \Id_{\partial M_2}$ for all $t$, and $h_1 \circ f = f'$. There are no morphisms in $\Embc$ between manifolds of different dimensions.

\begin{ex}\label{ex:embc-map-boundary}
Let $M$ be a compact manifold and let $N \subseteq \partial M$ be a union of connected components of $\partial M$. Let $U \subseteq M$ be an open collaring around $N$, that is, an open subset of $M$ such that there exists a diffeomorphism $\varphi: N \times [0,1) \to U$. In that case, $\varphi: N \times [0, 1/2] \to U \subseteq M$ defines a morphism in $\Embc$. Moreover, such a morphism does not depend on the chosen collaring as any two collars are ambient diffeotopic. For this reason, we will denote this map $N \times [0, 1/2] \to M$ in $\Embc$ just by $N \to M$.
\end{ex}

Consider a sheaf $\cS: \Embc \to \Cat$ with $\cS(\emptyset)=\emptyset$. By a sheaf we mean a contravariant functor such that, given a compact manifold $M$ and a covering $\left\{W_\alpha\right\}_{\alpha \in \Lambda}$ of $M$ by compact submanifolds, for any collection $s_\alpha \in \cS(W_\alpha)$ with $\cS(i_{\alpha,\beta})(s_\alpha) = \cS(i_{\beta,\alpha})(s_\beta)$ for all $\alpha,\beta \in \Lambda$ (where, $i_{\alpha,\beta}: W_\alpha \cap W_\beta \to W_\alpha$ is the inclusion) there exists an unique $s \in \cS(M)$ such that $\cS(i_\alpha)(s) = s_\alpha$, $i_\alpha: W_\alpha \to M$ being the inclusion maps.

\begin{rmk}
\begin{itemize}
	\item In contrast with usual sheaves, the coverings considered here are closed. However, they also give an open covering since the embeddings $i_\alpha$ have to extend to a small open set around $M_\alpha$ by the tameness condition.
	\item We should think about these sheaves as a way of endowing the bordisms with and geometric extra structure. The formulation of these extra structures as sheaves goes back to \cite{Madsen-Weiss:2007, Galatius-Tillmann-Madsen-Weiss:2009}.
\end{itemize}
\end{rmk}

Given such a sheaf $\cS$, we define the \emph{category of embeddings over $\cS$}, or just $\cS$-embeddings, $\EEmbc{\cS}$. The objects of this category are pairs $(M, s)$ with $M \in \Embc$ and $s \in \cS(M)$ and the initial object $\emptyset$. Observe that, if $M \neq \emptyset$ has $\cS(M)=\emptyset$, then $M$ does not appear as object of $\EEmbc{\cS}$. Given objects $(M_1, s_1), (M_2, s_2) \in \EEmbc{\cS}$, a morphism between them is a pair $(f, \alpha)$ where $f: M_1 \to M_2$ is a morphism in $\Embc$ and $\alpha: s_1 \to \cS(f)(s_2)$ is a morphism of $\cS(M_1)$. If $\alpha = 1_{s_1}$ (so $\cS(f)(s_2)=s_1$) we will just denote the morphism by $f: (M_1, s_1) \to (M_2, s_2)$. The composition of two morphisms $(f, \alpha)$ and $(f', \alpha')$ is $(f' \circ f, \cS(f)(\alpha') \circ \alpha)$.

Moreover, given $n \geq 1$, we define the \emph{category of $n$-bordisms over $\cS$}, $\EBord{n}{\cS}$. It is a $2$-category given by the following data:
\begin{itemize}
	\item Objects: The objects of $\EBord{n}{\cS}$ is the subclass of objects $(M, s) \in \EEmbc{\cS}$ with $M$ a $(n-1)$-dimensional closed manifold. 
	\item $1$-morphisms: Given objects $(M_1, s_1)$, $(M_2, s_2)$ of $\EBord{n}{\cS}$, a morphism $(M_1, s_1) \to (M_2, s_2)$ is a pair $(W, s) \in \EEmbc{\cS}$ where $W$ is a compact $n$-dimensional manifold such that $\partial W = M_1 \sqcup M_2$ (i.e.\ an unoriented bordism between $M_1$ and $M_2$). It also has to satisfy that $\cS(i_1)(s)=s_1$ and $\cS(i_2)(s)=s_2$, being $i_k: M_k \to W$ the inclusions for $k =1,2$.
	%Two bordisms $(W, s)$ and $(W', s')$ are in the same class if there exists a boundary preserving diffeomorphism $f: W \to W'$ such that $\cS(f)(s')=s$.
	\\
With respect to the composition, given $(W, s): (M_1, s_1) \to (M_2, s_2)$ and $(W', s'): (M_2, s_2) \to (M_3, s_3)$, we define $(W',s') \circ (W,s)$ as the morphism $(W \cup_{M_2} W', s \cup s'): (M_1, s_1) \to (M_3, s_3)$ where $W \cup_{M_2} W'$ is the usual gluing of bordisms along $M_2$ and $s \cup s' \in \cS(W \cup_{M_2} W')$ is the object given by gluing $s$ and $s'$ with the sheaf property.
	\item $2$-morphisms: Given two $1$-morphisms $(W, s), (W',s'): (M_1, s_1) \to (M_2, s_2)$, a $2$-cell $(W, s) \Rightarrow (W',s')$ is a morphism $(f, \alpha): (W, s) \to (W',s')$ of $\EEmbc{\cS}$ with $f$ a boundary preserving diffeomorphism. Composition of $2$-cells is just composition in $\EEmbc{\cS}$.
\end{itemize}

In this form, $\EBord{n}{\cS}$ is not a category since, for some $(M,s) \in \EBord{n}{\cS}$, it might happen that $\Hom_{\EBord{n}{\cS}} ((M,s),(M,s))$ has no unit morphism. This problem can be solved by slightly weakening the notion of bordism, allowing that $M$ itself can be seen as a bordism $M: M \to M$. With this modification, $(M,s): (M,s) \to (M,s)$ is the desired unit and it is a straightforward check to see that $\EBord{n}{\cS}$ is a (strict) $2$-category. Furthermore, it has a natural monoidal structure by means of disjoint union of objects and bordisms.

%Let $R$ be a fixed commutative and unitary ring. Another important category is the \emph{category of $R$-modules with twists}, $\Modt{R}$. The objects of $\Modt{R}$ are $R$-modules. A $1$-morphism between modules $M$ and $N$ is given by a triple $(L, f, g)$, where $L$ is a $R$-module and $f: M \to L$ and $g: L \to N$ are $R$-modules homomorphisms.

Another important category is the category of algebras with twists. Let $R$ be a fixed commutative and unitary ring. Let $M$ and $N$ be two $R$-algebras and suppose that $f,g: M \to N$ are two homomorphisms as $R$-modules. Suppose that there exist decompositions of $f, g$ given by $R$-algebras $L, L'$, $R$-algebra homomorphisms $f^{alg}: M \to L, g^{alg}: M \to L'$ and $R$-module homomorphisms $f^{mod}: L \to N, g^{mod}: L' \to N$ such that $f = f^{mod} \circ f^{alg}$ and $g=g^{mod} \circ g^{alg}$. We say that $g$ is an \emph{immediate twist} of $f$ if there exist a $R$-algebra homomorphism $\alpha^{alg}: L \to L'$ and a $R$-module homomorphism $\alpha^{mod}: L' \to L$ such that $g^{alg} = \alpha^{alg} \circ f^{alg}$ and $g^{mod} = f^{mod} \circ \alpha^{mod}$.
\begin{comment}
homomorphisms $f_1: M \to L$, $f_2: L \to N$ and $\psi: L \to L$ such that $f = f_2 \circ f_1$ and $g = f_2 \circ \psi \circ f_1$.
\[
\begin{displaystyle}
   \xymatrix
   {	M \ar[r]^{f_1} \ar@/_1pc/[rr]_g & L \ar[r]^{f_2} \ar@(ul,ur)^{\psi} & N
   }
\end{displaystyle}   
\]
\end{comment}
\[
\begin{displaystyle}
   \xymatrix
   { && L \ar[rrd]^{f^{mod}} \ar@/_0.8pc/[dd]_{\alpha^{alg}}  && \\
   M \ar[rru]^{f^{alg}}\ar[rrd]_{g^{alg}} &&  && N \\
   && L'\ar[rru]_{g^{mod}} \ar@/_0.8pc/[uu]_{\alpha^{mod}} &&
   }
\end{displaystyle}   
\]
The data of an immediate twist is given by the tuple $(L, L', f^{alg}, f^{mod}, g^{alg}, g^{mod}, \alpha^{alg}, \alpha^{mod})$.
In general, given $f,g: M \to N$ two $R$-module homomorphisms, we say that $g$ is a \emph{twist} of $f$ if there exists a finite sequence $f=f_0, f_1, \ldots, f_r = g: M \to M$ of homomorphisms such that $f_{i+1}$ is an immediate twist of $f_i$ for $0 \leq i \leq r-1$.

To these sequences, we add the relation that if the immediate twist between $f_{k-1}$ and $f_{k}$ is given by $(L, L', f_{k-1}^{alg}, f_{k-1}^{mod}, f_{k}^{alg}, f_{k}^{mod}, \alpha^{alg}, \alpha^{mod})$ and the one between $f_k$ and $f_{k+1}$ is given by $(L', L'', f_{k}^{alg}, f_{k}^{mod}, f_{k+1}^{alg}, f_{k+1}^{mod}, \beta^{alg}, \beta^{mod})$ then the sequence $f_0, \ldots, f_{k-1}, f_k, f_{k+1}, \ldots f_r$ is equivalent to $f_0, \ldots, f_{k-1}, f_{k+1}, \ldots, f_r$ with the immediate twist between $f_{k-1}$ and $f_{k+1}$ given by $(L, L'', f_{k-1}^{alg}, f_{k-1}^{mod}, f_{k+1}^{alg}, f_{k+1}^{mod}, \beta^{alg} \circ \alpha^{alg}, \alpha^{mod} \circ \beta^{mod})$.

In that case, we define the \emph{category of $R$-algebras with twists}, $\Modt{R}$, as the category whose objects are $R$-algebras, its $1$-morphisms are $R$-modules homomorphisms and, given homomorphisms $f$ and $g$, a $2$-morphism $f \Rightarrow g$ is a twist from $f$ to $g$. Composition of $2$-cells is juxtaposition of twists. With this definition, $\Modt{R}$ has a $2$-category structure. Moreover it is a monoidal category with the usual tensor product. Thus, the forgetful functor $\Modt{R} \to \Mod{R}$ to the usual category of $R$-modules is monoidal.

\begin{defn}
Let $\cS: \Embc \to \Cat$ be a sheaf and let $R$ be a ring. A \emph{(lax monoidal) Topological Quantum Field Theory over $\cS$}, shortened (lax monoidal) $\cS$-TQFT, is a (lax) monoidal symmetric $2$-functor $$
Z: \EBord{n}{\cS} \to \Modt{R}.
$$
\end{defn}

\begin{rmk}
Recall that a functor between monoidal categories $F: (\cC, \otimes_\cC, I_\cC) \to (\cD, \otimes_\cD, I_\cD)$ is said to be lax monoidal if there exists:
\begin{itemize}
		\item A morphism $\alpha: I_{\cD} \to F(I_{\cC})$.
		\item A natural transformation $\Delta: F(-) \otimes_{\cD} F(-) \Rightarrow F(- \otimes_\cC -)$.		
	\end{itemize}
If $\alpha$ and $\Delta$ are isomorphisms, $F$ is said to be \emph{pseudo-monoidal} (or just monoidal) and, if they are identity morphisms, $F$ is called \emph{strict monoidal}.
\end{rmk}

The important point is that there exists a physically inspired construction of lax monoidal TQFTs that allows us to construct them from simpler data. The idea of the construction is to consider an auxiliar category $\cC$ with pullbacks and final object, that is going to play the role of a category of fields (in the physical sense). Then, we are going to split our functor $Z$ as a composition
$$
	\Bord{n} \stackrel{\Fld{}}{\longrightarrow} \Span{\cC} \stackrel{\Qtm{}}{\longrightarrow} \Mod{R},
$$
where $\Span{\cC}$ is the category of spans of $\cC$ (see Section \ref{sec:quantisation}). The first arrow, $\Fld{}$, is the field theory and we will describe how to built it in Section \ref{sec:field-theory}. The second arrow, $\Qtm{}$, is the quantisation part. It will be constructed by means of an algebraic tool called a $\cC$-algebra (see Section \ref{sec:quantisation}). A $\cC$-algebra can be thought as collection of rings parametrized by $\cC$ with a pair of induced homomorphisms from every morphism of $\cC$.

\subsection{$\cC$-algebras}
\label{sec:C-algebras}
Let $\cC$ be a cartesian monoidal category (i.e.\ a category with finite products where the monoidal structure is precisely to take products) and let $A: \cC \to \Rng$ be a contravariant functor. If $\star \in \cC$ is the final object, the ring $A(\star)$ plays an special role. Note that, for every $a \in \cC$, we have an unique map $c_a: a \to \star$ which gives rise to a ring homomorphism $A(c_a): A(\star) \to A(a)$. Hence, we can see $A(a)$ as a $A(\star)$-algebra in a natural way. Such a module structure is the one considered in the first condition of Definition \ref{def:C-algebra} below.

Given $a, b,d \in \cC$ with morphisms $a \to d$ and $b \to d$, let $p_1: a \times_d b \to a$ and $p_2: a \times_d b \to b$ be the corresponding projections. We define the \emph{external product} over $d$
$$
	\boxtimes_d: A(a) \otimes_{A(d)} A(b) \to A(a \times_d b),
$$ by $z \boxtimes_d w = A(p_1)(z) \cdot A(p_2)(w)$ for $z \in A(a)$ and $w \in A(b)$. The external product over the final object will be denoted just by $\boxtimes = \boxtimes_\star: A(a) \otimes_{A(\star)} A(b) \to A(a \times b)$. It gives a natural transformation $\boxtimes: A(-) \otimes_{A(\star)} A(-) \Rightarrow A(- \times -)$.

\begin{rmk}\label{rmk:module-exterior-prod}
For $b = \star$, we have that $p_1 = \lambda: a \times \star \to a$ is the unital isomorphism of the monoidal structure so it gives raise to an isomorphism $A(\lambda): A(a) \to A(a \times \star)$ of $A(\star)$-algebras. Under this isomorphism, the external product $\boxtimes: A(a) \otimes A(\star) \to \cA(a \times \star) \cong A(a)$ coincides with the given $A(\star)$-module structure on $A(a)$.
\end{rmk}

\begin{defn}\label{def:C-algebra}
Let $(\cC, \times, \star)$ be cartesian monoidal category. A \emph{$\cC$-algebra}, $\cA$, is a pair of functors
$$
	A: \cC^{op} \to \Rng, \hspace{1cm} B: \cC \to \Mod{A(\star)},
$$
such that:
\begin{itemize}
	\item They agree on objects, that is, $A(a) = B(a)$ for all $a \in \cC$, as $A(\star)$-modules.
	\item They satisfy the Beck-Chevaley condition (also known as the base change formula), that is, given $a_1, a_2, b \in \cC$ and a pullback diagram
\[
\begin{displaystyle}
   \xymatrix
   {
	d \ar[r]^{\;g'}\ar[d]_{f'} & a_1\ar[d]^{f} \\
	a_2 \ar[r]_{g} & b
   }
\end{displaystyle}
\]
we have that $A(g) \circ B(f) = B(f') \circ A(g')$.
	\item The external product $\boxtimes: B(-) \otimes_{A(\star)} B(-) \Rightarrow B(- \times -)$ is a natural transformation with respect to $B$.
\end{itemize}
\end{defn}

\begin{rmk}\label{rmk:notation-grothendieck-six}
\begin{itemize}
	\item A $\cC$-algebra can be thought as a collection of algebras parametrized by $\cC$. For this reason, we will denote $\cA_a = A(a)$ for $a \in \cC$.
	\item The Beck-Chevalley condition appears naturally in the context of Grothendieck's yoga of six functors $f_*, f^*, f_!, f^!, \otimes$ and $\mathbb{D}$ in which $(f_*, f^*)$ and $(f_!, f^!)$ are adjoints, and $f^*$ and $f_!$ satisfy the Beck-Chevalley condition. 
In this context, we can take $A$ to be the functor $f \mapsto f^*$ and $B$ the functor $f \mapsto f_!$. Moreover, in order to get in touch with this framework, we will denote $A(f) = f^*$ and $B(f) = f_!$.
\end{itemize}
\end{rmk}

Using the covariant functor $B$, for every object $a \in \cC$, we obtain a $\cA_\star$-module homomorphism $(c_a)_!: \cA_a \to \cA_\star$. The special element $\mu(a) = (c_a)_!(1) \in \cA_\star$, where $1 \in \cA_a$ is the unit of the ring, will be called the \emph{characteristic} of $a$.

\begin{ex}\label{ex:algebra-sheaves}
Given a locally compact Hausdorff topological space $X$, let us consider the category $\Sh{X}$ of sheaves (of rational vector spaces) on $X$ with sheaf transformations between them. It is an abelian monoidal category with monoidal structure given by tensor product of sheaves. Given a continuous map $f: X \to Y$, we can induce two special maps at the level of sheaves. The first one is the inverse image $f^*: \Sh{Y} \to \Sh{X}$ and it is an exact functor. We also have the direct image with compact support functor, $f_!: \Sh{X} \to \Sh{Y}$. In this case, $f_!$ is only left exact, so we can consider its derived functor $Rf_!: \Sh{X} \to \Dp{\Sh{Y}}$ whose stalks are $(R^kf_!(\cF))_y = H_c^k(f^{-1}(y), \cF)$ for $y \in Y$ and $\cF$ a sheaf on $X$. In this context, the base change theorem with compact support \cite[Theorem 2.3.27]{Dimca:2004} implies that, for any pullback diagram of locally Hausdorff spaces
\[
\begin{displaystyle}
   \xymatrix
   {
	X \times_Z Y \ar[r]^{\;\;\;\;\;\;g'} \ar[d]_{f'} & X \ar[d]^{f} \\
	Y \ar[r]_{g} & Z
   }
\end{displaystyle}
\]
there is a natural isomorphism $g^* \circ Rf_! \cong Rf'_! \circ g'^*$. Even more, $Rf_!$ preserves the external product.

This situation can be exploited to obtain a $\cC$-algebra. However we must surround the difficulty that, for a general topological space, $R^kf_!\cF$ might not vanish for arbitrary large $k$. In order to solve this problem, let us restrict to the full subcategory $\Top_0$ of the category of locally compact Hausdorff topological spaces that have finite cohomological dimension (or, even simpler, to the category of smooth manifolds). In this subcategory we do have $Rf_!: \Sh{X} \to \Db{\Sh{Y}}$, for $f: X \to Y$ continuous, so it induces a map in $K$-theory $f_!: \K{\Sh{X}} \to \K{\Sh{Y}}$ and analogously for $f^*: \K{\Sh{Y}} \to \K{\Sh{X}}$. Even more, $f^*$ is a ring homomorphism and $f_!$ is a module homomorphism over $\K{\Sh{\star}} = \K{(\Vect{\QQ})} = \ZZ$. By the previous properties, $\K{\Sh{-}}$ is a $\Top_0$-algebra with $A: f \mapsto f^*$ and $B: f \mapsto f_!$. Observe that the unit object in $\K{\Sh{X}}$ is the image of constant sheaf $\underline{\QQ}_X$ on $X$ with stalk $\QQ$.

Moreover, given $X \in \Top_0$, let $c_X: X \to \star$ the projection onto the singleton set. Then, the characteristic of $X$ is the object $(c_X)_!(\underline{\QQ}_X)$ which is a sheaf whose unique stalk is $(c_X)_!(\underline{\QQ}_X)_\star = [H_c^\bullet(X; \QQ)] = \chi_c(X) \in \K{(\Vect{\QQ})} = \ZZ$. Hence, the characteristic of the object $X$ is nothing but the Euler characteristic of $X$ (with compact support). This justifies the name `characteristic'.
\end{ex}

\subsection{Quantisation}
\label{sec:quantisation}

Given a category $\cC$ with pullbacks, we can construct the $2$-category of \emph{spans of $\cC$}, $\Span{\cC}$. As described in \cite{Benabou}, the objects of $\Span{\cC}$ are the same as the ones of $\cC$. A morphism $a \to b$ in $\Span{\cC}$ is an span, that is, a triple $(d, f,g)$ of morphisms
$$
\xymatrix{
&\ar[dl]_{f} d \ar[dr]^{g}&\\
a&& b
}
$$
where $d \in \cC$. Given two spans $(d_1, f_1, g_1): a \to b$ and $(d_2, f_2, g_2): b \to c$, we define the composition $(d_2, f_2, g_2) \circ (d_1, f_1,g_1) = (d_1 \times_b d_2, f_1 \circ f_2', g_2 \circ g_1')$, where $f_2', g_1'$ are the morphisms in the pullback diagram
$$
\xymatrix{
&&\ar[dl]_{f'_2}d_1 \times_b d_2\ar[dr]^{g'_1}&&\\
&\ar[dl]_{f_1}d_1\ar[dr]^{g_1}&&\ar[dl]_{f_2}d_2\ar[dr]^{g_2}&\\
a&&b&&c
}
$$
Finally, a $2$-morphism $(d, f,g) \Rightarrow (d', f',g')$ between $a \stackrel{\hspace{3pt}f}{\leftarrow} d \stackrel{g}{\rightarrow} b$ and $a \stackrel{\hspace{5pt}f'}{\leftarrow} d' \stackrel{g'}{\rightarrow} b$ is a morphism $\alpha: d' \to d$ (notice the inverted arrow!) such that the following diagram commutes
\[
\begin{displaystyle}
   \xymatrix
   {	& d' \ar[rd]^{g'} \ar[ld]_{f'} \ar[dd]^\alpha & \\
   		a &  & b\\
	&d  \ar[ru]_{g} \ar[lu]^{f}& 
   }
\end{displaystyle}   
\]
Moreover, if $\cC$ is a monoidal category, $\Span{\cC}$ inherits a monoidal structure by tensor product on objects and morphisms.

\begin{prop}\label{prop-quantisation}
Let $\cA$ be a $\cC$-algebra, with $\cC$ a category with final object $\star$ and pullbacks. Then, there exists a lax monoidal $2$-functor $\Qtm{\cA}: \Span{\cC} \to \Modt{\cA_\star}$ such that
$$
	\Qtm{\cA}(a) = \cA_a \hspace{1.5cm} \Qtm{\cA}(d, f, g) = g_! \circ f^*: \cA_a \to \cA_b,
$$
for $a, b \in \cC$ and a span $a \stackrel{f}{\leftarrow} d \stackrel{g}{\to} b$. The functor $\Qtm{\cA}$ is called the \emph{quantisation} of $\cA$.
\begin{proof}
More detailed, the functor $\Qtm{\cA}: \Span{\cC} \to \Modt{\cA_\star}$ is given as follows:
\begin{itemize}
	\item For any $a \in \cC$ we define $\Qtm{\cA}(a) = \cA_a$ as $\cA_\star$-algebra.
	\item Fixed $a, b \in \cC$, we define the functor
	$$
		(\Qtm{\cA})_{a,b}: \Hom_{\Span{\cC}}(a, b) \to \Hom_{\Modt{\cA_\star}}(\cA_a, \cA_b)$$ by:
	\begin{itemize}
		\item For a $1$-morphism $a \stackrel{f}{\leftarrow} d \stackrel{g}{\rightarrow} b$ we define $\Qtm{\cA}(d,f,g) = g_! \circ f^*: \cA_a \to \cA_d \to \cA_b$.
		\item Let $\alpha: (d, f,g) \Rightarrow (d', f',g')$ be a $2$-morphism given by a morphism $\alpha: d' \to d$ of $\cC$. Since $\alpha$ is a $2$-cell in $\Span{\cC}$ we have that $\alpha^* \circ f^* = f'^*$ and $g_! \circ \alpha_! = g'_!$. This implies that $(\cA_{d}, \cA_{d'}, f^*, g_!, f'^*, g'_!, \alpha^*, \alpha_!)$ is an immediate twist $g_! \circ f^* \Rightarrow g'_! \circ f'^*$.
%Since $\alpha$ is a $2$-cell in $\Span{\cC}$ we have that $g_! \circ \psi \circ f^* = (g')_! \circ (f')^*$ so $\psi$ is a twist $g_! \circ f^* \Rightarrow (g')_! \circ (f')^*$. Observe that, if $\alpha$ is an isomorphism, then $\psi = 1_{\cA_d}$.
%\item For a $2$-morphism $\alpha: (d, f,g) \Rightarrow (d', f',g')$ given by $\alpha: d' \to d$, we define the endomorphism $\psi = \alpha_! \alpha^*: \cA_d \to \cA_{d'}$. Since $\alpha$ is a $2$-cell in $\Span{\cC}$ we have that $g_! \circ \psi \circ f^* = (g')_! \circ (f')^*$ so $\psi$ is a twist $g_! \circ f^* \Rightarrow (g')_! \circ (f')^*$. Observe that, if $\alpha$ is an isomorphism, then $\psi = 1_{\cA_d}$.
	\end{itemize}
	\item For $(\Qtm{\cA})_{1_a}$, $a \in \cC$, we take the identity $2$-cell.
	\item Given $1$-morphisms $(d_1, f_1, g_1): a \to b$ and $(d_2, f_2, g_2): b \to c$ we have $(\Qtm{\cA})_{b,c}(d_2, f_2,g_2) \circ (\Qtm{\cA})_{a,b}(d_1, f_1,g_1) = (g_2)_!(f_2)^*(g_1)_!(f_1)^*$. On the other hand, by definition of composition in $\Span{\cC}$ we have $(\Qtm{\cA})_{a,c}((d_2, f_2,g_2) \circ (d_1, f_1, g_1))  =  (g_2)_! (g_1')_!  (f_2')^*  (f_1)^*$. Here $g_1'$ and $f_2'$ are the maps in the pullback
	$$
\xymatrix{
&&\ar[dl]_{f'_2}d_1\times_b d_2\ar[dr]^{g'_1}&&\\
&\ar[dl]_{f_1}d_1\ar[dr]^{g_1}&&\ar[dl]_{f_2}d_2\ar[dr]^{g_2}&\\
a&&b&&c.
}
$$
By the Beck-Chevalley condition we have $(g_1')_!(f_2')^* = (f_2)^*(g_1)_!$ and the two morphisms agree. Thus, we can take the $2$-cell $(\Qtm{\cA})_{a,b,c}$ as the identity.
\end{itemize}

This proves that $\Qtm{\cA}$ is a $2$-functor. Furthermore, $\Qtm{\cA}$ is also lax monoidal taking $\Delta_{a, b}= \boxtimes: \cA_a \otimes \cA_b \to \cA_{a \times b}$ to be the external product as in Section \ref{sec:C-algebras}. In order to check that, suppose that we have spans $(d, f, g): a \to b$ and $(d', f', g'): a' \to b'$. Then, since both the pullback and the pushout maps preserve the tensor product we have that, for all $z \in \cA_{a}$ and $w \in \cA_{a'}$, $g_!f^*z \boxtimes g'_!f'^*w = (g \times g')_!(f^*z \boxtimes f'^*w) = (g \times g')_!(f \times f')^*(z \boxtimes w)$. Therefore, the following diagram commutes
\[
\begin{displaystyle}
   \xymatrix
   {
   	\cA_a \otimes \cA_{a'} \ar[r]^{\boxtimes} \ar[d]_{g_!f^* \otimes g'_!f'^*} & \cA_{a \times a'} \ar[d]^{(g \times g')_!(f \times f')^*}\\
   	\cA_{b} \otimes \cA_{b'} \ar[r]_{\boxtimes} & \cA_{b \times b'}
   }
\end{displaystyle}   
\]
Hence, $\Delta: \cA_{-} \otimes \cA_{-} \to \cA_{- \times -}$ is natural, as we wanted.
\end{proof}
\end{prop}

\begin{rmk}\label{rmk:non-invertible-2-cells}
\begin{itemize}
	\item If the functor $A$ defining the $\cC$-algebra is monoidal, then $\Qtm{\cA}$ is also strict monoidal.
%	\item If the Beck-Chevalley condition was satisfied up to natural isomorphism, then we would have equality after a pair of automorphisms in $\cA_{d_1}$ and $\cA_{d_2}$ which correspond to an invertible $2$-cell in $\Modt{\cA_\star}$. In this case, $\Qtm{\cA}$ would be a pseudo-functor. 
	\item If $\alpha: (d,f,g) \Rightarrow (d',f',g')$ is an invertible $2$-cell between spans, then $\Qtm{\cA}(d,f,g)= g_! \alpha_!\alpha^*f^* = (g')_!(f')^* = \Qtm{\cA}(d',f',g')$ since $\alpha_!\alpha^* =  (1_{d})^*  (1_{d})_!= 1_{\cA_d}$. However, if $\alpha$ is not invertible, this may no longer hold and $\alpha$ induces a non-trivial twisting. This is the reason of keeping track this cumbersome structure.
\end{itemize}
\end{rmk}

The construction described in this section is strongly related to the sheaf theoretic formalism of \cite[Part III]{Gaitsgory-Rozenblyum}, specially Chapters 7 and 8. In that paper, it is proven that a bivalent functor $\Phi: \cC \to \cS$ satisfying the six functors formalism (which is the analogous of a $\cC$-algebra) is equivalent to a functor out of the $2$-category of correspondences \cite[Theorem 2.13 of Chapter 7]{Gaitsgory-Rozenblyum}.

\subsection{Field theory}
\label{sec:field-theory}

Let us fix a sheaf $\cS$. Let $(W_1 s_1)$ and $(W_2, s_2)$ be two $n$-dimensional objects of $\EEmbc{\cS}$ and let $(M, s)$ be a $(n-1)$-dimensional object. Suppose that there exist tame embeddings $f_1: (M,s) \to (W_1, s_2)$ and $f_2: (M, s) \to (W_2, s_2)$ of boundaries, as in Example \ref{ex:embc-map-boundary}.
In that case, the pushout of $f_1$ and $f_2$ in $\EEmbc{\cS}$, $(W_1 \cup_M W_2, s_1 \cup s_2)$, exists and it is given by the gluing of $W_1$ and $W_2$ along $M$.

\begin{figure}[h]
	\begin{center}
	\includegraphics[scale=0.35]{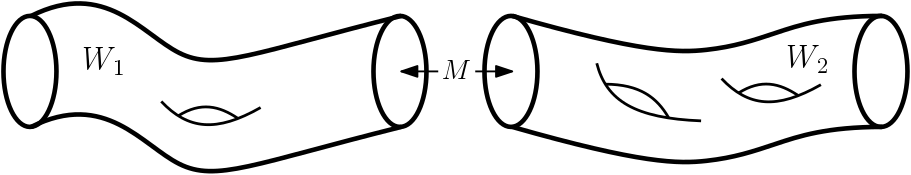}
	\end{center}
\vspace{-0.1cm}
\end{figure}

This situation will be called a \emph{gluing pushout}.

\begin{defn}\label{def:seifert-vankampen}
Given a category $\cC$ with final object, a contravariant functor $F: \EEmbc{\cS} \to \cC$ is said to have the \emph{Seifert-van Kampen property} if $F$ sends gluing pushouts in $\EEmbc{\cS}$ into pullbacks of $\cC$ and sends the initial object of $\EEmbc{\cS}$ (i.e.\ $\emptyset$) into the final object of $\cC$.
\end{defn}

\begin{prop}\label{prop:field-theory}
Let $\cS$ be a sheaf, let $\cC$ be a category with final object and pullbacks and let $F: \EEmbc{\cS} \to \cC$ be a contravariant functor satisfying the Seifert-van Kampen property. Then, there exists a monoidal $2$-functor $\cF_F: \EBord{n}{\cS} \to \Span{\cC}$ such that
$$
	\cF_F(M, s) = F(M, s), \hspace{1cm} F(M_1, s_1) \stackrel{F(i_1)}{\leftarrow} F(W, s) \stackrel{F(i_2)}{\rightarrow} F(M_2, s_2),
$$
for all objects $(M,s), (M_1, s_1), (M_2, s_2) \in \EBord{n}{\cS}$ and bordisms $(W, s): (M_1, s_1) \to (M_2, s_2)$ where $i_k: M_k \to W$ are the inclusions as boundaries. In this situation, the functor $F$ is called the \emph{geometrisation} and $\cF_F$ is called the \emph{field theory} of $F$.
\begin{proof}
The complete definition of $\cF$ is given as follows. Given $(M,s) \in \EBord{n}{\cS}$, we just define $\cF_F(M,s) = F(M,s)$. With respect to morphisms, given a bordism $(W,s): (M_1,s_1) \to (M_2, s_2)$ in $\EBord{n}{\cS}$, let $i_1: M_1 \to W$ and $i_2: M_2 \to W$ be the inclusions of $M_1$ and $M_2$ as boundaries of $W$. The functor $\cF_F$ assigns the span
\[
\begin{displaystyle}
   \xymatrix
   {
   	F(M_1, s_1) && F(W, s) \ar[rr]^{F(i_2)}\ar[ll]_{\;\;\;\;\;\;\;F(i_1)} && F(M_2, s_2) 
   }
\end{displaystyle}   
\]

This assignment is a functor. In order to check it, let $(W,s): (M_1, s_1) \to (M_2, s_2)$ and $(W', s'): (M_2, s_2) \to (M_3, s_3)$ be two bordisms with inclusions $i_k: M_k \to W$ and $i_k': M_k \to W'$. By construction $(W' \circ W, s \cup s')$ is the gluing pushout in $\EEmbc{\cS}$
\[
\begin{displaystyle}
   \xymatrix
   {
  	(M_2, s_2) \ar[d]_{i_2'}\ar[r]^{i_2} & (W,s) \ar[d]^{j} \\
  	(W',s') \ar[r]_{j'\;\;\;\;\;\;\;}  & (W' \circ W, s \cup s')
   }
\end{displaystyle}   
\]
Since $F$ sends gluing pushouts into pullbacks, the following diagram is a pullback in $\cC$
\[
\begin{displaystyle}
   \xymatrix
   {
  	F(W' \circ W, s \cup s') \ar[r]^{\;\;\;\;\;\;\;\;\;F(j)} \ar[d]_{F(j')} & F(W,s) \ar[d]^{F(i_2)} \\
  	F(W',s') \ar[r]_{F(i_2')} & F(M_2,s_2)
   }
\end{displaystyle}   
\]
Therefore, $\cF_F(W', s') \circ \cF_F(W, s)$ is given by the span
	$$
\xymatrix{
&&\ar[dl]_{F(j)} F(W' \circ W, s \cup s') \ar[dr]^{F(j')}&&\\
&\ar[dl]_{F(i_1)}F(W, s)\ar[dr]^{F(i_2)}&&\ar[dl]_{F(i_2')}F(W',s')\ar[dr]^{F(i_3')}&\\
F(M_1, s_1)&&F(M_2, s_2)&&F(M_3, s_3).
}
$$
Since $i_1 \circ j$ and $i_3' \circ j'$ are the inclusions onto $W' \circ W$ of $M_1$ and $M_3$ respectively, the previous span is also $\cF_F(W' \circ W)$, as we wanted.

For the monoidality, let $(M_1, s_1), (M_2, s_2) \in \EBord{n}{\cS}$. As the coproduct $(M_1 \sqcup M_2, s_1 \sqcup s_2)$ can be seen as a gluing pushout along $\emptyset$, $F$ also sends coproducts in $\EBord{n}{\cS}$ into products on $\cC$. Hence, $F(M_1 \sqcup M_2, s_1 \sqcup s_2) = F(M_1, s_1) \times F(M_2, s_2)$ and, since the monoidal structure on $\Span{\cC}$ is given by products on $\cC$, monoidality holds for objects. For morphisms, the argument is analogous.

For the $2$-functor structure, suppose that $(f, \alpha)$ is a $2$-cell between $1$-morphisms $(W, s), (W', s'): (M_1, s_1) \to (M_2, s_2)$ of $\EBord{n}{\cS}$. In that case $(f, \alpha)$ is also a morphism in $\EEmbc{\cS}$ so we obtain a morphism $F(f, \alpha): (W', s') \to (W, s)$ fitting in the commutative diagram
\[
\begin{displaystyle}
   \xymatrix
   {
   	& F(W', s') \ar[rd]^{F(i_2')} \ar[dd]_{F(f,\alpha)} \ar[ld]_{F(i_1)} & \\
	F(M_1, s_1) & & F(M_2, s_2) \\
	& F(W, s) \ar[ru]_{F(i_2)} \ar[lu]^{F(i_1)} &
   }
\end{displaystyle}
\]
This produces the desired $2$-morphism in $\Span{\cC}$.
\end{proof}
\end{prop}

\begin{rmk}
We thank D. Ben-Zvi for suggesting us the name `field theory' based on the physical interpretation of TQFTs.
\end{rmk}

Composing the functors from Proposition \ref{prop-quantisation} and Proposition \ref{prop:field-theory}, we obtain the main result of this section.

\begin{thm}\label{thm:physical-constr-sheaf}
Let $\cC$ be a category with final object $\star$ and pullbacks. Given a functor $F: \EEmbc{\cS} \to \cC$ with Seifert-van Kampen property and a $\cC$-algebra $\cA$, there exists a lax monoidal Topological Quantum Field Theory over $\cS$
$$
	Z_{F, \cA}: \EBord{n}{\cS} \to \Modt{\cA_\star}.
$$
\end{thm}

\begin{rmk}
\label{rmk:expression-TQFT}
From the explicit construction given in the proofs of Propositions \ref{prop-quantisation} and \ref{prop:field-theory}, the functor $Z_{F, \cA}$ satisfies:
\begin{itemize}
	\item For an object $(M, s) \in \EBord{n}{\cS}$, it assigns $Z_{F, \cA}(M,s) = \cA_{F(M,s)}$.
	\item For a bordism $(W,s): (M_1, s_1) \to (M_2, s_2)$, it assigns $Z_{F, \cA}(W,s) = F(i_2)_! \circ F(i_1)^*: \cA_{F(M_1,s_1)} \to \cA_{F(M_2, s_2)}$.
	\item For a closed $n$-dimensional manifold $(W,s)$, seen as a bordism $(W, s): \emptyset \to \emptyset$, the homomorphism $Z_{F, \cA}(W, s): \cA_\star \to \cA_\star$ is given by multiplication by the characteristic $\mu({F(W,s)}) \in \cA_\star$. This follows from the fact that, since $F(\emptyset)=\star$, the inclusion $i: \emptyset \to (W,s)$ gives the projection $c = F(i): F(W,s) \to \star$. Hence, we have $Z_{F, \cA}(W)(1_\star) = c_! \circ c^*(1_\star) = c_!(1_{F(W,s)}) = \mu(F(W,s))$. Here, $1_\star$ and $1_{F(W,s)}$ denote the units in $\cA_\star$ and $\cA_{F(W,s)}$ respectively.
\end{itemize}
\end{rmk}

\begin{rmk}In the context of bordisms with sheaves, the requirement of the whimsical twisting structure on $\Modt{\cA_\star}$ becomes evident. The existence of a non-invertible $2$-cell $(f,\alpha): (W, s) \Rightarrow (W, s')$ reflects a non-invertible morphism $\alpha: s \to \cS(f)(s')$ that can be interpreted as a restriction of the extra structure sheaf. As explained in Remark \ref{rmk:non-invertible-2-cells}, these non-invertible cells become non-trivial twists in $\Modt{\cA_\star}$ that compare how the homomorphism changes under morphisms of sheaves.
\end{rmk}

\begin{rmk}\label{rmk:field-theory-very-lax}

We can also consider the case in which the geometrisation functor $F: \EEmbc{\cS} \to \cC$ no longer has the Seifert-van Kampen property, but it still maps the initial object into the final object. In that case the image of a gluing pushout under $F$ is not a pullback but, as for any functor, it is a cone. Suppose that we have two bordisms $(W,s): (M_1, s_1) \to (M_2, s_2)$ and $(W', s'): (M_2, s_2) \to (M_3, s_3)$ that fit in the gluing pushout in $\EEmbc{\cS}$
\[
\begin{displaystyle}
   \xymatrix
   {
   		(M_2, s_2) \ar[r]\ar[d] & (W,s)\ar[d] \\
   		(W', s') \ar[r] & (W' \circ W, s \cup s')
   }
\end{displaystyle}   
\]
By definition, there exists an unique morphism $\phi: F(W' \circ W, s \cup s') \to F(W', s') \times_{F(M_2, s_2)} F(W,s)$ in $\cC$ such that the induced diagram commutes
\[
\begin{displaystyle}
   \xymatrix
   {
   		F(W' \circ W, s \cup s') \ar@{--{>}}[rd]^{\phi} \ar@/^1pc/[rrd] \ar@/_1pc/[ddr] & & \\
   		&F(W',s') \times_{F(M_2, s_2)} F(W, s)  \ar[r]\ar[d]  & F(W', s') \ar[d] \\
   		&F(W,s)  \ar[r] & F(M_2, s_2)
   }
\end{displaystyle}  
\]
This morphism $\phi$ induces a $2$-morphism $\cF_F(W', s') \circ \cF_F(W,s) \Rightarrow \cF_F(W' \circ W, s \cup s')$. Therefore, in this case, $\cF: \EBord{n}{\cS} \to \Span{\cC}$ is no longer a functor but a lax $2$-functor. Thus, the induced functor
$$
	Z_{F, \cA}: \EBord{n}{\cS} \to \Modt{\cA_{\star}}
$$
is a lax monoidal symmetric lax $2$-functor. We will call such a functors \emph{very lax Topological Quantum Field Theories}.
\end{rmk}

As a final remark, it is customary in the literature to focus on the field theory as a functor $\Fld{}: \Bord{n} \to \Span{\cC}$ and on the quantisation as a functor $\Qtm{}: \Span{\cC} \to \Mod{R}$, and to forget about the geometrisation and the $\cC$-algebra, despite that they underlie the whole construction (see, for example, \cite{Freed-Hopkins-Lurie-Teleman:2010, Freed-Hopkins-Teleman:2010}). In this form, the Seifert-van Kampen property and the base change condition are hidden, encoded in the functoriality of $\Fld{}$ of $\Qtm{}$ respectively.

However, in this paper we want to emphasize the role of the geometrisation $F$ and the $\cC$-algebra $\cA$ in the construction. The reason is we are constructing TQFTs with a view towards the creation of new effective computational methods of algebraic invariants. In this way, $\cA$ determines the algebraic invariant under study and $F$ determines the object for which we are going to compute the invariant. This principle opens the door to the development of new computational methods based on TQFTs further than the scope of this paper.

\subsection{Some useful sheaf structures}
\label{sec:some-extra-struct}

In this section, we will describe in detail some examples of extra structures induced by sheaves. The first example is
%the trivial sheaf $\cS_{t}: \Embc \to \Cat$ that assigns, to every compact manifold, the singleton category (i.e.\ the category with an unique object and only the identity as morphism) and, to every morphism the identity morphism. In that case, $\EBord{n}{\cS_t}$ is the usual category of unoriented bordisms. 
the sheaf $\cSp: \Embc \to \Cat$ of unordered configurations of points, also called the sheaf of pairs. For a compact manifold $M$, the category $\cSp(M)$ has, as objects, finite subsets $A \subseteq M$ meeting every connected component and every boundary component of $M$. Given two finite subsets $A_1, A_2 \subseteq M$, a morphism $A_1 \to A_2$ in $\cSp(M)$ is an inclusion $A_1 \subseteq A_2$. With respect to morphisms, if $f: M \to N$ is an embedding the functor $\cSp(f): \cSp(N) \to \cSp(M)$ is given as follows. For an object $A \subseteq N$ it assigns $\cSp(f)(A) = f^{-1}(A) \in \cSp(M)$ and, if we have an inclusion $A_1 \subseteq A_2$, then it gives the inclusion $f^{-1}(A_1) \subseteq f^{-1}(A_2)$ as a morphism in $\cSp(M)$. It is straighforward to check that $\cSp$ is a sheaf.

The associated categories $\EEmbc{\cSp}$ and $\EBord{n}{\cSp}$ will be called the \emph{category of embeddings of pairs} and the \emph{category of $n$-bordisms of pairs} and will be denoted by $\Embpc$ and $\Bordp{n}$, respectively. A (lax monoidal) $\cSp$-TQFT will be referred to as a (lax monoidal) \emph{Topological Quantum Field Theory of pairs}.

Another important sheaf for applications is the so-called sheaf of parabolic structures. The starting point is a fixed set $\Lambda$ that we will call the \emph{parabolic data}. We define $\cSpar{\Lambda}: \Embc \to \Cat$ as the following functor. For a compact manifold $W$, $\cSpar{\Lambda}(W)$ is the category whose objects are (maybe empty) finite sets $Q = \left\{(S_1, \lambda_1), \ldots, (S_r, \lambda_r)\right\}$, with $\lambda_i \in \Lambda$, called \emph{parabolic structures} on $W$. The $S_i$ are pairwise disjoint compact submanifolds of $W$ of codimension $2$ with a co-orientation such that $S_i \cap \partial M = \partial S_i$ transversally. A morphism $Q \to Q'$ between two parabolic structures in $\cSpar{\Lambda}(W)$ is just an inclusion $Q \subseteq Q'$.

Moreover, suppose that we have a tame embedding $f: W_1 \to W_2$ in $\Embc$ and let $Q \in \cSpar{\Lambda}(W_2)$. Given $(S, \lambda) \in Q$, if $S \cap f(\partial W_1)$ transversally then the intersection has the expected dimension and, thus, $f^{-1}(S)$ is a codimension $2$ submanifold of $W_1$. Furthermore, the co-orientation of $S$ induces a co-orientation on $f^{-1}(S)$ by pullback. Hence, we can define $\cSpar{\Lambda}(f)(Q)$ as the set of pairs $(f^{-1}(S), \lambda)$ for $(S, \lambda) \in Q$ with $S \cap f(\partial W_1)$ transversal. For short, we will denote $Q|_{W_1} = \cSpar{\Lambda}(f)(Q)$. Obviously, if $Q \subseteq Q'$ then $Q|_{W_1} \subseteq Q'|_{W_1}$ so $\cSpar{\Lambda}(f): \cSpar{\Lambda}(W_2) \to \cSpar{\Lambda}(W_1)$ is a functor. With this definition, $\cSpar{\Lambda}$ is a sheaf, called the \emph{sheaf of parabolic structures over $\Lambda$}.

\begin{ex}\label{ex:parab-struct}In the case of surfaces, a parabolic structure over a surface $M$ is a set $Q = \left\{(p_1, \lambda_1), \ldots, (p_r, \lambda_r)\right\}$ with $\lambda_i \in \Lambda$ and $p_i \in M$ points with a preferred orientation of a small disc around them.
\end{ex}

As for the sheaf of unordered points, the sheaf $\cSpar{\Lambda}$ gives us categories $\EEmbc{\cSpar{\Lambda}}$ and $\EBord{n}{\cSpar{\Lambda}}$, that we will shorten $\Embparc{\Lambda}$ and $\Bordpar{n}{\Lambda}$. Even more, we can combine the two previous sheaves and to consider the sheaf $\cS_{p, \Lambda} = (\cSpar{\Lambda} \times \cSp) \circ \Delta: \Embc \to \Cat$, where $\Delta: \Embc \to \Embc \times \Embc$ is the diagonal functor. In that case, we will denote by $\Embpparc{\Lambda}$ and by $\Bordppar{n}{\Lambda}$ the categories of $\cS_{p, \Lambda}$-embeddings and $\cS_{p, \Lambda}$-bordisms, respectively. In this case, a (lax monoidal) $\cS_{p, \Lambda}$-TQFT will be called a (lax monoidal) \emph{parabolic Topological Quantum Field Theory of pairs}.

\subsection{Reduction of a TQFT}
\label{sec:reduction-TQFT}

Let $\cC$ be a category with final object and pullbacks, let $F: \Embc \to \cC$ be a contravariant functor with the Seifert-van Kampen property and let $\cA$ be a $\cC$-algebra. By Theorem \ref{thm:physical-constr-sheaf}, these data give rise to a lax monoidal TQFT $Z = Z_{F, \cA}: \Bord{n} \to \Mod{\cA_\star}$. However, for $M \in \Bord{n}$, the module $Z(M) = \cA_{F(M)}$ may be very complicated as it might be infinitely generated. 

This problem can be mitigated if the objects $F(M)$ have some kind of symmetry that may be exploited. For example, suppose that there is a natural action of a group $G$ on $F(M)$ so that the categorical quotient object $F'(M) = F(M)/G \in \cC$ can be defined. In the general case, such a symmetry can be modelled by an assignment $\tau$ that, for any $M \in \Bord{n}$, gives a pair $(F'(M), \tau_{M})$ where $F'(M) \in \cC$ and $\tau_{M}: F(M) \to F'(M)$ is a morphism in $\cC$. Suppose also that $F'(\emptyset) = F(\emptyset)$.

In that case, we can `reduce' the field theory and to consider $\Fld{F, \tau}: \Bord{n} \to \Span{\cC}$ such that $\Fld{F, \tau}(M)=F'(M)$ for $M \in \Bord{n}$. For a bordism $W: M_1 \to M_2$, $\Fld{F, \tau}(W)$ is the span
\[
\begin{displaystyle}
   \xymatrix
   {	
	F'(M_1) && F(W) \ar[rr]^{\tau_{M_2} \circ F(i_2)} \ar[ll]_{\tau_{M_1} \circ F(i_1)} && F'(M_2)
   }
\end{displaystyle}.
\]
With this field theory, we form $Z_\tau = \Qtm{\cA} \circ \Fld{F}'$, called the \emph{prereduction} of $Z$ by $\tau$.

However, even if $F$ had the Seifert-van Kampen property, $\Fld{F, \tau}$ may not be a functor so $Z_\tau$ will be in general a very lax TQFT. Indeed, the prereduction $Z_\tau$ is simpler than $Z$ but, as it is not a functor, we can no longer use it for a computational method. In this section we will show that, under some mild conditions, we can slightly modify $Z_\tau$ in order to obtain an almost-TQFT \cite{GPLM-2017}, $\cZ_\tau$, called the \emph{reduction}. The reduction has essentially the same complexity as $Z_\tau$ so it can be used to give an effective computational method.

In order to do so, consider the wide subcategory $\ETub{n}{\cS}$ of $\EBord{n}{\cS}$ of \emph{$\cS$-tubes}. A morphism $(M,s) \in \EBord{n}{\cS}$ is an object of $\ETub{n}{\cS}$ if there exists a compact $n$-dimensional manifold $W$ such that $\partial W = M$. A bordism $(W,s): (M_1, s_1) \to (M_2, s_2)$ of $\EBord{n}{\cS}$ is called a \emph{strict tube} if $M_1$ and $M_2$ are connected or empty. In this way, the morphisms of $\ETub{n}{\cS}$ are disjoint unions of strict tubes, the so-called tubes. A monoidal functor $Z: \ETub{n}{\cS} \to \Mod{R}$ is called an \emph{almost-TQFT over $\cS$}.

Now, let $Z_\tau: \EBord{n}{\cS} \to \Mod{\cA_\star}$ be a prereduction by $\tau$. For any $(M, s) \in \ETub{n}{\cS}$ connected or empty, we will denote by $\cV_{(M,s)} \subseteq \cA_{F'(M, s)}$ the submodule generated by the images $Z_\tau(W_r, s_r) \circ \ldots \circ Z_\tau(W_1, s_1)(1) \in \cA_{F(M, s)}$ of all the sequences of strict tubes $(W_k, s_k): (M_{k-1}, s_{k-1}) \to (M_{k}, s_k)$ with $M_0 = \emptyset$. Notice that, by definition, the submodules $\cV_{(M,s)}$ are invariant for strict tubes.

\begin{lem}\label{lem:existence-descend-reduction}
Let $\cC$ be a category, let $a, a', b, b' \in \cC$ and consider morphisms $f_!: a \to a'$, $f^*: a' \to a$, $g: b \to b'$ and $h: a \to b$. Suppose that the morphism $\eta = f_! \circ f^*: a' \to a'$ is invertible, then there exists a unique morphism $h': a' \to b'$ such that $h' \circ f_! = g \circ h$.
\[
\begin{displaystyle}
   \xymatrix
   {	
	a \ar[r]^h \ar[d]_{f_!} & b \ar[d]^g \\
	a' \ar@{--{>}}[r]_{h'} \ar@/_1pc/[u]_{f^*} & b'
   }
\end{displaystyle}   
\]
\begin{proof}
First, let us prove that $h'$ is unique. Pre-composing with $f^* \circ \eta^{-1}: a' \to a$ we obtain that
$$
	g \circ h \circ f^* \circ \eta^{-1} = h' \circ f_! \circ f^* \circ \eta^{-1} = h' \circ \eta \circ \eta^{-1} = h'.
$$
Actually, this calculation shows the existence of that morphism since we must take $h' = g \circ h \circ f^* \circ \eta^{-1}$, and it is a straightforward check that $h'$ has the desired property.
\end{proof}
\end{lem}

\begin{prop}\label{prop:existence-reduction-TQFT}
Let $Z: \EBord{n}{\cS} \to \Mod{\cA_\star}$ be a lax monoidal TQFT and let $\tau$ be a reduction. Suppose that, for all $(M,s) \in \EBord{n}{\cS}$, we have that $(\tau_{(M,s)})_! \circ (\tau_{(M,s)})^*(\cV_{(M,s)}) \subseteq \cV_{(M,s)}$ and the morphisms $\eta_{(M,s)} = (\tau_{(M,s)})_! \circ (\tau_{(M,s)})^*: \cV_{(M,s)} \to \cV_{(M,s)}$ are invertible. Then, there exists an almost-TQFT
$$
	\cZ_\tau: \ETub{n}{\cS} \to \Mod{\cA_\star}
$$
such that $\cZ_\tau(M, s) = \cV_{(M,s)}$ for all $(M, s) \in \ETub{n}{\cS}$ connected or empty and $\cZ_\tau(W, s) \circ (\tau_{(M_1,s_1)})_! = (\tau_{(M_2, s_2)})_! \circ Z(W, s)$ for all strict tubes $(W, s): (M_1, s_1) \to (M_2, s_2)$. This TQFT is called the \emph{reduction} of $Z$ via $\tau$.
\begin{proof}
Recall that, in order to define an almost-TQFT, it is enough to define it on strict tubes and to extend it to a general tube by tensor product \cite{GPLM-2017}. Therefore, for $(M, s) \in \ETub{n}{\cS}$ connected or empty, we assign $\cZ_\tau(M,s) = \cV_{(M,s)}$. For a strict tube $(W, s): (M_1, s_1) \to (M_2, s_2)$, by Lemma \ref{lem:existence-descend-reduction} there exists a unique morphism $\cZ_\tau(W, s): \cV_{(M_1, s_1)} \to \cV_{(M_2, s_2)}$ such that the following diagram commute
\[
\begin{displaystyle}
   \xymatrix
   {	
		\tau_{(M_1, s_1)}^*\cV_{(M_1, s_1)} \ar[rr]^{Z(W, s)\;\;\;\;\;\;\;\;\;}\ar[d]_{(\tau_{(M_1, s_1)})_!}  && Z(W, s)\left(\tau_{(M_1, s_1)}^*{\cV_{(M_1, s_1)}}\right) \ar[d]^{(\tau_N)_!}\\
		\cV_{(M_1, s_1)} \ar@{--{>}}[rr]_{\cZ_\tau(W, s)} && \cV_{(M_2, s_2)}
   }
\end{displaystyle}   
\]

In order to prove that $\cZ_\tau$ is a functor, suppose that $(W, s): (M_1, s_1) \to (M_2, s_2)$ and $(W', s'): (M_2, s_2) \to (M_3, s_3)$ are strict tubes. Then, by the previous proposition, we have a commutative diagram
\[
\begin{displaystyle}
   \xymatrix
   {	
		\tau_{(M_1, s_1)}^*\cV_{(M_1, s_1)} \ar[r]^{Z(W, s)\;\;\;\;\;\;\;\;\;}\ar[d]_{\left(\tau_{(M_1, s_1)}\right)_!}  & Z(W, s)\left(\tau_{(M_1, s_1)}^*\cV_{(M_1, s_1)}\right) \ar[d]^{\left(\tau_{(M_2, s_2)}\right)_!} \ar[r]^{Z(W', s')\;\;\;\;\;\;\;\;\;}& Z(W' \circ W, s \cup s')\left(\tau_{(M_1, s_1)}^*\cV_{(M_1, s_1)}\right) \ar[d]^{\left(\tau_{(M_3, s_3)}\right)_!}\\
		\cV_{(M_1, s_1)} \ar[r]_{\cZ_\tau(W,s)} & \cV_{(M_2, s_2)} \ar[r]_{\cZ_\tau(W',s')} & \cV_{(M_3, s_3)}
   }
\end{displaystyle}   
\]
Therefore, we have $\cZ_\tau(W', s') \circ \cZ_\tau(W, s) \circ (\tau_{(M_1,s_1)})_! = (\tau_{(M_3, s_3)})_! \circ Z(W', s') \circ Z(W, s) = (\tau_{(M_3, s_3)})_! \circ Z(W' \circ W, s \cup s')$ and, by uniqueness, this implies that $\cZ_\tau(W', s') \circ \cZ_\tau(W, s) = \cZ_\tau(W' \circ W, s \cup s')$, as we wanted to prove.
\end{proof}
\end{prop}

\begin{rmk}\label{rmk:reduction-TQFT}
\begin{itemize}
	\item The almost-TQFT, $\cZ_\tau: \ETub{n}{\cS} \to \Mod{\cA_\star}$, satisfies that, for all closed $n$-dimensional manifolds $W$ and $s \in \cS(W)$
$$
	\cZ_\tau(W,s)(1) = \cZ_\tau(W,s) \circ (\tau_\emptyset)_!(1) = (\tau_\emptyset)_! \circ Z(W,s)(1) = Z(W,s)(1),
$$
where we have used that $\tau_\emptyset = 1_\star$. Hence, $\cZ_\tau$ and $Z$ compute the same invariant.
	\item It may happen, and it will be the case for representation varieties, that $\eta_{(M,s)}: \cV_{(M,s)} \to \cV_{(M,s)}$ are not invertible as $\cA_\star$-module. However, it could happen that, for some fixed multiplicative system $S \subseteq \cA_\star$, all the extensions to the localizations $\eta_{(M,s)}: S^{-1}\cV_{(M,s)} \to S^{-1}\cV_{(M,s)}$ are invertible. In that case, we can fix the problem by localizing all the modules and morphisms of the original TQFT to obtain another TQFT, $Z: \Bord{n} \to \Mod{S^{-1}\cA_\star}$, to which we can apply the previous construction.
%A tipical choice for $S \subseteq \cA_\star$ is to take the multiplicative system of non-zero divisors of $\cA_\star$ (so $S$ is the largest system with $\cA_\star \to S^{-1}\cA_\star$ a monomorphism). In that case, for any $\cA_\star$-module $N$, it is customary to denote $S^{-1}N=N_0$ and to call it the \emph{total quotient} of $N$. Other option is to take $S = \cA_\star - \left\{0\right\}$, but, in that case, $N \to S^{-1}N$ may collapse for zero divisors.
\end{itemize}
\end{rmk}

\begin{ex}\label{ex:reduction-TQFT-surfaces}
For $n=2$ (surfaces) and no sheaf, the unique non-empty connected object of $\Tub{2}$ is $S^1$. Hence, it is enough to consider $\cV = \cV_{S^1}$ and $\tau = \tau_{S^1}: F(S^1) \to F'(S^1)$. Actually, $\cV$ is the submodule generated by the elements $Z_\tau(L)^g \circ Z_\tau(D)(1)$ for $g \geq 0$, where $L: S^1 \to S^1$ is the holed torus and $D: \emptyset \to S^1$ is the disk. That is, it is generated by the images of all the compact surfaces with a disc removed. In that case, the only condition we need to check is that $\eta = \tau_! \circ \tau^*: \cV \to \cV$ is invertible.
\end{ex}

\begin{rmk}\label{rmk:left-reduction}
Instead of the diagram of Lemma \ref{lem:existence-descend-reduction}, we could look for a morphism $\cZ_\tau'(W, s): \cV_{(M_1,s_1)} \to \cV_{(M_2,s_2)}$ such the diagram
	\[
\begin{displaystyle}
   \xymatrix
   {	
		\tau_{{(M_1,s_1)}}^*\cV_{{(M_1,s_1)}} \ar[r]^{Z(W,s)\;\;\;\;\;\;\;\;\;}  & Z(W,s)\left(\tau_{{(M_1,s_1)}}^*{\cV_{{(M_1,s_1)}}}\right) \\
		\cV_{(M_1,s_1)} \ar[u]^{\tau_{(M_1,s_1)}^*} \ar@{--{>}}[r]_{\cZ_\tau'(W,s)} & \cV_{(M_2,s_2)} \ar[u]_{\tau_{(M_2,s_2)}^*}
   }
\end{displaystyle}   
\]
commutes. In that case, in the same conditions as in Proposition \ref{prop:existence-reduction-TQFT}, $\cZ_\tau'$ exists, it is an almost TQFT and $\cZ_\tau'(W,s)=\eta^{-1}_{(M_2,s_2)} \circ Z_\tau(W,s)$. We will call $\cZ_\tau'$ the \emph{left $\tau$-reduction}.
\end{rmk}

\section{Grothendieck rings of algebraic varieties}
\label{sec:grothendieck-rings}

Along this section, we will work over a fixed field $k$. Let $X$ be an algebraic variety and consider the slice category $\Varrel{X}$ of varieties over $X$. An object of this category is a pair $(Y, f)$, where $Y$ is an algebraic variety and $f$ is a regular morphism $f: Y \to X$. Its morphisms are intertwining regular maps between them. Moreover, given a regular morphism $h: X \to X'$ we have functors $h_!: \Varrel{X} \to \Varrel{X'}$ and $h^*: \Varrel{X'} \to \Varrel{X}$. They are given by post-composition, $h_!(Y, f)=(Y, h\circ f)$ for $(Y,f) \in \Varrel{X}$, and pullback, $h^*(Z, g)=(Z \times_{X'} X, g')$ for $(Z,g) \in \Varrel{X'}$, where $g'$ is the regular map fitting in the diagram

\[
\begin{displaystyle}
   \xymatrix
   {
	Z \times_{X'} X \ar[r] \ar[d]_{g'} & Z \ar[d]^{g}\\
	X \ar[r]_{h} & X'
   }
\end{displaystyle}
\]

In particular, if we take $X = \star$ we have that $\Varrel{\star} = \Var{k}$, the usual category of $k$-algebraic varieties. The functors induced by the final morphism $c: X \to \star$ are $c_!: \Varrel{X} \to \Var{k}$, the map that forgets the accompanying morphism, and $c^*: \Var{k} \to \Varrel{X}$, the usual cartesian product.

\subsection{Grothendieck rings as $\Var{k}$-algebras}
\label{sec:grothendieck-rings-Var-algebras}

The category $\Varrel{X}$ can be endowed with the structure of a semi-ring with the disjoint union and the fibered product. In this context, for a regular morphism $h: X \to X'$, the functor $h_!$ preserves the disjoint union and $h^*$ preserves both the disjoint union and the fibered product. We can form the associated Grothendieck ring (also known as the $K$-theory ring), $\KVarrel{X}$. It is a commutative unitary ring with unit denoted $T_X = (X, \Id_X) \in \KVarrel{X}$.

Furthermore, $\KVarrel{X}$ inherits a natural module structure over the ring $\KVar{k} = \KVarrel{\star}$ given by cartesian product (i.e. fibered product over the final object). In this way, we have that $h^*: \KVarrel{X'} \to \KVarrel{X}$ is a ring homomorphism and $h_!: \KVarrel{X} \to \KVarrel{X'}$ is a $\KVar{k}$-module homomorphism.

These data allow us to construct a $\Var{k}$-algebra, denoted $\KVarrel{}$, by assigning to an object $X \in \Var{k}$, the ring $\KVarrel{X}$. The induced morphisms are the homomorphisms $h_!, h^*$ described above. The Beck-Chevalley condition is a consequence of the usual base change formula for regular morphisms.

Therefore, by Proposition \ref{prop-quantisation}, this $\Var{k}$-algebra gives rise to a lax monoidal functor $\Qtm{} = \Qtm{\KVarrel{}}: \Span{\Var{k}} \to \Mod{\K{\Var{k}}}$. By construction, to a span of the form $S: \star \leftarrow X \rightarrow \star$, it assigns the homomorphism $\Qtm{}(S): \K{\Var{k}} \to \K{\Var{k}}$ such that $\Qtm{}(S)(\star) = [X]$, where $[X]$ denotes the virtual class of $X$ (i.e. the image of $X$ in $\KVar{k}$). 

\begin{rmk}
\begin{itemize}
	\item By definition, if $f: X \to B$ is a regular morphism, we have that $[X, f] = f_!(T_X) \in \KVarrel{B}$, where $T_X \in \KVarrel{X}$ is the unit of the ring. In order to get in touch with the notation of \cite{LMN}, and if we want to emphasize the base space, such element will be denoted $\RMc{f}{X}{B}=[X,f]$, or just $\RM{X}{B}$ if the morphism $f$ is clear from the context. As we will show in the computations of Sections \ref{sec:Hodge-monodromy-covering-spaces} and \ref{sec:genus-tube}, it plays the role of the Hodge monodromy representation introduced in \cite{LMN}. For that reason, we will also call $\RM{X}{B}$ the Hodge monodromy of $X$ over $B$.
	\item In particular, if $B = \star$, then $\RM{B}{\star}=[B] \in \KVar{k}$. In this sense, the Hodge monodromy generalizes the usual virtual classes in the Grothendieck ring.
	\item Given a regular morphism $X \to B$, since the $\KVar{k}$-module structure on $B$ is given by cartesian product, then $\RM{F \times X}{B}=[F]\RM{X}{B}$. Moreover, if $Y \to X$ is a regular morphism with fiber $F$ that is a locally trivial fibration in the Zariski, we can decompose $Y$ and $X$ into their trivializing pieces to still get $\RM{Y}{B} = [F]\RM{X}{B}$. This multiplicative property will be widely used in the computations of Section \ref{sec:sl2-repr-var} and is the analogous of Proposition 2.4 of \cite{LMN} in the context of the Grothendieck ring.
\end{itemize}
\end{rmk}

In order to get in touch with the definition of Hodge monodromy representation of \cite{LMN}, we need to consider Saito's theory of mixed Hodge modules \cite{Saito:1990}. As described in \cite{Gonzalez-Prieto:Thesis}, for the case $k=\CC$ and $B \in \Var{\CC}$, we can substitute ring $\KVarrel{B}$ by the Grothendieck ring of the category of mixed Hodge modules over $B$, $\KM{B}$. For a regular morphism $h: B \to B'$, there are also well-defined operations of pushout, $h_!: \KM{B} \to \KM{B'}$, and pullback, $h^*: \KM{B'} \to \KM{B}$.

This implies that $\KM{}$ defines a $\Var{\CC}$-algebra with underlying ring $\KM{\star} = \K{\MHSq}$, the Grothendieck ring of rational (mixed) Hodge structures. For a detailed construction, see \cite{Gonzalez-Prieto:Thesis}. In this way, the $\Var{\CC}$-algebra of virtual varieties, $\KVarrel{B}$, can be seen as a promotion of $\KM{}$ to compute, not only Hodge structure, but the whole virtual class. Further relations between these two constructions are analyzed in \cite{GP-2019}.

In particular, if we denote by $q = \QQ(-1)$ the $(-1)$-Tate structure of weight $2$, then $\K{\MHM{B}}$ has a natural $\ZZ[q,q^{-1}]$-module structure inherited from the one as $\KVar{\CC}$-module. Now, suppose that we have a representation of mixed Hodge modules $\rho: \pi_1(X) \to \Aut(V)$, where $V$ is a mixed Hodge structure of balanced type i.e.\ $V^{p,q}=0$ if $p \neq q$. In that case, $\rho$ naturally induces an element of $\KM{X}$ as an admissible variation of Hodge structure. Moreover, we have that $\rho = \sum_p  \rho^{p,p}\,q^p$ where $\rho^{p,p}: \pi_1(B) \to \Aut(V^{p,p})$ is the restriction of $\rho$ to the Hodge pieces of $V$.

In particular, consider a regular morphism $h: X \to B$ that is a locally trivial fibration in the analytic topology, and whose fiber $F$ carries a balanced Hodge structure. It induces monodromy representations on the compactly supported cohomology $\rho: \pi_1(B) \to \Aut(H_c^\bullet(F))$. Since the derived category is flattened in the Grothendieck ring via an alternating sum, we have that, as an element of $\KM{B}$,
$$
	\RM{X}{B} = \sum_k (-1)^k\,\rho|_{H_c^k(F)} = \sum_{k,p}(-1)^k\, H_c^{k;p,p}(F;\QQ)\,q^p.
$$
This can be seen as an element of the representation ring of $\pi_1(B)$ tensorized with $\ZZ[q^{\pm 1}]$. This is precisely the definition of \cite{LMN} of Hodge monodromy representations. In particular, reinterpreting Hodge monodromy representations as mixed Hodge modules, all the computations of \cite{LMN, Martinez:2017, MM} can be seen as calculations of mixed Hodge modules and, at the end of the day, as calculations of virtual classes of algebraic varieties.

\subsection{Virtual covering spaces}
\label{sec:Hodge-monodromy-covering-spaces}

Suppose that $X$ and $B$ are smooth complex varieties and that $f: X \to B$ is a regular morphism which is a covering space in the analytic topology with finite fiber $F$ and degree $d$. In that case, the monodromy action $\pi_1(B) \to \GL{d}(\QQ)$ determines completely the covering. For this reason, the class $\RM{X}{B} = [X, f] \in \KVarrel{B}$ only depends on the monodromy of $f$, which justifies the name Hodge monodromy.

In particular, for the purposes of this paper, we will focus on the following case. Fix different points $p_1, \ldots, p_s \in \CC$ and take $B = \CC-\left\{p_1, \ldots, p_s\right\}$. For any $p_i$, we will consider the variety
$$
	X_{p_i} = \left\{(x, y) \in B \times \CC \,|\, y^2 = x-p_i\right\},
$$
with projection $X_{p_i} \to B$, $f(x,y)=x$. It defines a double covering whose monodromy action is trivial for the loops around $p_j$ with $j \neq i$ and the loop around $p_i$, $\gamma_{p_i}$, permutes both fibers
$$
	\gamma_{p_i} \mapsto \begin{pmatrix}
	0 & 1\\
	1 & 0\\
\end{pmatrix} \sim \begin{pmatrix}
	1 & 0\\
	0 & -1\\
\end{pmatrix}.
$$
For this reason, it makes sense to define $S_{p_i} = \RM{X_{p_i}}{B} - T_{B}$. It may be understood as the non-trivial $\ZZ_2$-representation over $\gamma_{p_i}$.

\begin{ex}
Observe that $X_{p_i}$ is isomorphic to an affine parabola minus $2s-1$ points. Therefore, we have that $[X_{p_i}] = [\CC] - 2s-1$ and $[B] = [\CC] - s$. Thus, if $c: B \to \star$ is the final morphism, we have that
$$
	[\CC] - 2s-1 = c_! \RM{X_{p_i}}{B} = c_!(T_{B}) + c_!(S_{p_i}) = [\CC] - s + c_!(S_{p_i}).
$$
This implies that the characteristic is $c_!(S_{p_i}) = -s-1$. This agrees with discussion after \cite[Theorem 6]{MM}. Actually, a modification of this argument shows that also that $c_!(S_{p_{i_1}} \times \ldots \times S_{p_{i_r}}) = -s-1$ for any $r \geq 0$.
\end{ex}

\begin{ex}\label{ex:monodromy-plane-minus-points}
Let $X = \CC^* - \left\{\pm 1\right\}$ and $B = \CC - \left\{\pm 2\right\}$ with the double cover $t: X \to B$ given by $t(\lambda) = \lambda + \lambda^{-1}$ for $\lambda \in \CC^* - \left\{\pm 1\right\}$. Let us consider the variety
$$
	X_{2,-2} = \left\{(x,y) \in B \times \CC \,|\, y^2 = x^2-4\right\},
$$
with projection $X_{2,-2} \to B$, $(x,y) \mapsto x$. We have an isomorphism $\CC^* - \left\{\pm 1\right\} \cong X_{2,-2}$ given by $\lambda \mapsto (\lambda + \lambda^{-1}, \lambda - \lambda^{-1})$. Hence $\RM{\CC^* - \left\{\pm 1\right\}}{\CC - \left\{\pm 2\right\}} = [X_{2,-2}]$, which can be shown to be equal to $T_{\CC - \left\{\pm 2\right\}} + S_2 \times S_{-2}$.

This is compatible with the interpretation of $S_2, S_{-2}$ as one-dimensional monodromies. Indeed, if we take $\pi_1(\CC - \left\{\pm 2\right\}) = \langle \gamma_2, \gamma_{-2}\rangle$, with $\gamma_{\pm 2}$ the positive small loops around $\pm 2$, then the monodromy action is given by
$$
	\gamma_2 \mapsto \begin{pmatrix}
	0 & 1\\
	1 & 0\\
\end{pmatrix} \sim \begin{pmatrix}
	1 & 0\\
	0 & -1\\
\end{pmatrix}, \hspace{2cm} \gamma_{-2} \mapsto \begin{pmatrix}
	0 & 1\\
	1 & 0\\
\end{pmatrix} \sim \begin{pmatrix}
	1 & 0\\
	0 & -1\\
\end{pmatrix}.
$$
That is, the monodromy action decomposes as a trivial representation and the product of the non-trivial representations around $\pm 2$.
\end{ex}

\section{TQFTs for representation varieties}
\label{sec:repr-varieties}

In this section we will review the definition of representation and character varieties. They are the central objects of this paper. We will work over an algebraically closed field $k$.

\begin{defn}
Let $\Gamma$ be a finitely generated groupoid i.e.\ a groupoid with finite many objects and whose vertex groups are all finitely generated \cite{GPLM-2017}, and let $G$ be an algebraic group over $k$. The set of groupoid representations of $\Gamma$ into $G$
$$
	\Rep{G}(\Gamma)=\Hom(\Gamma, G)
$$
is called the \emph{representation variety}. 
\end{defn}

As its name suggests, $\Rep{G}(\Gamma)$ has a natural algebraic structure. First of all, suppose that $\Gamma$ is a group. In that case, let $\Gamma = \langle \gamma_1, \ldots, \gamma_r\;|\; R_\alpha(\gamma_1, \ldots, \gamma_r)=1\rangle$ be a presentation of $\Gamma$ with finitely many generators, where $R_\alpha$ are the relations (possibly infinitely many). We define the injective map $\psi: \Hom(\Gamma, G) \to  G^r$ given by $\psi(\rho)=(\rho(\gamma_1), \ldots, \rho(\gamma_r))$. The image of $\psi$ is the algebraic subvariety of $G^r$ 
$$
	\img{\psi} = \left\{(g_1, \ldots, g_r) \in G^r\;\right|R_\alpha(g_1, \ldots, g_r)=1\left.\right\}.
$$
Hence, we can impose an algebraic structure on $\Hom(\Gamma, G)$ by declaring that $\psi$ is a regular isomorphism on its image. Observe that this algebraic structure does not depend on the chosen presentation.

In the general case in which $\Gamma$ is a groupoid with $n$ nodes, pick a set $J = \left\{a_1, \ldots, a_s\right\}$ of objects of $\Gamma$ such that every connected component of $\Gamma$ contains exactly one element of $J$. Moreover, for any object $a$ of $\Gamma$, pick a morphism $f_a: a \to a_i$ where $a_i$ is the object of $J$ in the connected component of $a$. Hence, if $\rho: \Gamma \to G$ is a groupoid homomorphism, it is uniquely determined by the group representations $\rho_i: \Gamma_{a_i} \to G$ for $a_i \in J$, where $\Gamma_{a_i} = \Hom_{\Gamma}(a_i, a_i)$ is the vertex group at $a_i$, together with the elements $g_a \in G$ corresponding to the morphisms $f_a$ for any object $a$. Since the elements $g_a$ can be chosen without any restriction, we have a natural identification
$$
	\Hom(\Gamma, G) \cong \Hom(\Gamma_{a_1},G) \times \ldots \times \Hom(\Gamma_{a_s}, G) \times G^{n-s},
$$
and each of these factors has a natural algebraic structure as representation variety. This endows $\Hom(\Gamma, G)$ with an algebraic structure.

The representation variety $\Rep{G}(\Gamma)$ has a natural action of $G$ by conjugation. Recall that two representations $\rho, \rho'$ are said to be isomorphic if $\rho' = g \cdot \rho$ for some $g \in G$.

\begin{defn}
Let $\Gamma$ be a finitely generated groupoid and $G$ an algebraic reductive group. The Geometric Invariant Theory quotient
$$
	\Char{G}(\Gamma) = \Rep{G}(\Gamma) \sslash G,
$$
is called the \emph{character variety}.
\end{defn}

\begin{rmk}
The Geometric Invariant Theory quotient, also known as GIT quotient, is a kind of quotient that make sense for algebraic varieties. For a complete introduction to this topic, see \cite{MFK:1994, Newstead:1978}.
\end{rmk}

\begin{ex}\label{ex:gamma-for-repr}
Let $M$ be a compact connected manifold with fundamental group $\Gamma = \pi_1(M)$. The fundamental group of such a manifold is finitely generated since a compact manifold has the homotopy type of a finite CW-complex. Hence, we can form its representation variety, that we will shorten $\Rep{G}(M) = \Rep{G}(\pi_1(M))$. The corresponding character variety is called the character variety of $M$ and it is denoted by $\Char{G}(M) = \Rep{G}(M) \sslash G$. More generally, if $A \subseteq M$ is a finite set, the fundamental groupoid $\Pi(M, A)$ is a finitely generated groupoid. Recall that $\Pi(M, A)$ is the groupoid whose elements are homotopy classes of paths in $M$ between points in $A$. The corresponding representation variety is denoted $\Rep{G}(M,A)$ and the associated character variety $\Char{G}(M,A)$.
\end{ex}

A step further in the construction of character varieties can be done by considering an extra structure on them, called a parabolic structure. A \emph{parabolic structure} on $\Rep{G}(\Gamma)$, $Q$, is a finite set of pairs $(\gamma, \lambda)$, where $\gamma \in \Gamma$ and $\lambda \subseteq G$ is a locally closed subset which is closed under conjugation. Given such a parabolic structure, we define the \emph{parabolic representation variety}, $\Rep{G}(\Gamma, Q)$, as the subset of $\Rep{G}(\Gamma)$
$$
	\Rep{G}(\Gamma, Q) = \left\{\rho \in \Rep{G}(\Gamma)\,\,|\,\, \rho(\gamma) \in \lambda\,\,\textrm{for all } (\gamma, \lambda) \in Q\right\}.
$$
As in the non-parabolic case, $\Rep{G}(\Gamma, Q)$ has a natural algebraic variety structure since, using suitable generators, we have an identification $\Rep{G}(\Gamma, Q) = \Rep{G}(\Gamma) \cap \left(G^r \times \lambda_1 \times \ldots \times \lambda_s\right)$. The conjugacy action of $G$ on $\Rep{G}(\Gamma)$ restricts to an action on $\Rep{G}(\Gamma, Q)$ since the subsets $\lambda_i$ are closed under conjugation. The GIT quotient of the representation variety by this action,
$$
	\Char{G}(\Gamma, Q) = \Rep{G}(\Gamma, Q) \sslash G,
$$
is called the \emph{parabolic character variety}.

\begin{rmk}
As a particular choice for the subset $\lambda \subseteq G$, we can take the conjugacy classes of an element $h \in G$, denoted $[h]$. Observe that $[h] \subseteq G$ is locally closed since, by \cite[Lemma 3.7]{Newstead:1978}, it is an open subset of its Zariski closure.
\end{rmk}

\begin{ex}\label{ex:par-struc-repr}
Let $\Sigma = \Sigma_g - \left\{p_1, \ldots, p_s\right\}$ with $p_i \in \Sigma_g$ distinct points, called the punctures or the marked points. In that case, we have a presentation of the fundamental group of $\Sigma$ given by
$$
	\pi_1(\Sigma) = \left\langle \alpha_1, \beta_1 \ldots, \alpha_{g}, \beta_g, \gamma_1, \ldots, \gamma_s\;\;\left|\;\; \prod_{i=1}^g [\alpha_i, \beta_i]\prod_{j=1}^s \gamma_s = 1 \right.\right\rangle,
$$
where the $\gamma_i$ are the positive oriented simple loops around the punctures. As parabolic structure, we will take $Q = \left\{(\gamma_1, \lambda_1), \ldots, (\gamma_s, \lambda_s)\right\}$. The corresponding parabolic representation variety is
$$
	\Rep{G}(\pi_1(\Sigma), Q) = \left\{(A_1, B_1 \ldots, A_{g}, B_g, C_1, \ldots, C_s) \in G^{2g+s}\;\;\left|\;\; \begin{matrix}
	\displaystyle \prod_{i=1}^g [A_i, B_i]\prod_{j=1}^s C_s = 1 \\
	C_j \in \lambda_j
\end{matrix}	
	\right.\right\}.
$$
\end{ex}

Let $M$ be a compact differentiable manifold and let $S \subseteq M$ be a codimension $2$ closed connected submanifold with a co-orientation. Embed the normal bundle as a small tubular neighborhood $U \subseteq M$ around $S$. Fixed $s \in S$, consider an oriented local trivialization $\psi: V \times \RR^2 \to U$ of the normal bundle around an open neighborhood $V \subseteq S$ of $s$. In that situation, the loop $\gamma(t) = \psi(s, (\cos{t}, \sin{t})) \in \pi_1(M - S)$ is called the \emph{positive meridian} around $s$.

\begin{figure}[h]
	\begin{center}
	\includegraphics[scale=0.25]{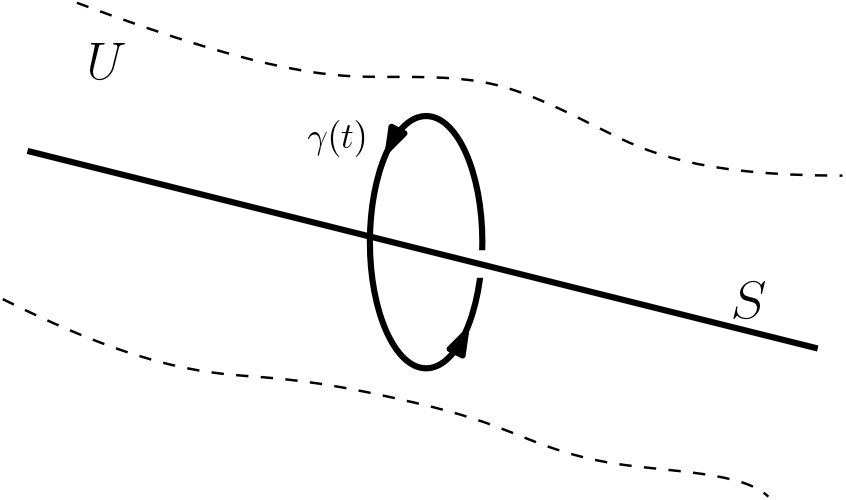}
	\end{center}
\vspace{-0.5cm}
\end{figure}

\begin{rmk}
\begin{itemize}
	\item Two positive meridians around $s, s' \in S$ are conjugate to each other on $\pi_1(M-S)$.
	\item If $S$ is not connected, then the conjugacy classes of meridians are in correspondence with the connected components of $S$.
	\item A loop $\gamma \in \pi_1(M-S)$ is a generalized knot. In that sense, the kernel of $\pi_1(M - S) \to \pi_1(M)$ are the loops in $M-S$ that `surround' $S$. In \cite{Smith:1978}, it is proven that this kernel is the smallest normal subgroup containing all the meridians around the connected components of $S$.
\end{itemize}
\end{rmk}

Take $\Lambda$ to be a collection of locally closed subsets of $G$ that are invariant under conjugation. Suppose that $Q=\left\{(S_1, \lambda_1), \ldots, (S_s, \lambda_s)\right\}$ is a parabolic structure over $\Lambda$, in the sense of Section \ref{sec:some-extra-struct}. Let us denote $S = \bigcup_i S_i$. To this parabolic structure, we can build the parabolic structure on $\Rep{G}(\Pi\left(M - S, A\right))$, where $A \subseteq M$ is any finite set, also denoted by $Q$. It is the parabolic structure $Q = \left\{(\gamma_{1,1}, \lambda_{1}), \ldots, (\gamma_{1,m_1}, \lambda_{1}), \ldots, (\gamma_{s, 1}, \lambda_s), \ldots, (\gamma_{s, m_s}, \lambda_s)\right\}$, where $\gamma_{i, 1}, \ldots, \gamma_{i, m_i}$ is a generating set of positive meridians around the connected components of $S_i$ (based on any point of $A$). As for the non-parabolic case, we will shorten the corresponding parabolic character variety by $\Rep{G}(M, A, Q) = \Rep{G}(\Pi\left(M-S, A\right), Q)$.

\begin{rmk}
Suppose that $M=\Sigma_g$ is a closed oriented surface. A parabolic structure is given by $Q=\left\{(p_1, \lambda_1), \ldots, (p_s, \lambda_s)\right\}$, with $p_i \in \Sigma_g$ points with a preferred orientation of a small disk around them (see Example \ref{ex:parab-struct}). In that case, the meridian of $p_i$ is given by a small loop encycling $p_i$ positively with respect to the orientation of the small disk around it. Therefore, the associated parabolic structure of representation variety is
$$
	\Rep{G}(\Sigma_g, Q) = \left\{(A_1, B_1 \ldots, A_{g}, B_g, C_1, \ldots, C_s) \in G^{2g+s}\;\;\left|\;\; \begin{matrix}
	\displaystyle \prod_{i=1}^g [A_i, B_i]\prod_{j=1}^s C_s = 1 \\
	C_j^{\epsilon_j} \in \lambda_j
\end{matrix}	
	\right.\right\},	
$$
where $\epsilon_j = 1$ if the orientation of the disk around $p_j$ agrees with the global orientation and $\epsilon_j = -1$ if it does not. Notice that they agree with the ones of Example \ref{ex:par-struc-repr}.
\end{rmk}

\subsection{Standard TQFT for representation varieties}

\label{sec:TQFT-for-repr}

In this section, we sketch the construction of a lax monoidal parabolic TQFT of pairs computing virtual classes of representation varieties. This construction appeared for the first time in \cite{GPLM-2017} in the context of computation of virtual Hodge structures. In this paper, we will extend it to compute their virtual classes as algebraic varieties and we will propose several variants that will be useful for computational purposes.

Fix $G$ an algebraic group and take as parabolic data $\Lambda$ a collection of subvarieties of $G$ that are closed under conjugation. As the category of fields for this construction, we take the category of algebraic varieties, $\Var{k}$. The geometrisation functor $\Rep{G}: \Embpparc{\Lambda} \to \Var{k}$ is given as follows. For $(M, A, Q) \in \Embpparc{\Lambda}$, where $M$ is a compact manifold, $A \subseteq M$ is a finite set and $Q$ is a parabolic structure on $M$, it assigns $\Rep{G}(M, A, Q)$, the parabolic representation variety of $(M, A)$ on $G$ with parabolic structure $Q$.

On the other hand, to a morphism $(f, \alpha): (M, A, Q) \to (M', A', Q')$ in $\Embpparc{\Lambda}$, we associate the regular morphism $\Rep{G}(f): \Rep{G}(M',A',Q') \to \Rep{G}(M, A, Q)$ induced by the groupoid homomorphism $f_*: \Pi(M, A) \to \Pi(M',A')$. Observe that, as the morphism $\alpha$ is just inclusion (of points and parabolic structures), the morphism $f_*$ preserves the parabolic structures. By the Seifert-van Kampen theorem for fundamental groupoids \cite{Brown:1967}, the functor $\Rep{G}$ has the Seifert-van Kampen property so it gives rise to a field theory $\Fld{\Rep{G}}: \Bordppar{n}{\Lambda} \to \Span{\Var{k}}$.

Finally, as $\Var{k}$-algebra, we will take the one constructed from the Grothendieck rings of algebraic varieties, $\KVarrel{}$, as described in Section \ref{sec:grothendieck-rings-Var-algebras}.
%Saito's mixed Hodge modules. We are going to take $\KM{} = (A,B)$. On objects, both functors assign to $X \in \CVar$ the Grothendieck group of the category mixed Hodge modules on $X$, $A(X) = B(X) = \KM{X}$, seen as a $\KM{\star} = \KVar{k}$-algebra. With respect to morphisms, given a regular map $f: X \to Y$ between complex algebraic varieties we assign $A(f)=f^*: \KM{Y} \to \KM{X}$ and $B(f)=f_!: \KM{X} \to \KM{Y}$. By the results of Section \ref{sec:grothendieck-rings-mhm}, $f^*$ is a ring homomorphism and $f_!$ is a $\KVar{k}$-module homomorphism and they commute with external product, as needed for $\CVar$-algebras. Finally, they satisfy the Beck-Chevalley condition by item \ref{thm:funct-induced-functors-mhm:beck-chevalley} of Theorem \ref{thm:funct-induced-functors-mhm}. Hence $\KM{}$ is a $\CVar$-algebra.
Therefore, Theorem \ref{thm:physical-constr-sheaf} implies the following result.

\begin{thm}
Let $G$ be an algebraic group and $n \geq 1$. There exists a lax monoidal TQFT
$$
	Z_{\Rep{G}, \KVarrel{}}: \Bordppar{n}{\Lambda} \to \Mod{\K{\Var{k}}},
$$
computing virtual classes of parabolic $G$-representation varieties. 
\end{thm}
For short we will denote $\Zs{G} = Z_{\Rep{G}, \KVarrel{}}$ and we will call it the \emph{standard TQFT}. Using the properties explained in Remark \ref{rmk:expression-TQFT}, this TQFT satisfies that, for any $n$-dimensional connected closed manifold $W$, any non-empty finite subset $A \subseteq W$ and any parabolic structure $Q$ in $W$, we have
$$
	\Zs{G}(W, A, Q)\left(\star\right) = \coh{\Rep{G}(W, Q)} \times \coh{G}^{|A|-1}
$$ 
where $\star = [\star] \in \KVar{k}$ is the unit of the ring. 

\begin{rmk}\label{rmk:Deligne-Hodge-soft-TQFT}
The virtual class $\coh{G} \in \KVar{k}$ is usually known for the standard groups.
\end{rmk}

\begin{figure}[h]
	\begin{center}
	\includegraphics[scale=0.65]{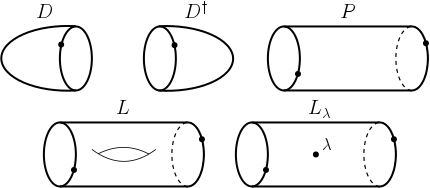}
	\caption{}
	\label{img:gen-tubes-pairs}
	\end{center}
\vspace{-0.3cm}
\end{figure}

In the case $n=2$ (surfaces), this TQFT was explicitly described in \cite[Section 4.3]{GPLM-2017} but using mixed Hodge modules instead of motivic classes. Observe that the morphisms of $\Tubppar{2}{\Lambda} = \ETub{2}{\cS_{p,\Lambda}}$ are generated by the set $\Delta = \left\{D, D^\dag, L, P\right\} \cup \bigcup_{\lambda \in \Lambda} \left\{L_\lambda\right\}$, as depicted in Figure \ref{img:gen-tubes-pairs}.

Unraveling the definitions of the TQFT, the associated morphisms of the discs $D$ and $D^\dag$ under $\Zs{G}$ are
$$
	\Zs{G}(D) = i_!: \KVar{k}=\KVarrel{1} \to \KVarrel{G}, \hspace{1cm} \Zs{G}(D^\dag) = i^*: \KVarrel{G} \to \KVarrel{1} = \KVar{k}.
$$
where $i: 1 \hookrightarrow G$ is the inclusion. For the holed torus $L: (S^1,\star) \to (S^1, \star)$ the situation is a bit more complicated. The associated field theory is the span
$$
\begin{matrix}
	G & \stackrel{p}{\longleftarrow} & G^4 & \stackrel{q}{\longrightarrow} & G \\
	g & \mapsfrom & (g, g_1, g_2, h) & \mapsto & hg[g_1,g_2]h^{-1}
\end{matrix}
$$
Hence, the image of $L$ under the TQFT is the morphism $\Zs{G}(L) = q_! \circ p^*: \KVarrel{G} \to \KVarrel{G}$.
For the morphism $P$, the associated field theory is the span
$$
\begin{matrix}
	G & \stackrel{u}{\longleftarrow} & G^2 & \stackrel{v}{\longrightarrow} & G \\
	g & \mapsfrom & (g, h) & \mapsto & hgh^{-1}
\end{matrix}
$$
Thus, we have that $\Zs{G}(P) = v_! \circ u^*: \KVarrel{G} \to \KVarrel{G}$.
Finally, for the tube with parabolic datum $\lambda$, $L_\lambda$, we have that its image under the field theory is the span
$$
\begin{matrix}
	G & \stackrel{r}{\longleftarrow} & G^2 \times \lambda & \stackrel{s}{\longrightarrow} & G \\
	g & \mapsfrom & (g, h, \xi) & \mapsto & hg\xi h^{-1}
\end{matrix}
$$
Therefore, its image under the TQFT is the morphism $\Zs{G}(L_\lambda) = s_! \circ r^*: \KVarrel{G} \to \KVarrel{G}$.
	
With this description, we have proven the following result.

\begin{thm}\label{thm:almost-tqft-parabolic}
Let $\Sigma_g$ be the closed oriented surface of genus $g$ and let $Q$ be a parabolic structure on $\Sigma_g$ with $s$ marked points with data $\lambda_1, \ldots, \lambda_s \in \Lambda$. Then, 
$$
	\coh{\Rep{G}(\Sigma_g, Q)} = \frac{1}{\coh{G}^{g+s}}{\Zs{G}(D^\dag) \circ \Zs{G}(L_{\lambda_s}) \circ \ldots \circ \Zs{G}(L_{\lambda_1}) \circ \Zs{G}(L)^g \circ \Zs{G}(D) (\star)}.
$$
\end{thm}

\subsection{Piecewise algebraic varieties}
\label{subsec:piecewise-alg-var}
In this section, we extend slightly the notion of an algebraic variety in order to allow topological spaces which are not algebraic varieties, but disjoint union of algebraic varieties. This is similar to the construction of the Grothendieck ring $\KVar{k}$, but not adding formal additive inverses and tracking the morphisms. The idea is to allow us to consider a variation of the TQFT above that, instead of GIT quotients, uses usual quotients as orbit spaces.

\begin{defn}
Let $k$ be an algebraically closed field. We define the \emph{category of piecewise varieties over $k$}, $\PVar{k}$, as the category given by:
\begin{itemize}
	\item Objects: The objects of $\PVar{k}$ is the Grothendieck semi-ring of algebraic varieties over $k$. This is the semi-ring generated by the symbols $[X]$, for $X$ an algebraic variety, with the relation $[X] = [U] + [Y]$ if $X = U \sqcup Y$ with $U \subseteq X$ open and $Y \subseteq X$ closed. It is a semi-ring with disjoint union as sum and cartesian product as product.
	\item Morphisms: Given $[X], [Y] \in \PVar{k}$ a morphism $f: [X] \to [Y]$ is given by a function $f: X \to Y$ that decomposes as a disjoint union of regular maps. More precisely, it has to exist decompositions $[Y] = \sum_i [Y_i]$ and $[X] = \sum_{i,j} [X_{ij}]$ into algebraic varieties such that $f = \bigsqcup\limits_{i,j} f_{ij}$ with $f_{ij}: X_{ij} \to Y_i$ a regular morphism.
	\item Composition is given by the usual composition of maps. Observe that, given $f: [X] \to [Y]$ and $g: [Y] \to [Z]$, by Chevalley theorem \cite{EGAI}, we can find a common decomposition $Y = \sqcup_{i,j} Y_{ij}$ such that we get a decomposition into regular maps
	$$
	f = \bigsqcup_{i,j,k} f_{ijk}: \bigsqcup_{i,j,k} X_{ijk} \to \bigsqcup_{i,j} Y_{ij}, \hspace{1cm} g = \bigsqcup_{i,j} g_{ij}: \bigsqcup_{i,j} Y_{ij} \to \bigsqcup_i Z_i.
$$
Hence, $g \circ f = \bigsqcup\limits_{i,j,k} \left(g_{ij} \circ f_{ijk}\right)$ is a decomposition into regular maps of $g \circ f$.
\end{itemize}
\end{defn}

\begin{rmk}\label{rmk:forgetful-piecewise-var}
\begin{itemize}
	\item There is a functor $\Var{k} \to \PVar{k}$ that sends $X \mapsto [X]$ and analogously for regular maps. On the other way around, there is a forgetful functor $\PVar{k} \to \Sets$ that recovers the underlying set of a piecewise variety.
	\item The functor $\Var{k} \to \PVar{k}$ is not an embedding of categories. For example, let $X = \left\{y^2 = x^3\right\}$ be a cuspidal cubic affine plane curve. Observe that, removing the origin in $X$, we have a decomposition $[X] = [X - \star] + [\star] = [\mathbb{A}^1 - \star] + [\star] = [\mathbb{A}^1]$. Hence, the images of $X$ and $\mathbb{A}^1$ under this functor agree, even though they are not isomorphic.
\end{itemize}
\end{rmk}
	
\begin{ex}
The affine line with a double point at the origin $X$ is not an algebraic variety. However, it is a piecewise variety since we can decompose $[X] = [\mathbb{A}^1] + [\star]$ with both pieces locally closed subsets.
\end{ex}

Let $G$ be an algebraic reductive group acting on a variety $X$. Maybe after restricting to the open subset of semi-stable points, we can suppose that a (good) GIT quotient $\pi: X \to X \sslash G$ exists. On $X$, we find the subset $X_1 \subseteq X$ of poly-stable points whose orbits have maximum dimension. It is an open subset by an adaptation of \cite[Proposition 3.13]{Newstead:1978} and it is non-empty since the set of poly-stable points is so. On the poly-stable points, the GIT quotient is an orbit space so we have a $G$-invariant regular map $\pi_1: X_1 \to X_1 / G$.
Now, let $Y = X - X_1$. As $Y$ is closed on $X$ and the action of $G$ restricts to an action on $Y$ and we can repeat the argument to obtain a regular $G$-invariant map $\pi_2: X_2 \subseteq Y \to X_2 / G$, where $X_2 \subseteq Y$ is the set of poly-stable points of $Y$ of maximum dimension. Repeating this procedure, we obtain a stratification $X = X_1 \sqcup \ldots \sqcup X_r$ and a set of regular maps $\pi_i: X_i \to X_i/G$ where each $X_i / G$ has a natural algebraic structure.

\begin{defn}\label{defn:piecewise-quotient}
Let $G$ be a reductive group and let $X$ be an algebraic variety. With the decomposition above the \emph{piecewise quotient} of $X$ by $G$, denoted by $[X / G]$, is the object of $\PVar{k}$ given by
$$
	[X/G] = [X_1/G] + \ldots + [X_r / G].
$$
We also have a piecewise quotient morphism $\pi = \sqcup_i \pi_i: [X] \to [X/G]$ where $\pi_i: X_i \to X_i/G$ are the projections on the orbit space of each of the strata.
\end{defn}

\begin{rmk}
Since all the strata considered are made of poly-stable points over the previous stratum, no orbits are identified under the quotient. In this sense, the previous construction says that the orbit space of $X$ by $G$ has the structure of a piecewise variety.
\end{rmk}

\begin{ex}
Let $G$ be an algebraic group. Given a pair of topological spaces $(X, A)$ such that $\Pi(X,A)$ is finitely generated and a parabolic structure $Q$ on it, we define the \emph{piecewise character variety} as the object of $\PVar{k}$
$$
	\CharW{G}(X, A, Q) = [\Rep{G}(X,A, Q) / G].
$$ 
Here, $G$ acts on $\Rep{G}(X,A, Q)$ by conjugation.
\end{ex}

The piecewise quotients have similar properties than usual GIT quotients but in the category $\PVar{k}$. More precisely, let $G$ be a reductive group acting on a variety $X$ and let $f: [X] \to [Y]$ be a morphism of piecewise varietes which is $G$-invariant. Using the universal property of categorical quotients, we have that there exists an unique piecewise regular morphism $\tilde{f}: [X / G] \to [Y]$ such that $\tilde{f} \circ \pi = f$.
	\[
\begin{displaystyle}
   \xymatrix
   {	[X] \ar[r]^f \ar[d]_{\pi} & [Y] \\
   		[X/G] \ar@{--{>}}[ru]_{\tilde{f}}
   }
\end{displaystyle}   
\]

The category $\PVar{k}$ is just an intermediate step in the path from the category $\Var{k}$ to its Grothendieck ring $\KVar{k}$. For this reason, the $\Var{k}$-algebra $\KVarrel{}$ of Section \ref{sec:grothendieck-rings-Var-algebras} extends in a natural way to a $\PVar{k}$-algebra. That will be very useful in the following Section \ref{sec:geometric-reduced-TQFT}.

\subsection{Geometric and reduced TQFT}
\label{sec:geometric-reduced-TQFT}

Using the reduction method of Section \ref{sec:reduction-TQFT}, here we will describe an step further in the standard TQFT of Section \ref{sec:TQFT-for-repr}. 

From an algebraic group $G$, we considered the representation variety geometrisation $\Rep{G}: \Embpparc{\Lambda} \to \Var{k}$. After the functor $\Var{k} \to \PVar{k}$ of Remark \ref{rmk:forgetful-piecewise-var}, we may consider the reduction $\pi$ that assigns, to $(M, A, Q) \in \Bordppar{n}{\Lambda}$, the piecewise character variety $\CharW{G}(M, A,Q) = [\Rep{G}(M, A, Q) / G]$ together with the piecewise quotient morphism $\pi_{(M, A, Q)}: \Rep{G}(M, A, Q) \to \CharW{G}(M, A, Q)$. As explained in Section \ref{sec:reduction-TQFT} this reduction can be used to modify the standard TQFT, $\Zs{G}$, of Proposition \ref{prop:existence-reduction-TQFT} in order to obtain a very lax TQFT $\Zg{G} = (\Zs{G})_{\pi}: \Bordppar{n}{\Lambda} \to \Modt{\KVar{k}}$. 

Using the notation of Section \ref{sec:reduction-TQFT}, suppose that, for all $(M, A, Q) \in \Tubppar{n}{\Lambda}$, the morphism $\eta_{(M, A, Q)} = (\pi_{(M, A, Q)})_! \circ (\pi_{(M, A, Q)})^*: (\cV_{(M, A, Q)})_0 \to (\cV_{(M, A, Q)})_0$ is invertible as a $\Ko{\Var{k}}$-module homomorphism. In that case, the morphisms $\cZg{G}(W,A,Q) = \Zg{G}(W,A,Q) \circ \eta_{(M, A,Q)}^{-1}$ give rise to a strict TQFT, the $\pi$-reduction of $\Zs{G}$ (Proposition \ref{prop:existence-reduction-TQFT}). Here, $(\cV_{(M, A, Q)})_0$ stands for the localization $S^{-1}\cV_{(M, A, Q)}$ by an appropriate multiplicative system $S \subseteq \KVar{k}$.

\begin{defn}
Let $G$ be an algebraic group and $\Lambda$ a collection of conjugacy closed subvarieties of $G$. The Topological Quantum Field Theories
$$
	\Zg{G}: \Bordppar{n}{\Lambda} \to \Modt{\KVar{k}}, \hspace{1cm} \cZg{G}: \Bordppar{n}{\Lambda} \to \Modt{\Ko{\Var{k}}},
$$
are called the \emph{geometric TQFT} and the \emph{reduced TQFT} for representation varieties, respectively. 
\end{defn}
These TQFTs will be extensively used along Section \ref{sec:sl2-repr-var} since, as we will see, they allow us to compute virtual classes of representation varieties in a simpler way than the standard TQFT.

\begin{ex}
In the case of surfaces, the morphisms of $\Zg{G}$ can be explicitly written down from the description of Section \ref{sec:TQFT-for-repr}. Consider the set of generators $\Delta = \left\{D, D^\dag, L, P\right\} \cup \left\{L_\lambda\right\}_{\lambda \in \Lambda}$ of $\Tubppar{2}{\Lambda}$ as depicted in Figure \ref{img:gen-tubes-pairs}. Since, by definition, $\Zg{G}(W, A, Q) = \pi_! \circ Z(W, A, Q) \circ \pi^*$ for any bordism $(W, A, Q)$, then we have that
\begin{align*}
	&\Zg{G}(D) = (\pi \circ i)_!: \KVar{k} \to \KVarrel{[G/G]}, \hspace{1cm} \Zg{G}(D^\dag) =  (\pi \circ i)^*: \KVarrel{[G/G]} \to \KVar{k}, \\
	&\Zg{G}(L) = \hat{q}_! \circ \hat{p}^*: \KVarrel{[G/G]} \to \KVarrel{[G/G]}, \hspace{0.4cm} \Zg{G}(P) = \hat{v}_! \circ \hat{u}^*: \KVarrel{[G/G]} \to \KVarrel{[G/G]},\\
 & \hspace{3.6cm}\Zg{G}(L_{\lambda}) = \hat{s}_! \circ \hat{r}^*: \KVarrel{[G/G]} \to \KVarrel{[G/G]},
\end{align*}
where $\hat{q} = \pi \circ q, \hat{p} = \pi \circ p: G^4 \to [G/G]$, $\hat{v} = \pi \circ v, \hat{u} = \pi \circ u: G^2 \to [G/G]$ and $\hat{r} = \pi \circ r, \hat{s} = \pi \circ s: G^2 \times \lambda \to [G/G]$. Observe that $\pi_{\emptyset}$ is the identity map since $\Rep{G}(\emptyset) = \CharW{G}(\emptyset) = 1$ and that $\Zg{G}(S^1,\star)= [\Rep{G}(S^1, \star)/G] = [G/G]$.

Finally, as mentioned in Example \ref{ex:reduction-TQFT-surfaces}, the only relevant submodule is $\cV = \cV_{(S^1, \star)} \subseteq \KVarrel{\CharW{G}(S^1, \star)} =  \KVarrel{[G/G]}$, so we can restrict the previous maps to $\cV$. Moreover, if the morphism $\eta = \eta_{(S^1, \star)}: \cV_0 \to \cV_0$ is invertible, the reduction $\cZg{G}$ can also be constructed. 
\end{ex}

\section{$\SL{2}(\CC)$-representation varieties}
\label{sec:sl2-repr-var}

For the rest of the paper, we will focus on the case of surfaces ($n = 2$), the ground field $k=\CC$ and $G = \SL{2}(\CC)$, the complex special linear group of order two. In order to enlighten the notation, we shall write $\Var{\CC} = \Var{}$, $\Zs{\SL{2}(\CC)}=\Zs{}$, $\Zg{\SL{2}(\CC)}=\Zg{}$ and $\cZg{\SL{2}(\CC)}=\cZg{}$. As an application of the previous TQFTs, we will compute the virtual classes of some parabolic representations varieties for arbitrary genus. For that, we will give explicit expressions of the homomorphisms of Section \ref{sec:geometric-reduced-TQFT}.

In $\SL{2}(\CC)$, there are five special types of elements, namely the matrices
$$
	\Id = \begin{pmatrix}
	1 & 0\\
	0 & 1\\
\end{pmatrix}, \hspace{0.2cm}
	-\Id = \begin{pmatrix}
	-1 & 0\\
	0 & -1\\
\end{pmatrix}, \hspace{0.2cm}
	J_+ = \begin{pmatrix}
	1 & 1\\
	0 & 1\\
\end{pmatrix}, \hspace{0.2cm}
	J_- = \begin{pmatrix}
	-1 & 1\\
	0 & -1\\
\end{pmatrix}, \hspace{0.2cm}
D_\lambda = \begin{pmatrix}
	\lambda & 0\\
	0 & \lambda^{-1}\\
\end{pmatrix},
$$
with $\lambda \in \CC-\left\{0, \pm 1\right\}$. Observe that any element of $\SL{2}(\CC)$ is conjugated to one of those elements. Such a distinguished representant is unique up to the fact that $D_\lambda$ and $D_{\lambda^{-1}}$ are conjugated for all $\lambda \in \CC-\left\{0, \pm 1\right\}$.

Hence, if we denote the orbit of $A \in \SL{2}(\CC)$ under conjugation by $[A]$, we have a stratification
\begin{equation} \label{eqn:dec-sl2}
	\SL{2}(\CC) = \left\{\Id\right\} \sqcup \left\{-\Id\right\} \sqcup [J_+] \sqcup [J_-] \sqcup D,
\end{equation}
where
$$
D = \bigcup_{\lambda \in \CC-\left\{0, \pm 1\right\}} [D_\lambda] = \left\{A \in \SL{2}(\CC)\,|\, \tr A \neq 2\right\}.
$$

\begin{rmk}
The notation $[A]$ for the orbit of $A \in \SL{2}(\CC)$ is unfortunate since it collapses with the notation for virtual classes. However, we decided to use it to get in touch with the usual notation in the literature. We tried to use them in such a way that it should be possible to disambiguate which one we are referring to from the context.
\end{rmk}

The virtual classes of these strata can be easily computed. Here, we will sketch the computations. For more details, we refer to \cite{LMN}. Recall that $\coh{\SL{2}(\CC)} = q^3-q$, where we have shorten the Lefschetz motif as $q = \coh{\CC} \in \KVarc$. In order to check it, observe that the surjective regular map $\pi: \SL{2}(\CC) \to \CC^2-\left\{(0,0)\right\}$ given by $\pi(A) = A(1,0)$ defines a fibration
$
	\CC \to \SL{2}(\CC) \stackrel{\pi}{\to} \CC^2 - \left\{(0,0)\right\}
$.
Moreover, this fibration is locally trivial in the Zariski topology so, in particular, it is trivial in the Grothendieck ring. Hence, we have $\coh{\SL{2}(\CC)} = \coh{\CC} \cdot \coh{\CC^2- \left\{(0,0)\right\}} = q(q^2-1) = q^3-q$.

For the orbits of the Jordan type elements $J_\pm$, their stabilizers under the conjugacy action are $\Stab{J_\pm} \cong \CC$ so we have a fibration $
	\CC = \Stab{J_\pm} \to \left\{J_\pm \right\} \times \SL{2}(\CC) {\to} [J_\pm]$.
This fibration is locally trivial in the Zariski topology, so we have
$
	[J_\pm] = {\coh{\SL{2}(\CC)}}/{\coh{\CC}} = \frac{q^3-q}{q} = q^2-1$.
Analogously, for the orbit of $[D_\lambda]$ we have a locally trivial fibration $\CC^* \to \left\{D_\lambda\right\} \times \SL{2}(\CC) \to [D_\lambda]$ so $[D_\lambda] = (q^3-q)/(q-1) = q^2+q$.

For the computation of $[D]$, the monodromy representation of the trace map, $\tr: D \to \CC - \left\{\pm 2\right\}$, is the key. Any element of $D$ is given by a matrix
$$
	A=\begin{pmatrix}
	t/2 - a & b\\
	c & t/2+a\\
\end{pmatrix},
$$
with $a,b,c \in \CC$, $t \in \CC - \left\{\pm 2\right\}$ and satifying $t^2/4-a^2-bc=1$. With these coordinates, the trace map is $\tr(a,b,c,t)=t$. For $b \neq 0$, we can solve for $c$ and we get that $D \cap \left\{b \neq 0\right\} = \CC \times \CC^* \times \left(\CC-\left\{\pm 2\right\}\right)$, and the trace map is projection onto the last component so $\RM{D \cap \left\{b \neq 0\right\}}{\CC - \left\{\pm 2\right\}} = q(q-1)T_{\CC - \left\{\pm 2\right\}}$.

On the other hand, for $b = 0$ we have that $D \cap \left\{b=0\right\} = \CC \times \left\{t^2/4-a^2=1\right\}$, which is isomorphic to $\CC \times X_{2,-2}$ by the map $(c,a,t) \mapsto (c, t, 2a)$. Hence, $\RM{D \cap \left\{b = 0\right\}}{\CC - \left\{\pm 2\right\}} = qT_{\CC - \left\{\pm 2\right\}} + qS_2 \times S_{-2}$. Therefore, putting together both computations we get $\RM{D}{\CC - \left\{\pm 2\right\}} = q^2T_{\CC - \left\{\pm 2\right\}} + q S_2 \times S_{-2}$. In particular, using the projection map $c: \CC - \left\{\pm 2\right\} \to \star$, we obtain that $[D] = c_!\left(\RM{D}{\CC - \left\{\pm 2\right\}}\right) = q^2c_!(T_{\CC - \left\{\pm 2\right\}}) + q c_!(S_2 \times S_{-2}) = q^3-2q^2-q$.

Finally, with respect to the action by conjugation, observe that the GIT quotient is given by the trace map $\tr: \SL{2}(\CC) \to \SL{2}(\CC) \sslash \SL{2}(\CC) = \CC$. Moreover, the decomposition (\ref{eqn:dec-sl2}) corresponds, stratum by stratum, with the decomposition used in Example \ref{defn:piecewise-quotient} for the definition of the piecewise quotient. Hence, as $D \sslash \SL{2}(\CC) = \CC - \left\{ \pm 2\right\}$, we have that the virtual class of the quotient is
\begin{equation*}
	[\SL{2}(\CC) / \SL{2}(\CC)] = \left\{\Id\right\} + \left\{-\Id\right\} + \left\{[J_+]\right\} + \left\{[J_-]\right\} + \left[\CC-\left\{\pm 2\right\}\right].
\end{equation*}

\subsection{The core submodule}
\label{subsec:core-submodule}

Using the description of the geometric TQFT given in Section \ref{sec:geometric-reduced-TQFT}, we have obtained that
\begin{align*}
	\Zg{}(S^1, \star) &= \KVarrel{[\SL{2}(\CC) / \SL{2}(\CC)]} = \KVarrel{\left\{\Id\right\}} \oplus \KVarrel{\left\{-\Id\right\}} \\
	&\oplus \KVarrel{\left\{[J_+]\right\}} \oplus \KVarrel{\left\{[J_-]\right\}} \oplus \KVarrel{\CC - \left\{\pm 2\right\}}.
\end{align*}
The first four summands are relative varieties over a point and, thus, they are naturally isomorphic to $\KVarc$. For short, in the following we will denote $\Bt = \CC - \left\{\pm 2\right\}$. 

Let $T_1 \in \KVarrel{\left\{\Id\right\}}$, $T_{-1} \in \KVarrel{\left\{-\Id\right\}}$, $T_+ \in \KVarrel{\left\{[J_+]\right\}}$ and $T_- \in \KVarrel{\left\{[J_-]\right\}}$ be the units of these rings. On $\KVarrel{\Bt}$, we consider the unit $T_{\Bt}$ and $S_2, S_{-2} \in \KVarrel{\Bt}$ the elements described in Section \ref{sec:Hodge-monodromy-covering-spaces}. Pull them to $\KVarrel{[\SL{2}(\CC) / \SL{2}(\CC)]}$ via the natural inclusion. In this situation, the set
$$
\cS = \left\{T_1, T_{-1}, T_+, T_-, T_{\Bt}, S_2, S_{-2}, S_2 \times S_{-2}\right\} \subseteq \KVarrel{[\SL{2}(\CC) / \SL{2}(\CC)]}
$$
will be called the set of \emph{core elements}. The submodule of $\KVarrel{[\SL{2}(\CC) / \SL{2}(\CC)]}$ generated by $\cS$, $\cW$, is called the \emph{core submodule}.
\begin{comment}
\begin{defn}
The set $\cS = \left\{T_1, T_{-1}, T_+, T_-, T_{\Bt}, S_2, S_{-2}, S_2 \times S_{-2}\right\} \subseteq \KVarrel{[\SL{2}(\CC) / \SL{2}(\CC)]}$ is called the set of \emph{core elements}. The submodule of $\KVarrel{[\SL{2}(\CC) / \SL{2}(\CC)]}$ generated by $\cS$, $\cW$, is called the \emph{core submodule}.
\end{defn}
\end{comment}
The submodule $\cW$ will be very important in the incoming computations since we will show that $\cW$ is the submodule generated by the elements $\Zg{}(L)^g \circ \Zg{}(D)(\star)$ for $g \geq 0$, that is $\cW = \cV_{(S^1, \star)}$ using the notations of Section \ref{sec:reduction-TQFT}.

Now, consider the quotient map $\tr: \SL{2}(\CC) \to [\SL{2}(\CC)/\SL{2}(\CC)]$. From this map, we can build the endomorphism $\eta = \tr_! \circ \tr^*: \KVarrel{[\SL{2}(\CC)/\SL{2}(\CC)]} \to \KVarrel{[\SL{2}(\CC)/\SL{2}(\CC)]}$. This morphism will be useful in the upcoming sections due to its role in Proposition \ref{prop:existence-reduction-TQFT}.

\begin{prop}\label{prop:eta-morph}
The submodule $\cW$ is invariant for the morphism $\eta$. Actually, we have
$$
\begin{tabular}{ccc}
	$\eta(T_{\pm 1}) = T_{\pm 1}$, &\hspace{0.5cm}\null& $\eta(T_{\pm}) = (q^2-1)T_{\pm}$,\\
	$\eta(T_{\Bt}) = q^2 T_{\Bt} + q S_2 \times S_{-2}$, && $\eta(S_2 \times S_{-2}) = q T_{\Bt} + q^2 S_2 \times S_{-2}$,\\
	$\eta(S_2) = q^2 S_2 + q S_{-2}$, && $\eta(S_{-2}) = q S_2 + q^2 S_{-2}$.
\end{tabular}
$$
\begin{proof}
First of all, observe that, if $i: X \hookrightarrow [\SL{2}(\CC) / \SL{2}(\CC)]$ is the inclusion of a locally closed subset and $\bar{X} = \tr^{-1}(X) \subseteq \SL{2}(\CC)$, they fit in a pullback diagram
\[
\begin{displaystyle}
   \xymatrix
   {	
	\bar{X} \ar[r]^{\tr|_{\bar{X}}} \ar@{^{(}-{>}}[d]& X \ar@{^{(}-{>}}[d]^i \\
	\SL{2}(\CC) \ar[r]_{\tr\;\;\;\;}& [\SL{2}(\CC)/\SL{2}(\CC)]
   }
\end{displaystyle}   
\]
Thus, $\tr_! \circ \tr^* \circ i_! = (\tr|_{\bar{X}})_! \circ (\tr|_{\bar{X}})^*$ that is, we can compute the image of $\eta$ in $\bar{X}$. For this reason, since $\tr|_{\left\{\pm \Id\right\}}: \left\{\pm \Id\right\} \to \left\{\pm \Id\right\}$ is the identity map, we have $\eta(T_{\pm 1}) = T_{\pm 1}$. Analogously, since $(\tr|_{[J_\pm]})^*$ is a ring homomorphism, it sends units into units so, for $T_{\pm}$, we have that $\eta(T_{\pm}) = (\tr|_{[J_\pm]})_!T_\pm = [J_+] = (q^2-1)T_{\pm}$.

For the first identity of the second row, we have $(\tr|_D)_!(\tr|_D)^*T_{\Bt} = (\tr|_D)_!T_D = \RM{D}{\Bt}$ and the result follows from the computations of the previous section. For $S_2$, consider the auxiliar variety $X_{2} = \left\{(t,y) \in \Bt \times \CC^*\,|\,y^2 = t-2\right\}$, for which $\RM{X_{2}}{\Bt} = T_{\Bt} + S_2$ (see Section \ref{sec:Hodge-monodromy-covering-spaces}), and the variety $X_{2}' = \left\{(A,y) \in D \times \CC^*\,|\,y^2 = \tr(A)-2\right\}$. They fit in a commutative diagram
	\[
\begin{displaystyle}
   \xymatrix
   {	
   		& X_{2}'\ar[d] \ar[r] \ar[dl]& X_{2} \ar[d] \\
   		\Bt & D \ar[r]_{\tr|_D} \ar[l]^{\;\;\tr|_D} & \Bt
   }
\end{displaystyle}   
\]
where the rightmost square is a pullback. Hence, we have that $\eta(T_{\Bt}) + \eta(S_2) = (\tr|_D)_! (\tr|_D)^*(T_{\Bt} + S_2) = \RM{X_{2}'}{\Bt}$ and, thus, $\eta(S_2) = \RM{X_{2}'}{\Bt} - q^2T_{\Bt}-qS_2 \times S_{-2}$ by the computation for $T_{\Bt}$.

In order to compute $\RM{X_{2}'}{\Bt}$, observe that writing an element $A \in X_{2}'$ with coordinates
$$
A = \begin{pmatrix}
	t/2-a & b\\
	c & t/2+a\\
\end{pmatrix},
$$
we get that $
	X_{2}' = \left\{(a,b,c,t,y) \in \CC^5\,|\, y^2=t-2, t^2/4-a^2-bc=1, t \neq \pm 2\right\}$. In this way, the projection map becomes projection over $t$. For $b \neq 0$, we can solve for $c$ and get that $X_{2}' \cap \left\{c \neq 0\right\} = \CC \times \CC^* \times X_2$, so $\RM{X_{2}' \cap \left\{c \neq 0\right\}}{\Bt} = q(q-1)\RM{X_2}{\Bt} = q(q-1)(T_{\Bt} + S_2)$. On the other hand, we have that
\begin{align*}
	X_{2}' \cap \left\{c = 0\right\} &= \CC \times \left\{(a,t,y) \in \CC^3\,|\, y^2=t-2, 4a^2= t^2-4, t \neq \pm 2\right\} \\
	& \cong \CC \times \left\{(a',t,y) \in \CC^3\,|\, y^2=t-2, (a')^2= t+2, t \neq \pm 2\right\} = \CC \times \left(X_2 \times_{\Bt} X_{-2}\right),
\end{align*}
where the last isomorphism is given by setting $a' = 2a/y$. Hence, $\RM{X_{2}' \cap \left\{c = 0\right\}}{\Bt} = q\left(T_{\Bt} + S_2 + S_{-2} + S_2 \times S_{-2}\right)$. Putting both computations together, we find
$$
	\RM{X_{2}'}{\Bt} = \RM{X_{2}' \cap \left\{c \neq 0\right\}}{\Bt}+\RM{X_{2}' \cap \left\{c = 0\right\}}{\Bt} = q^2(T_{\Bt} + S_2)+q(S_{-2} + S_{2} \times S_{-2}).
$$
From this, it follows that $\eta(S_2) = q^2S_2 + qS_{-2}$, as claimed. The calculations for $S_{-2}$ and $S_2 \times S_{-2}$ are analogous.
\end{proof}
\end{prop}

The morphism $\eta: \cW \to \cW$ is not invertible since, on $\KVarc$, the elements $q$, $q-1$ and $q+1$ have no inverse. We can solve this problem by considering the localization of $\KVarc$ over multiplicative system generated by $q$, $q-1$ and $q+1$, denoted $\Ko{\Varc}$. Extending the localization to $\cW$, we obtain the $\Ko{\Varc}$-module $\cW_0$. In that module, by the proposition above, we have an isomorphism $\eta: \cW_0 \to \cW_0$.

Let us consider the parabolic data $\Lambda = \left\{\left\{\Id\right\}, \left\{-\Id\right\}, [J_+], [J_-]\right\}$, let $\Zs{}: \Bordppar{n}{\Lambda} \to \Modt{\KVarc}$ be the associated standard TQFT of Section \ref{sec:TQFT-for-repr} and let $\Zg{}: \Bordppar{n}{\Lambda} \to \Modt{\KVarc}$ be the geometric TQFT of Section \ref{sec:geometric-reduced-TQFT}. As a result of the upcoming computations of Sections \ref{sec:tube-J+} and \ref{sec:genus-tube}, we will show that $\cW = \cV_{(S^1, \star)}$. Hence, applying Proposition \ref{prop:existence-reduction-TQFT}, we obtain the following result.

\begin{thm}
Let $\Zs{}: \Bordppar{n}{\Lambda} \to \Modt{\Ko{\Var{k}}}$ be the standard TQFT for $\SL{2}(\CC)$-representation varieties and let $\Zg{}: \Bordppar{n}{\Lambda} \to \Modt{\Ko{\Var{k}}}$ be the geometric TQFT. There exists an almost-TQFT, $\cZg{}: \Tub{2}(\Lambda) \to \Mod{\Ko{\Var{k}}}$, such that:
\begin{itemize}
	\item The image of the circle is $\cZg{}(S^1, \star) = \cW \subseteq \Ko{\Varrel{[\SL{2}(\CC)/\SL{2}(\CC)]}}$.
	\item For an strict tube $(W, A, Q): (S^1, \star) \to (S^1, \star)$ it assigns the morphism $\cZg{}(W, A, Q) = \Zg{}(W, A, Q) \circ \eta^{-1}: \cW_0 \to \cW_0$.
	\item For a closed surface $W: \emptyset \to \emptyset$, it gives $
		\cZg{}(W, A, Q)(\star) = \coh{\Rep{\SL{2}(\CC)}(W, Q)} \times (q^3-q)^{|A|-1}$.
\end{itemize}
\end{thm}

\subsection{Discs and first tube}
\label{sec:discs-tubes}

As a warm lap, in this section we will compute the morphisms $\Zg{}(L_{-\Id})$, $\Zg{\SL{2}(\CC)}(D)$ and $\Zg{}(D^\dag)$.
For the first case the image under the field theory of $L_{-\Id}$ is the span
$$
\begin{matrix}
	[\SL{2}(\CC)/\SL{2}(\CC)] & \stackrel{\tr \circ \pi_1 }{\longleftarrow} & \SL{2}(\CC)^2 & \stackrel{-\tr \circ \pi_1}{\longrightarrow} & [\SL{2}(\CC)/\SL{2}(\CC)]
\end{matrix},
$$
where $\pi_1: \SL{2}(\CC)^2 \to \SL{2}(\CC)$ is the projection onto the first variable.
Observe that $(\pi_1)_! \circ (\pi_1)^*: \KVarrel{\SL{2}(\CC)} \to \KVarrel{\SL{2}(\CC)}$ is just multiplication by $\coh{\SL{2}(\CC)} = q^3-q$. Hence, we have $\Zg{}(L_{-\Id}) = (-\tr)_! \circ (\pi_1)_! \circ (\pi_1)^* \circ \tr^* = (q^3-q)\sigma_! \circ \eta$, where $\sigma: [\SL{2}(\CC)/\SL{2}(\CC)] \to [\SL{2}(\CC)/\SL{2}(\CC)]$ is the reflection that sends $\sigma([A]) = [-A]$ for $A \in \SL{2}(\CC)$. It is straighforward to check that
$$
	\sigma_! (T_{\pm 1}) = T_{\mp 1} \hspace{0.9cm} \sigma_! (T_{\pm}) = T_{\mp} \hspace{0.9cm} \sigma_! (T_{\Bt}) = T_{\Bt}
$$
$$\sigma_! (S_{\pm 2}) = S_{\mp 2} \hspace{0.9cm} \sigma_! (S_2 \times S_{-2}) = S_2 \times S_{-2}. 
$$
In particular, $\sigma_!(\cW) \subseteq \cW$ so, by Proposition \ref{prop:eta-morph}, we have that $\Zg{}(L_{-\Id})(\cW) \subseteq \cW$. Hence, $\cZg{}(L_{-\Id}) = \Zg{}(L_{-\Id}) \circ \eta^{-1} = (q^3-q)\sigma_!$. Thus, the matrix of $\cZg{}(L_{-\Id})$ in the generators $\cS$ of $\cW_0$ is
$$
	\cZg{}(L_{-\Id}) = (q^3-q)\begin{pmatrix}
	0 & 1 & 0 & 0 & 0 & 0 & 0 & 0\\
	1 & 0 & 0 & 0 & 0 & 0 & 0 & 0\\
	0 & 0 & 0 & 1 & 0 & 0 & 0 & 0\\
	0 & 0 & 1 & 0 & 0 & 0 & 0 & 0\\
	0 & 0 & 0 & 0 & 1 & 0 & 0 & 0\\
	0 & 0 & 0 & 0 & 0 & 0 & 1 & 0\\
	0 & 0 & 0 & 0 & 0 & 1 & 0 & 0\\
	0 & 0 & 0 & 0 & 0 & 0 & 0 & 1\\
	\end{pmatrix}.
$$

Finally, with respect to the discs $D: \emptyset \to (S^1,\star)$ and $D^\dag: (S^1,\star) \to \emptyset$ recall that $\Zg{}(D) = i_!$ and $\Zg{}(D^\dag) = i^*$, where $i: \left\{\Id\right\} \hookrightarrow [\SL{2}(\CC)/\SL{2}(\CC)]$ is the inclusion. The image of the former map is very easy to identify since, by definition $\Zg{}(D)(\star) = i_!(\star)=T_1$. For the later one, observe that $i^*T_1=\star$ but $i^*T_{-1} = i^*T_\pm = i^*T_{\Bt} = i^*S_{\pm 2} = 0$, so $\Zg{}(D^\dag): \cW \to \KVar{}$ is just the projection onto $T_1$. Since $\eta$ fixes $T_1$, we also have $\Zg{}(D)= \cZg{}(D)$ and $\Zg{}(D^\dag)=\cZg{}(D^\dag)$.

\subsection{The tube with a Jordan type marked point}
\label{sec:tube-J+}

In this section, we will focus on the computation of the image of the tube $L_{[J_+]}: (S^1, \star) \to (S^1, \star)$. Again, the field theory for $\Zg{}$ on $L_{[J_+]}$ is the span
$$
\begin{matrix}
	[\SL{2}(\CC)/\SL{2}(\CC)] & \stackrel{\hat{r}}{\longleftarrow} & \SL{2}(\CC)^2 \times [J_+] & \stackrel{\hat{s}}{\longrightarrow} & [\SL{2}(\CC)/\SL{2}(\CC)]\\
	\tr(A) & \mapsfrom & (A, B, C) & \mapsto & \tr(BACB^{-1}) = \tr(AC)
\end{matrix}
$$
so $\Zg{}(L_{[J_+]}) = \hat{s}_! \circ \hat{r}^*: \KVarrel{[\SL{2}(\CC)/\SL{2}(\CC)]} \to \KVarrel{[\SL{2}(\CC)/\SL{2}(\CC)]}$.

Notice that, if $i: X \hookrightarrow [\SL{2}(\CC) / \SL{2}(\CC)]$ is an inclusion of a locally closed subset, then we have a commutative diagram
\[
\begin{displaystyle}
   \xymatrix
   {	
	\hat{r}^{-1}(X) \ar[r]^{\;\;\;\hat{r}|_{\hat{r}^{-1}(X)}} \ar@{^{(}-{>}}[d]& X \ar@{^{(}-{>}}[d]^i \\
	\SL{2}(\CC)^2 \times [J_+] \ar[r]_{\hat{r}\;}& [\SL{2}(\CC)/\SL{2}(\CC)]
   }
\end{displaystyle}   
\]
Hence, for any $S \in \KVarrel{X}$, $\hat{s}_!\circ \hat{r}^* \circ i_!S = \hat{s}_!\circ (\hat{r}|_{\hat{r}^{-1}(X)})^*S$. Moreover, if we decompose $[\SL{2}(\CC)/\SL{2}(\CC)] = \sum_k Y_k$ with inclusions $j_k: Y_k \hookrightarrow [\SL{2}(\CC)/\SL{2}(\CC)]$ then $\hat{s}_! = \sum_k (j_k)_!(j_k)^*\hat{s}_! = \sum_k (j_k)_!(\hat{s}|_{\hat{s}^{-1}(Y_k)})_!$, since $\sum_k (j_k)_!(j_k)^* = 1$. Therefore,
$$
	\Zg{}(L_{[J_+]}) (i_!S) = \hat{s}_!\circ \hat{r}^* \circ i_!S = \sum_k (j_k)_! \circ (\hat{s}|_{\hat{r}^{-1}(X) \cap \hat{s}^{-1}(Y_k)})_! \circ (\hat{r}|_{\hat{r}^{-1}(X) \cap \hat{s}^{-1}(Y_k)})^* S
$$
Notice that we can decompose
\begin{align*}
	\hat{r}^{-1}(X) \cap \hat{s}^{-1}(Y) &= \left\{(A, B, C) \in \hat{r}^{-1}(X)  \,\,|\,\, \tr(AC) \in Y\right\} \\
	&\cong \left\{A \in \tr^{-1}(X) \,\,|\,\, \tr(AJ_+) \in Y\right\} \times \left(\SL{2}(\CC)/\CC\right) \times \SL{2}(\CC),
\end{align*}
where the last isomorphism is given by sending $(A,P,B) \in \tr^{-1}(X) \times \left(\SL{2}(\CC)/\CC\right) \times \SL{2}(\CC)$, in the variety downstairs, to $(PAP^{-1}, B, PJ_+P^{-1})$, in the variety upstairs. For short, we shall denote $Z_{X,Y} = \hat{r}^{-1}(X) \cap \hat{s}^{-1}(Y)$ and $Z_{X,Y}^0 = \left\{A \in \tr^{-1}(X) \,\,|\,\, \tr(AJ_+) \in Y\right\}$.
Using as stratification the decomposition $[\SL{2}(\CC)/\SL{2}(\CC)] = \left\{\pm \Id\right\} + \left\{[J_\pm]\right\} + \Bt$, we can compute:
\begin{itemize}
	\item For $T_{1}$, we have $\hat{r}^{-1}(\Id) = Z_{\Id, [J_+]} = [J_+] \times \SL{2}(\CC)$ and $\hat{s}$ is just the projection onto a point $\hat{s}: [J_+] \times \SL{2}(\CC) \to \left\{[J_+]\right\}$. Hence, $\Zg{}(L_{[J_+]})(T_1) = (q^2-1)(q^3-q) T_+$.
	\item For $T_{-1}$, since $[-J_+] = [J_-]$, we have $\hat{r}^{-1}(-\Id) = Z_{-\Id, [J_-]} = [J_-] \times \SL{2}(\CC)$. Again, $\hat{s}$ is just the projection $[J_-] \times \SL{2}(\CC) \to \left\{[J_-]\right\}$ so $\Zg{}(L_{[J_+]})(T_{-1}) = (q^2-1)(q^3-q) T_-$.
	\item For $T_+$, the situation becomes more involved since we have a non-trivial decomposition of $\hat{r}^{-1}([J_+])$. We analize each stratum separately:
	\begin{itemize}
		\item For $\Id$, we have $Z_{[J_+], \Id} = [J_+] \times \SL{2}(\CC)$ since all the elements of $Z_{[J_+], \Id}$ have the form $(A^{-1}, B, A)$. Thus, $\hat{s}: [J_+]  \times \SL{2}(\CC) \to \left\{\Id\right\}$ is the projection onto a point so $\Zg{}(L_{[J_+]})(T_{+})|_{\Id} = (q^2-1)(q^3-q)T_1$.
		\item For $-\Id$, observe that $Z_{[J_+], -\Id} = \emptyset$ so $\Zg{}(L_{[J_+]})(T_{+})|_{-\Id} = 0$.
		\item For $[J_+]$, a straightforward check shows that
		$$
			Z_{[J_+], [J_+]}^0 = \left\{\left.\begin{pmatrix}
	1 & b\\
	0 & 1\\
\end{pmatrix}\,\right|\,\, b \in \CC-\left\{0, -1\right\}
\right\}.
		$$
	Hence, $Z_{[J_+], [J_+]} = \left(\CC-\left\{0, -1\right\}\right) \times \left(\SL{2}(\CC)/\CC\right)  \times \SL{2}(\CC)$ and $\hat{s}$ is a projection onto the singleton $\left\{[J_+]\right\}$. Thus, $\Zg{}(L_{[J_+]})(T_{+})|_{[J_+]} = (q-2)(q^2-1)(q^3-q)T_+$.
		\item For $[J_-]$, we have that 
		$$
			Z_{[J_+], [J_-]}^0 = \left\{\left.\begin{pmatrix}
	1+a & \frac{a^2}{2}\\
	-4 & 1-a\\
\end{pmatrix}\,\right|\,\, a \in \CC
\right\}.
		$$
		Therefore, $Z_{[J_+], [J_-]} = \CC \times \left(\SL{2}(\CC)/\CC\right)  \times \SL{2}(\CC)$ so $\Zg{}(L_{[J_+]})(T_{+})|_{[J_-]} = q(q^2-1)(q^3-q)T_-$.
		\item For $\Bt$, we have that
$$
		Z_{[J_+], \Bt}^0 = \left\{\left.\begin{pmatrix}
	1+a & -\frac{a^2}{c}\\
	c & 1-a\\
\end{pmatrix}\,\right|\,\, a \in \CC, c \in \CC^* - \left\{-4\right\}
\right\},
		$$
with the projection $Z_{[J_+], \Bt}^0 \to \Bt$ given by $A \mapsto c + 2$. Hence $\hat{s}: Z_{[J_+], \Bt} = Z_{[J_+], \Bt}^0 \times \left(\SL{2}(\CC)/\CC\right)  \times \SL{2}(\CC) \to \Bt$ applies $\hat{s}(A, B, P) = c + 2$. In particular, $\hat{s}$ is trivial, so $\Zg{}(L_{[J_+]})(T_+)|_{\Bt} = \RM{Z_{[J_+], \Bt}}{\Bt} = q(q^2-1)(q^3-q)T_{\Bt}$.
	\end{itemize}
Summarizing, we have
$$
	\Zg{}(L_{[J_+]})(T_{+}) = (q^2-1)(q^3-q)\left(T_1 + (q-2)T_+ + qT_- + qT_{\Bt}\right).
$$
	\item For $T_-$, the calculation is very similar to the one of $T_+$, since $Z_{[J_-], Y}^0 = Z_{[J_+], -Y}^0$ for every stratum $Y$. Hence, we have:
		\begin{itemize}
		\item For $\Id$, $Z_{[J_-], \Id} = \emptyset$, so $\Zg{}(L_{[J_+]})(T_{-})|_{\Id} = 0$. Analogously, for $-\Id$ we have $Z_{[J_-], -\Id} = -[J_+] \times \SL{2}(\CC)$, so $\Zg{}(L_{[J_+]})(T_{-})|_{-\Id} = (q^2-1)(q^3-q)T_{-1}$.
		\item For $[J_+]$, $Z_{[J_-], [J_+]}^0 = Z_{[J_+], [J_-]}^0 \cong \CC$, so $\Zg{}(L_{[J_+]})(T_{-})|_{[J_+]} = q(q^2-1)(q^3-q)T_+$.
		\item For $[J_-]$, we have $Z_{[J_-], [J_-]}^0 = Z_{[J_+], [J_+]}^0 \cong \CC - \left\{0, -1\right\}$. Hence, $\Zg{}(L_{[J_+]})(T_{-})|_{[J_-]} = (q-2)(q^2-1)(q^3-q)T_-$.
		\item For $\Bt$, since $-\Bt=\Bt$, we have $Z_{[J_-], \Bt}^0 = Z_{[J_+], -\Bt}^0 \cong \CC \times \left(\CC^*-\left\{-4\right\}\right)$. Hence, $\hat{s}$ is also trivial on $Z_{[J_-], \Bt}$ so $\Zg{}(L_{[J_+]})(T_-)|_{\Bt} = q(q^2-1)(q^3-q)T_{\Bt}$.
	\end{itemize}
Summarizing, this calculation says that
$$
	\Zg{}(L_{[J_+]})(T_{-}) = (q^2-1)(q^3-q)\left(T_{-1} + qT_+ + (q-2)T_- + qT_{\Bt}\right).
$$
	\item For $T_{\Bt}$, observe that we have $Z_{\Bt, \pm \Id} = \emptyset$ so $\Zg{}(L_{[J_+]})(T_{\Bt})$ has no components on $\pm \Id$. In this section, we will compute the image of $\Zg{}(L_{[J_+]})(T_{\Bt})$ on $[J_\pm]$. Its image on $\Bt$ is much harder and it will be described later. For $[J_\pm]$ observe that $Z_{\Bt, [J_\pm]}^0 \cong Z_{[J_\pm], \Bt}^0 \cong \CC \times (\CC^* - \left\{-4\right\})$. Hence, $\Zg{}(L_{[J_+]})(T_{\Bt})|_{[J_\pm]} = q(q-2)(q^2-1)(q^3-q)T_\pm$.
	\item For $S_2$, $\Zg{}(L_{[J_+]})(S_2)$ has no components on $\pm \Id$. Again, we will only focus on its image on $[J_\pm]$. We consider the auxiliar variety
$$
	X_2 = \left\{(t,y) \in \Bt \times \CC^*\,\,|\,\, y^2 = t-2\right\},
$$
with the fibration $X_2 \to \Bt$ given by $(t,y) \mapsto t$ whose monodromy representation is $\RM{X_2}{\Bt} = T_{\Bt} + S_2$. In that case, we have a commutative diagram
\[
\begin{displaystyle}
   \xymatrix
   {	
	& \hat{X}_2\ar[d]\ar[ld] \ar[r]& X_2 \ar[d]\\
	[J_+] & Z_{\Bt, [J_+]} \ar[l]^{\hat{s}} \ar[r]_{\;\;\;\;\hat{r}} & \Bt 
   }
\end{displaystyle}   
\]
where $\hat{X}_2 = X_2 \times_{\Bt} Z_{\Bt, [J_+]}$. Computing directly, we obtain that $
	\hat{X}_2 \cong X_2 \times \CC \times \left(\SL{2}(\CC)/\CC\right) \times \SL{2}(\CC)$. Hence, $\RM{\hat{X}_2}{[J_+]} = q(q^2-1)(q^3-q)\coh{X_2}T_+ = q(q^2-1)(q^3-q)(q-3)T_+$. Since
$$
	\RM{\hat{X}_2}{[J_+]} = \hat{s}_!\hat{r}^* \RM{X_2}{\Bt} = \Zg{}(L_{[J_+]})(T_{\Bt})|_{[J_+]} + \hat{s}_!\hat{r}^* S_2,
$$
	we obtain that $\Zg{}(L_{[J_+]})(S_2)|_{[J_+]} =\hat{s}_!\hat{r}^* S_2 = -q(q^2-1)(q^3-q)T_+$. Furthermore, as $Z_{\Bt,[J_+]} \cong Z_{\Bt,[J_-]}$, we also have $\Zg{}(L_{[J_+]})(S_2)|_{[J_-]} = -q(q^2-1)(q^3-q)T_-$.
\item For $S_{-2}$ the calculation is analogous to the one of $S_2$ so $\Zg{}(L_{[J_+]})(S_{-2})|_{[J_+]} = -q(q^2-1)(q^3-q)T_+$ and $\Zg{}(L_{[J_+]})(S_{-2})|_{[J_-]} = -q(q^2-1)(q^3-q)T_-$. The same holds for $S_2 \times S_{-2}$ so $\Zg{}(L_{[J_+]})(S_2 \times S_{-2})|_{[J_+]} = -q(q^2-1)(q^3-q)T_+$ and $\Zg{}(L_{[J_+]})(S_2 \times S_{-2})|_{[J_-]} = -q(q^2-1)(q^3-q)T_-$.
\end{itemize}

\subsubsection*{The components in $\Bt$}

The calculation of the components of the image of $T_{\Bt}, S_2, S_{-2}$ and $S_2 \times S_{-2}$ on $\Bt$ is harder than the previous ones and it requires a subtler analysis. For this reason, we compute it separately. First of all, observe that
$$
	Z_{\Bt, \Bt}^0 = \left\{\left.\begin{pmatrix}
	\frac{t}{2} + a & b\\
	t'-t & \frac{t}{2} - a\\
\end{pmatrix}\,\right|\,\, \begin{matrix}
	t, t' \in \Bt, a,b \in \CC\\
	\frac{t^2}{4}-a^2-(t'-t)b = 1
\end{matrix}
\right\}.
$$
The projections onto each copy of $\Bt$ are $\hat{r}(a, b, t, t') = t$ and $\hat{s}(a, b, t, t') = t'$. We decompose $Z_{\Bt, \Bt}^0 = Z_1 \sqcup Z_2$ where $Z_1 = Z_{\Bt, \Bt}^0 \cap \left\{t=t'\right\}$ and $Z_2 = Z_{\Bt, \Bt}^0 \cap \left\{t \neq t'\right\}$.
For $Z_1$, we have the explicit expression $Z_1 = \left\{4a^2 = (t+2)(t-2)\right\} \times \CC$, so $\RM{Z_1}{\Bt} = q(T_{\Bt} + S_2 \times S_{-2})$. Thus, the contribution of $Z_1$ to $\Zg{}(L_{[J_+]})(T_{\Bt})|_{\Bt}$ is $q(q^2-1)(q^3-q)(T_{\Bt} + S_2 \times S_{-2})$.
	
	With respect to $S_2$, we again consider $X_2 = \left\{y^2 = t-2\right\}$, for which we have a commutative diagram
\[
\begin{displaystyle}
   \xymatrix
   {	
	& \hat{X}_2\ar[d]\ar[ld] \ar[r]& X_2 \ar[d]\\
	\Bt & Z_1 \ar[l]^{\hat{s}} \ar[r]_{\hat{r}} & \Bt 
   }
\end{displaystyle}   
\]
where $\hat{X}_2 = X_2 \times_{\Bt} Z_1$ is the pullback of $Z_1$ and $X_2$ over $\Bt$. Computing directly
$$
	\hat{X}_2 = \left\{(t, a, y) \in \Bt \times \CC^* \times \CC \times \CC^* \,\left|\, \begin{matrix}4a^2=(t+2)(t-2)\\y^2=t-2\end{matrix}\right.\right\} \times \CC,
$$
with projection $\hat{X}_2 \to \Bt$ given by $(t, a, y, b) \mapsto t$. This is isomorphic to $X_2 \times_{\Bt} X_{-2} \to \Bt$, so we have $\RM{\hat{X}_2}{\Bt} = q(T_{\Bt} + S_2 + S_{-2} + S_2 \times S_{-2})$. Since $\RM{\hat{X}_2}{\Bt} = \RM{Z_1}{\Bt} + \hat{s}_!\hat{r}^* S_2$, we finally obtain that $\hat{s}_!\hat{r}^*S_2 = q(S_2 + S_{-2})$. Therefore, the contribution of $Z_1$ to $\Zg{}(L_{[J_+]})(S_2)|_{\Bt}$ is $q(q^2-1)(q^3-q)(S_2 + S_{-2})$. By symmetry, the contribution to $\Zg{}(L_{[J_+]})(S_{-2})|_{\Bt}$ is $q(q^2-1)(q^3-q)(S_2 + S_{-2})$. An analogous computation for $S_2 \times S_{-2}$ shows that the contribution of $Z_1$ to $\Zg{}(L_{[J_+]})(S_2 \times S_{-2})|_{\Bt}$ is $q(q^2-1)(q^3-q)(T_{\Bt} + S_2 \times S_{-2})$.

The most involved stratum is $Z_2$. For it, we have the explicit expression
$$
	Z_2 = \left\{\left.\begin{pmatrix}
	\frac{t}{2} + a & \frac{t^2 - 4a^2 -4}{4(t'-t)}\\
	t'-t & \frac{t}{2} - a\\
\end{pmatrix}\,\right|\,\, \begin{matrix}
	t, t' \in \Bt, t' \neq t\\
	a \in \CC
\end{matrix}
\right\} \cong \CC \times \left[(\CC-\left\{\pm 2\right\})^2- \Delta\right],
$$
where $\Delta \subseteq (\CC-\left\{\pm 2\right\})^2$ is the diagonal. The projections are $\hat{r}(a,t,t')=t$ and $\hat{s}(a,t,t')=t'$, so $\hat{s}$ is a trivial fibration. Hence, $\RM{Z_2}{\Bt} = q(q-3)T_{\Bt}$ and the contribution of $Z_2$ to $\Zg{}(L_{[J_+]})(T_{\Bt})|_{\Bt}$ is $q(q-3)(q^2-1)(q^3-q)T_{\Bt}$.

For $S_2$, again consider the auxiliar variety $X_2 = \left\{y^2 = t-2\right\}$ for $t \neq \pm 2$. In that case, the pullback of $X_2$ and $Z_2$ over $\Bt$ is
\begin{align*}
	\hat{X}_2 = \CC \times \left\{(t, t', y) \in \left[(\CC-\left\{\pm 2\right\})^2- \Delta\right] \times \CC^* \,\left|\, y^2=t-2\right.\right\}.
\end{align*}
If we ignore the condition $t \neq t'$, we obtain a variety $\overline{X}_2$ which is a trivial fibration over $\Bt$. Its fiber is a line times an affine parabola with three points removed, so $\RM{\overline{X}_2}{\Bt} = (q-3)qT_{\Bt}$.

Now, the difference $\overline{X}_2 - \hat{X}_{2}$ is a trivial line bundle times a double covering over $\Bt$. The fiber of the covering over $t' \in \Bt$ are the two points $\left\{(t',t', \pm \sqrt{t'-2})\right\}$, so $\RM{\overline{X}_2 - \hat{X}_{2}}{\Bt} = q(T_{\Bt} + S_2)$. Hence, $\RM{\hat{X}_2}{\Bt} = q(q-4)T_{\Bt} - q S_2$. Substracting the term from $T_{\Bt}$, we get that the contribution of $Z_2$ to $\Zg{}(L_{[J_+]})(S_2)|_{\Bt}$ is $-q(q^2-1)(q^3-q)\left(T_{\Bt}+S_2\right)$. For $S_{-2}$ and $S_2 \times S_{-2}$, the calculations are analogous and we get a contribution of $-q(q^2-1)(q^3-q)\left(T_{\Bt}+S_{-2}\right)$ and $-q(q^2-1)(q^3-q)\left(T_{\Bt}+S_2 \times S_{-2}\right)$, respectively.

Summarizing, if we put together the calculations for $Z_1$ and $Z_2$, we finally obtain that
\begin{align*}
 \Zg{}(L_{[J_+]})(T_{\Bt})|_{\Bt} & = (q^2-1)(q^3-q)\left(q (T_{\Bt} + S_2 \times S_{-2}) + q(q-3) T_{\Bt}\right) \\
& = (q^2-1)(q^3-q)\left((q^2-2q)T_{\Bt} + qS_2 \times S_{-2}\right),\\
\Zg{}(L_{[J_+]})(S_2)|_{\Bt} &= (q^2-1)(q^3-q)\left(q(S_2 + S_{-2}) - q(T_{\Bt} + S_2)\right) \\
&= (q^2-1)(q^3-q)\left(-qT_{\Bt} + qS_2\right),\\
\Zg{}(L_{[J_+]})(S_{-2})|_{\Bt} &= (q^2-1)(q^3-q)\left(q (S_2 + S_{-2}) - q(T_{\Bt} + S_{-2})\right) \\
&= (q^2-1)(q^3-q)\left(-qT_{\Bt} + qS_{-2}\right),
\\
\Zg{}(L_{[J_+]})(S_2 \times S_{-2})|_{\Bt} &= (q^2-1)(q^3-q)\left(q (T_{\Bt} + S_2 \times S_{-2}) - q(T_{\Bt} + S_2 \times S_{-2})\right) = 0.
\end{align*}

In particular, we have shown that the core submodule, $\cW$, is invariant under the morphism $\Zg{}(L_{[J_+]}): \KVarrel{[\SL{2}(\CC)/\SL{2}(\CC)]} \to \KVarrel{[\SL{2}(\CC)/\SL{2}(\CC)]}$. Hence, multiplying by $\eta^{-1}$ on the right and using the previous expression, we obtain the following theorem.

\begin{thm}
In the set of generators $\cS$, the matrix of $\cZg{}(L_{[J_+]}): \cW_0 \to \cW_0$ is
$$
\cZg{}(L_{[J_+]}) = (q^3-q)\begin{pmatrix}
0 & 0 & 1 & 0 & 0 & 0 & 0 & 0 \\
0 & 0 & 0 & 1 & 0 & 0 & 0 & 0 \\
q^{2} - 1 & 0 & q - 2 & q & (q-1)^2 & -q
+ 1 & -q + 1 & -2 \, q + 2 \\
0 & q^{2} - 1 & q & q - 2 & (q-1)^2& -q
+ 1 & -q + 1 & -2 \, q + 2 \\
0 & 0 & q & q & q^{2} - 2 \, q & -q + 1 & -q + 1
& -q + 2 \\
0 & 0 & 0 & 0 & 0 & -1 & q & 0 \\
0 & 0 & 0 & 0 & 0 & q & -1 & 0 \\
0 & 0 & 0 & 0 & q & 0 & 0 & -1
\end{pmatrix}
$$
\end{thm}

From this computation, we can calculate the image of the tube $L_{[J_-]}: (S^1, \star) \to (S^1, \star)$. The point is that, since $[J_-] = -[J_+]$, the span associated to $L_{[J_-]}$ fits in a commutative diagram
\[
\begin{displaystyle}
   \xymatrix
   {
   	& \SL{2}(\CC)^2 \times [J_-] \ar[d]^{\varphi} \ar[rd]\ar[ld] & \\
	\SL{2}(\CC) &\SL{2}(\CC)^2 \times [J_+] \ar[r]_{\;\;\;\;\;\;\;\;r} \ar[l]^{\varsigma \circ s\;\;\;\;\;\;\;} & \SL{2}(\CC) \\
   }
\end{displaystyle}   
\]
where $\varphi(A,B,C) = (A,B,-C)$ and $\varsigma: \SL{2}(\CC) \to \SL{2}(\CC)$ is given by $\varsigma(A)=-A$. Hence, as $\varphi$ is an isomorphism we have $\varphi_! \circ \varphi^* = 1$, so $\Zs{}(L_{[J_-]}) = \varsigma_! \circ s_! \circ r^*$. This implies that $\cZg{}(L_{[J_-]}) = \sigma_! \circ \cZg{}(L_{[J_+]})$, where $\sigma: [\SL{2}(\CC)/\SL{2}(\CC)] \to [\SL{2}(\CC)/\SL{2}(\CC)]$ is the factorization through the quotient of $\varsigma$. Hence, using the description of $\sigma$ of Section \ref{sec:discs-tubes}, we obtain that the matrix of $\cZg{}(L_{[J_-]})$ in $\cW_0$ is
$$
\cZg{}(L_{[J_-]}) = (q^3-q)\begin{pmatrix}
0 & 0 & 0 & 1 & 0 & 0 & 0 & 0 \\
0 & 0 & 1 & 0 & 0 & 0 & 0 & 0 \\
0 & q^{2} - 1 & q & q - 2 & q^{2} - 2 \, q + 1 & -q
+ 1 & -q + 1 & -2 \, q + 2 \\
q^{2} - 1 & 0 & q - 2 & q & q^{2} - 2 \, q + 1 & -q
+ 1 & -q + 1 & -2 \, q + 2 \\
0 & 0 & q & q & q^{2} - 2 \, q & -q + 1 & -q + 1
& -q + 2 \\
0 & 0 & 0 & 0 & 0 & q & -1 & 0 \\
0 & 0 & 0 & 0 & 0 & -1 & q & 0 \\
0 & 0 & 0 & 0 & q & 0 & 0 & -1
\end{pmatrix}
$$

\subsection{The genus tube}
\label{sec:genus-tube}

In this section, we discuss the case of the holed torus tube $L: (S^1, \star) \to (S^1, \star)$. 
As described there, the standard TQFT induces a morphism $\Zs{}(L): \KVarrel{\SL{2}(\CC)} \to \KVarrel{\SL{2}(\CC)}$.
This morphism has been implicitly studied in \cite{LMN, MM:2016, MM}. However, there, the approach was focused on computing Hodge monodromy representations of representation varieties since the TQFT formalism was not available. Despite of that, as we saw in Section \ref{sec:Hodge-monodromy-covering-spaces}, the monodromy representation has a direct interpretation in terms of virtual classes. 

Using this interpretation, in this section we will show that the computations of \cite{MM} is nothing but the calculation of the geometric TQFT as above, with some peculiarities. For a while, let us consider the general case of an arbitrary (reductive) group $G$. As explained in Section \ref{sec:TQFT-for-repr}, $\Zs{G}(L)=q_! \circ p^*$ with $p: G^4 \to G$ given by $p(a, g, h, b) = a$ and $q: G^4 \to G$ given by $q(a, g, h, b) = ba[g,h]b^{-1}$ (beware of the change of notation for the arguments). Hence, $\Zs{G}(L)^g = (q_!p^*)^g$. In order to compute it explicitly, recall that we have a pullback diagram
\[
\begin{displaystyle}
   \xymatrix
   {	
	G^7 \ar[r]\ar[d] & G^4 \ar[d]^{q} \\
	G^4 \ar[r]_p & G
   }
\end{displaystyle}   
\]
where the upper arrow is $(a, g, h, g', h', b, b') \mapsto (b^{-1}ab,b^{-1}gb,b^{-1}hb,b)$ and the leftmost arrow is $(a, g, h, g', h', b, b')  \mapsto (a[g,h], g', h', b')$. Repeating the procedure, we have that $\Zs{G}(L)^g = (\beta_g)_! \circ (\alpha_g)^*$ where $\alpha_g, \beta_g: G^{3g+1} \to G$ are given by
\begin{align*}
	\alpha_g(a, g_1, h_1, g_2, \ldots, g_g, h_g, b_1, \ldots, b_g) &= b_g\cdots b_1ab_1^{-1} \cdots b_g^{-1}, \\
	\beta_g(a, g_1, h_1, g_2, \ldots, g_g, h_g, b_1, \ldots, b_g) &= a\prod_{i = 1}^g [g_i,h_i].
\end{align*}
In particular, $\Zs{G}(L)^g \circ \Zs{G}(D)(\star) = \RM{G^{3g}}{G}$, where we consider the projection $\gamma_g: G^{3g} \to G$, $\gamma_g(g_1, h_1, \ldots, g_g, h_g, b_1, \ldots, b_g) = \beta_g(1, g_1, h_1, \ldots, g_g, h_g, b_1, \ldots, b_g) = \prod_{i\geq 1} [g_i,h_i]$.

Let us come back to the case $G = \SL{2}(\CC)$. The aim of the paper \cite{MM}, as initiated in \cite{LMN}, is to compute for $g \geq 0$ the Hodge monodromy $\RM{Y_g/\ZZ_2}{\Bt}$ where
$$
	Y_g = \left\{(A_1, B_1, \ldots, A_g, B_g, \lambda) \in \SL{2}(\CC)^{2g} \times \left(\CC^*-\left\{\pm 1\right\}\right)\,\left|\,\,\, \prod_{i=1}^g [A_i,B_i] = D_{\lambda} \right.\right\}.
$$
The action of $\ZZ_2$ on $Y_g$ is $-1\cdot(A_i,B_i, \lambda) = (P_0A_iP_0^{-1}, P_0B_iP_0^{-1}, \lambda^{-1})$, where $P_0 \in \SL{2}(\CC)$ is the matrix that permutes the two vectors of the standard basis of $\CC^2$, and the projection $w_g: Y_g/\ZZ_2 \to \Bt$ is $w_g(A_i,B_i, \lambda) = \lambda + \lambda^{-1}$. In \cite{MM}, it is proven that $\RM{Y_g/\ZZ_2}{\Bt} \in \cW|_{\Bt} \subseteq \KVarrel{[\SL{2}(\CC)/\SL{2}(\CC)]}$ for all $g \geq 0$, where $\cW|_{\Bt}$ denotes the submodule generated by $T_{\Bt}, S_{\pm 2}$ and $S_2 \times S_{-2}$. Moreover, they proved that there exists a module homomorphism $M: \cW|_{\Bt} \to \cW|_{\Bt}$ such that
$
	\RM{Y_{g}/\ZZ_2}{\Bt} = M\left(\RM{Y_{g-1}/\ZZ_2}{\Bt}\right),
$
for all $g \geq 1$.

\begin{rmk}In the paper \cite{MM}, the variety $Y_g$ is denoted by $\overline{X}_4^g$. However, we will not use this notation since it is confusing with ours.
\end{rmk}

\begin{lem}
Let $X_{g} = \left\{(A_i, B_i) \in \SL{2}(\CC)^{2g}\,|\, \prod_{i=1}^g [A_i,B_i]\in D\right\}$ and consider $\omega_g: X_g \to D$ given by $\omega_g(A_i,B_i) =\prod_{i=1}^g [A_i,B_i]$. Then, $\tr^*\RM{Y_g}{\Bt} = \RM{X_g}{D}$.
\begin{proof}
Let us denote by $\overline{Y_g/\ZZ_2} = Y_g/\ZZ_2 \times_{\Bt} D$ the pullback of $w_g: Y_g/\ZZ_2 \to \Bt$ under $\tr: D \to \Bt$, so $\tr^*\RM{Y_g/\ZZ_2}{\Bt} = \RM{\overline{Y_g/\ZZ_2}}{D}$. Now, consider the auxiliar varieties
\begin{align*}
	\overline{Y_g/\ZZ_2}' &= \left\{(\overline{A_i,B_i, \lambda, Q}) \in \left(Y_g \times \SL{2}(\CC)\right)/\ZZ_2 \,\left|\,\,\,\, \prod_{i=1}^g [A_i,B_i] = D_\lambda\right.\right\}, \\
	X_g' &= \left\{(A_i,B_i, \overline{Q}) \in \left(X_g \times \SL{2}(\CC)\right)/\ZZ_2 \,\left|\,\,\,\, Q^{-1}\prod_{i=1}^g [A_i,B_i]Q = D_\lambda \right.\right\},
\end{align*}
where $\ZZ_2$ acts on $\SL{2}(\CC)$ by $-1 \cdot Q = QP_0$. These varieties have morphisms $\overline{Y_g/\ZZ_2}' \to \overline{Y_g/\ZZ_2}$ and $X_g' \to X_g$ given by $(\overline{A_i,B_i, \lambda, Q}) \mapsto ((\overline{A_i,B_i, \lambda}), QD_\lambda Q^{-1})$ and $(A_i,B_i, \overline{Q}) \mapsto (A_i,B_i)$, both with trivial monodromy and fiber $\CC^*$.

Moreover, we have a regular isomorphism $\varphi: \overline{Y_g/\ZZ_2}' \to X_g'$ given by $\varphi(\overline{A_i,B_i, \lambda, Q}) = (QA_iQ^{-1}, QB_iQ^{-1}, Q)$. It fits in a commutative diagram with fibrations as rows
	\[
\begin{displaystyle}
   \xymatrix
   {	
	\CC^* \ar[r] & \overline{Y_g/\ZZ_2}' \ar[r]\ar@{{<}-{>}}[d]^\varphi & \overline{Y_g/\ZZ_2} \\
	\CC^* \ar[r] & X_g' \ar[r] & X_g
   }
\end{displaystyle}   
\]
Therefore, since isomorphisms do not change the monodromy representations, we have
$$
	\tr^*\RM{Y_g/\ZZ_2}{\Bt} = \RM{\overline{Y_g/\ZZ_2}}{D} = \frac{1}{q-1}\RM{\overline{Y_g/\ZZ_2}'}{D} = \frac{1}{q-1}\RM{X_g'}{D} = \RM{X_g}{D},
$$
as we wanted to prove.
\end{proof}
\end{lem}

\begin{cor}\label{prop:formula-M}
For $g \geq 0$, we have
$
	\Zg{}(L) \left(\RM{Y_{g-1}/\ZZ_2}{\Bt}\right)|_{\Bt} = (q^3-q)\,\eta \left(\RM{Y_g/\ZZ_2}{\Bt}\right)$.
\begin{proof}
If we compute the right hand side, by the previous Lemma we have $\eta \left(\RM{Y_g/\ZZ_2}{\Bt}\right) = \tr_! \circ \tr^* \left(\RM{Y_g/\ZZ_2}{\Bt}\right) = \RM{X_g}{\Bt}$ where the projection is $\tr \circ \omega_g: X_g \to \Bt$.

On the other hand, the left hand side can be rewritten as $\Zg{}(L) \left(\RM{Y_{g-1}/\ZZ_2}{\Bt}\right)|_{\Bt} = \tr_! \circ \Zs{}(L) \circ \tr^* \left(\RM{Y_{g-1}/\ZZ_2}{\Bt}\right)|_{\Bt} = \tr_! \circ \Zs{}(L)\left(\RM{X_{g-1}}{D}\right)|_{\Bt}$. Observe that $\gamma_{g}^{-1}(D) = X_{g} \times \SL{2}(\CC)^{g}$ so $\RM{X_g}{D} = \frac{1}{(q^3-q)^g}(\gamma_{g}|_{\gamma_{g}^{-1}(D)})_! T_{X_g} = \frac{1}{(q^3-q)^g}\left(\Zs{}(L)^g \circ Z(D)(\star)\right)|_{D}$. Using that, we finally obtain that
\begin{align*}
	\tr_! \circ \Zs{}(L)\left(\RM{X_g}{D}\right) &= \frac{1}{(q^3-q)^{g-1}}\,\tr_! \circ \Zs{}(L)\left(\Zs{}(L)^{g-1} \circ Z(D)(\star)\right)|_{D} \\
	&= \frac{1}{(q^3-q)^{g-1}}\,\tr_! \left(\Zs{}(L)^g \circ Z(D)()\right)|_{D} = (q^3-q)\,\tr_! \left(\RM{X_g}{D}\right) \\
	&= (q^3-q)\RM{X_g}{\Bt}.
\end{align*}
This proves the desired equality.
\end{proof}
\end{cor}

\begin{cor}
The morphism $M: \cW|_{\Bt} \to \cW|_{\Bt}$ satisfies that $(q^3-q)M = \eta^{-1} \circ \Zg{}(L)|_{\Bt}$.
\begin{proof}
By the results of \cite{MM}, the set $\left\{\RM{Y_g/\ZZ_2}{\Bt}\right\}_{g \geq 0}$ generate the same submodule as the elements $T_{\Bt}, S_2, S_{-2}$ and $S_2 \times S_{-2}$. Therefore, since $\RM{Y_g/\ZZ_2}{\Bt} = M\left(\RM{Y_{g-1}/\ZZ_2}{\Bt}\right)$, from Proposition \ref{prop:formula-M} we obtain that $(q^3-q) \eta M = \Zg{}(L)$ on the submodule generated by $T_{\Bt}, S_2, S_{-2}$ and $S_2 \times S_{-2}$.
\end{proof}
\end{cor}

Indeed, in \cite[Section 9]{MM}, a larger homomorphism $M: \cW \to \cW$ is defined for which the previous one is just its restriction to $\cW|_{\Bt}$. 
Analogous (and simpler) calculations can be done as in Corollary \ref{prop:formula-M} in order to show that $(q^3-q) \eta \circ M(T_{\pm 1}) = \Zg{}(L)(T_{\pm 1})$ and $(q^3-q)\eta \circ M(T_{\pm}) = \Zg{}(L)(T_{\pm})$. Hence, $(q^3-q)M = \eta^{-1} \circ \Zg{}(L)$ on the whole $\cW$. Using the explicit expression of $M$ from \cite{MM}, we have that, in the set of generators $\cS$, the matrix of $\Zg{}(L)$ is
$$
\scriptsize{
\fontsize{0.8}{10} \selectfont (q^3-q)^2      	
\begin{pmatrix}
q + 4 & 1 & q^{2} - 2 q - 3 & q^{2} + 3  q & \substack{q^{3}
- 2  q^{2}\\ - 3  q - 2} & -q^{2} - 4  q - 1 & 2  q^{2} - 7
 q - 1 & -5  q - 1 \\[0.15cm]
1 & q + 4 & q^{2} + 3  q & q^{2} - 2  q - 3 & \substack{q^{3}
- 2  q^{2}\\ - 3  q - 2} & 2  q^{2} - 7  q - 1 & -q^{2} - 4
 q - 1 & -5  q - 1 \\[0.15cm]
q^{2} - 2  q - 3 & q^{2} + 3  q & \substack{q^{4} + q^{3} \\+ 3  q + 3}
& q^{4} - 3  q^{2} - 6  q & \substack{q^{5} - 2  q^{4} - 3  q^{3}\\
+ q^{2} + 3  q} & \substack{-q^{4} + 2  q^{3}\\ - 4  q^{2} + 3  q} &
\substack{-q^{4} - q^{3}\\ - 4  q^{2} + 6  q} & -2  q^{3} - q^{2} + 3  q
\\[0.15cm]
q^{2} + 3  q & q^{2} - 2  q - 3 & q^{4} - 3  q^{2} - 6 
q & \substack{q^{4} + q^{3} \\+ 3  q + 3} & \substack{q^{5} - 2  q^{4} - 3  q^{3}
\\+ q^{2} + 3  q} & \substack{-q^{4} - q^{3} \\- 4  q^{2} + 6  q} & \substack{-q^{4}
+ 2  q^{3}\\ - 4  q^{2} + 3}  q & -2  q^{3} - q^{2} + 3  q \\[0.15cm]
q^{2} + 1 & q^{2} + 1 & q^{4} - 2  q^{2} & q^{4} - 2 
q^{2} & \substack{q^{5} - 2  q^{4} - q^{3}\\ + 2  q^{2} - 2} & \substack{-q^{4} - q^{3} \\+ q^{2} - q - 1} & \substack{-q^{4} - q^{3} \\+ q^{2} - q - 1} & \substack{-2 
q^{3} \\+ q^{2} - 2  q - 1} \\[0.15cm]
0 & 3  q & 3  q^{2} & -3  q & -3  q^{2} & 4
 q^{3} - 6  q^{2} & -4  q^{2} & -3  q^{2} \\[0.15cm]
3  q & 0 & -3  q & 3  q^{2} & -3  q^{2} & -4
 q^{2} & 4  q^{3} - 6  q^{2} & -3  q^{2} \\[0.15cm]
q & q & q^{3} & q^{3} & \substack{q^{4} - 2  q^{3}\\ - q^{2} - 2
 q} & -q^{3} - q^{2} - q & -q^{3} - q^{2} - q & q^{3} - 2
 q^{2} - q
\end{pmatrix}}
$$
Analogous calculation can be done to obtain $\cZg{}(L) = (q^3-q)\eta M\eta^{-1}$.

\begin{rmk}
The previous result can be restated as that $M: \cW_0 \to \cW_0$ is (up to a constant) the left $\tr$-reduction of the standard TQFT, as explained in Remark \ref{rmk:left-reduction}. However, we have chosen to consider right $\tr$-reductions in Section \ref{sec:geometric-reduced-TQFT} since they are more natural from the geometric point of view.
\end{rmk}

Let $\Sigma_g$ be the closed surface of genus $g \geq 0$. For the parabolic data $\Lambda = \left\{\left\{- \Id\right\},[J_+], [J_-]\right\}$ we consider a parabolic structure on $\Sigma_g$, $Q = \left\{(p_1, [C_1]), \ldots, (p_s, [C_s])\right\}$, where $C_i = J_+, J_-$ or $-\Id$ and $p_1, \ldots, p_s \in \Sigma_g$ with $s > 0$. Let us denote by $r_+$ be the number of $[J_+]$ in $Q$, $r_-$ the number of $[J_-]$ and $t$ the number of $\left\{-\Id\right\}$ (so that $r_+ + r_- + t = s$). Set $r = r_+ + r_-$ and $\sigma = (-1)^{r_- + t}$.

\begin{thm}\label{thm:Hodge-repr-sl2}
The virtual class of $\Rep{\SL{2}(\CC)}(\Sigma_g, Q)$ is
\begin{itemize}
	\item If $\sigma = 1$, then
\begin{align*}
	\coh{\Rep{\SL{2}(\CC)}(\Sigma_g, Q)} =&\, {\left(q^2 - 1\right)}^{2g + r - 1} q^{2g - 1} +
\frac{1}{2} \, {\left(q -
1\right)}^{2g + r - 1}q^{2g -
1}(q+1){\left({2^{2g} + q - 3}\right)} 
\\ &+ \frac{\left(-1\right)^{r}}{2} \,
{\left(q + 1\right)}^{2g + r - 1} q^{2g - 1} (q-1){\left({2^{2g} +q -1}\right)}.
\end{align*}
	\item If $\sigma = -1$, then
\begin{align*}
	\coh{\Rep{\SL{2}(\CC)}(\Sigma_g, Q)} = &\, {\left(q - 1\right)}^{2g + r - 1} (q+1)q^{2g - 1}{{\left( {\left(q + 1\right)}^{2 \,
g + r-2}+2^{2g-1}-1\right)} } \\
&+ \left(-1\right)^{r + 1}2^{2g - 1}  {\left(q + 1\right)}^{2g + r
- 1} {\left(q - 1\right)} q^{2g - 1}.
\end{align*}
\end{itemize}
\begin{proof}
By Theorem \ref{thm:almost-tqft-parabolic}, we have that
\small
$$
	\coh{\Rep{\SL{2}(\CC)}(\Sigma_g, Q)} = \frac{1}{(q^3-q)^{g+s}}\cZg{}(D^\dag) \circ \cZg{}(L_{[C_s]}) \circ \ldots \circ \cZg{}(L_{[C_1]}) \circ \cZg{}(L)^g \circ \cZg{}(D)(\star).
$$
\normalsize
Moreover observe that, as $\cZg{}$ is an almost-TQFT and the strict tubes commute, all these linear morphisms commute. Hence, we can group the $-\Id$ and $J_-$ tubes together and, using that $\cZg{}(L_{-\Id}) \circ \cZg{}(L_{-\Id}) = (q^3-q)^2 1_{\KVarrel{\SL{2}(\CC)}}$ and $\cZg{}(L_{-\Id}) \circ \cZg{}(L_{[J_-]}) = (q^3-q)\cZg{}(L_{[J_+]})$, we finally have:
\begin{itemize}
	\item If $\sigma = 1$, then all the $-\Id$ tubes cancel so we have
$$
	\coh{\Rep{\SL{2}(\CC)}(\Sigma_g, Q)} = \frac{1}{(q^3-q)^{g+r}}\cZg{}(D^\dag) \circ \cZg{}(L_{[J_+]})^r \circ \cZg{}(L)^g \circ \cZg{}(D)(\star).
$$
Now, observe that, as $\cZg{}(L)$ and $\cZg{}(L_{[J_+]})$ commute, they can be simultaneously diagonalized. Hence, there exists $P, A, B: \cW_0 \to \cW_0$, with $A$ and $B$ diagonal matrices, such that $PAP^{-1} = \cZg{}(L)$ and $PBP^{-1}=\cZg{}(L_{[J_+]})$. Therefore, we have that the virtual class is given by
$$
	\coh{\Rep{\SL{2}(\CC)}(\Sigma_g, Q)} = \frac{1}{(q^3-q)^{g+r}}\cZg{}(D^\dag) PB^rA^gP \cZg{}(D)(\star).
$$
From the matrices of $\cZg{}(L)$ and $\cZg{}(L_{[J_+]})$ for the set of generators $\cS$, the matrices $P,A$ and $B$ can be explicitly calculated with a symbolic computation software. Using these matrices, and recalling that $\cZg{}(D)$ and $\cZg{}(D^\dag)$ are the inclusion and projection onto $T_1$ respectively, the calculation follows.
	\item If $\sigma = -1$, then all the $-\Id$ tubes cancel except one, so we have
\small
$$
	\coh{\Rep{\SL{2}(\CC)}(\Sigma_g, Q)} = \frac{1}{(q^3-q)^{g+r+1}}\cZg{}(D^\dag) \circ \cZg{}(L_{-\Id}) \circ \cZg{}(L_{[J_+]})^r \circ \cZg{}(T)^g \circ \cZg{}(D)(\star).
$$
\normalsize
Again, decomposing $PAP^{-1} = \cZg{}(T)$ and $PBP^{-1}=\cZg{}(L_{[J_+]})$, we obtain that the virtual class is
$$
	\coh{\Rep{\SL{2}(\CC)}(\Sigma_g, Q)} = \frac{1}{(q^3-q)^{g+r+1}}\cZg{}(D^\dag) \cZg{}(L_{-\Id}) PB^rA^gP \cZg{}(D)(\star).
$$
Using a symbolic computation software, the result follows.
\end{itemize}
\end{proof}
\end{thm}

\begin{rmk}
The homomorphism $\cZg{}(L_{[J_+]}): \cW_0 \to \cW_0$ has non-trivial kernel. In this way, if $A$ is the diagonal matrix associated to $\cZg{}(L_{[J_+]})$, as in the proof of Theorem \ref{thm:Hodge-repr-sl2}, we have that $A^0 \neq \Id$ because of the vanishing eigenvalues. For this reason, we cannot directly set $r=0$ in the previous formula in order to obtain the virtual class of non-parabolic representation varieties. However, it may be traced back the extra term lost when we set $r=0$ and to obtain the formula for the non-parabolic case
\begin{align*}
	\coh{\Rep{\SL{2}(\CC)}(\Sigma_g)} &= \left.\coh{\Rep{\SL{2}(\CC)}(\Sigma_g, Q)}\right|_{r=0,\sigma=1} + q(q^2-1)^{2g-1} \\
	& = {\left(q^2 - 1\right)}^{2g - 1} q^{2g - 1} +
\frac{1}{2} \, {\left(q -
1\right)}^{2g - 1}q^{2g -
1}(q+1){\left({2^{2g} + q - 3}\right)} 
\\ &\;\;\;\;\,+ \frac{1}{2} \,
{\left(q + 1\right)}^{2g + r - 1} q^{2g - 1} (q-1){\left({2^{2g} +q -1}\right)} + q(q^2-1)^{2g-1}.
\end{align*}
Here $\left.\coh{\Rep{\SL{2}(\CC)}(\Sigma_g, Q)}\right|_{r=0,\sigma=1}$ denotes the virtual class obtained when setting $r=0$ in the formula for $\coh{\Rep{\SL{2}(\CC)}(\Sigma_g, Q)}$ in the case $\sigma = 1$ of Theorem \ref{thm:Hodge-repr-sl2}.

Analogous considerations can be done in the case of only one marked point with conjugacy class $[-\Id]$. In that case, if we denote $Q_0 = \left\{(p, [-\Id])\right\}$ with $p \in \Sigma_g$ we have 
\begin{align*}
	\coh{\Rep{\SL{2}(\CC)}(\Sigma_g, Q_0)} &= \left.\coh{\Rep{\SL{2}(\CC)}(\Sigma_g, Q)}\right|_{r=0,\sigma=-1} + q(q^2-1)^{2g-1} \\ &=\, {\left(q - 1\right)}^{2g + r - 1} (q+1)q^{2g - 1}{{\left( {\left(q + 1\right)}^{2 \,
g + r-2}+2^{2g-1}-1\right)} } \\
&\;\;\;\;\,+ \left(-1\right)^{r + 1}2^{2g - 1}  {\left(q + 1\right)}^{2g + r
- 1} {\left(q - 1\right)} q^{2g - 1} + q(q^2-1)^{2g-1}.
\end{align*}
Now, $\left.\coh{\Rep{\SL{2}(\CC)}(\Sigma_g, Q)}\right|_{r=0,\sigma=-1}$ denotes the virtual class of the case $\sigma = -1$ of Theorem \ref{thm:Hodge-repr-sl2} after setting $r=0$.
\end{rmk}

\begin{rmk}
When interpreting $q$ as a variable, the previous polynomials are also the Deligne-Hodge polynomials \cite{Hausel-Rodriguez-Villegas:2008} of the parabolic representations varieties with $q=uv$. This is due to the fact that $\DelHod{q} = \DelHod{\coh{\CC}} = uv$.

Under this point of view, the results of Theorem \ref{thm:Hodge-repr-sl2} generalize the previous work on Deligne-Hodge polynomials of parabolic representation varieties to the case of an arbitrary number of marked points. For genus $g=1,2$ and one marked point in $[J_\pm]$ (i.e.\ $r_+ = 1$ or $r_-=1$), this theorem agrees with the calculations of \cite[Sections 4.3, 4.4, 11 and 12]{LMN}. For arbitrary genus and zero or one marked points, this result agrees with \cite[Proposition 11]{MM}. Finally, for the case $g=1$ and two marked points in the classes $[J_\pm]$, this theorem agrees with the computations of \cite[Section 3]{LM}. 
\end{rmk}

%%%%%%%%%%%%%%%%%%%%%%%%%%%%%%%%%%%%%%%%%%%%%%%%%%%%%%%%%%%%%%%%
%%%%%%%%%%%%%%%%%%%%%%% BIBLIOGRAPHY %%%%%%%%%%%%%%%%%%%%%%%%%%%
%%%%%%%%%%%%%%%%%%%%%%%%%%%%%%%%%%%%%%%%%%%%%%%%%%%%%%%%%%%%%%%%

%\nocite{*}
\bibliography{bibliography.bib}{}
\bibliographystyle{abbrv}

\end{document}